\documentclass[aap,seceqn,MSNbibl,citesort,dvips]{arximspdf}
\usepackage{graphicx}

\doi{10.1214/10-AAP742}
\volume{21}
\issue{5}
\pubyear{2011}
\firstpage{1860}
\lastpage{1910}

\makeatletter

\newtheorem{theorem}{Theorem}
\newtheorem{lemma}{Lemma}
\newtheorem{remark}{Remark}

\newcommand{\rmtr}{\mathrm{tr}} % laipsnyje tr

\newcommand{\bbA}{\mathbf{A}}
\newcommand{\cbbA}{{\check\bbA}}
\newcommand{\tbbA}{{\tilde\bbA}}
\newcommand{\bbB}{\mathbf{B}}
\newcommand{\bbC}{\mathbf{C}}
\newcommand{\bbD}{\mathbf{D}}
\newcommand{\bbe}{\mathbf{e}}
\newcommand{\bbI}{\mathbf{I}}
\newcommand{\bbr}{\mathbf{r}}
\newcommand{\bbS}{\mathbf{S}}
\newcommand{\cbbS}{{\check\bbS}}
\newcommand{\bbt}{\mathbf{t}}
\newcommand{\bbU}{\mathbf{U}}
\newcommand{\bbX}{\mathbf{X}}
\newcommand{\tbbX}{\widetilde{\mathbf X}}
\newcommand{\hbbX}{\widehat{\mathbf X}}
\newcommand{\bbY}{\mathbf{Y}}

\newcommand{\uij}{\,ij}

\makeatother

\begin{document}
\begin{frontmatter}

\title{Central limit theorem for Hotelling's $T^2$ statistic under
large dimension}
\runtitle{Hotelling's $T^2$ statistic}

\begin{aug}
\author[A]{\fnms{G. M.} \snm{Pan}\thanksref{t1}\ead[label=e1]{gmpan@ntu.edu.sg}} and
\author[B]{\fnms{W.} \snm{Zhou}\corref{}\thanksref{t2}\ead[label=e2]{stazw@nus.edu.sg}}
\runauthor{G. M. Pan and W. Zhou}
\affiliation{Nanyang Technological University and National University
of Singapore}
\dedicated{Dedicated to Z. D. Bai on the occasion of his 65th birthday}
\address[A]{Division of Mathematical Sciences\\
School of Physical\\
\quad and Mathematical Sciences\\
Nanyang Technological University\\
Singapore 637371\\
\printead{e1}}
\address[B]{Department of Statistics\\
\quad and Applied Probability\\
National University of Singapore\\
Singapore 117546\\
\printead{e2}}
\end{aug}

\thankstext{t1}{Supported in part by Grant M58110052 at the Nanyang
Technological
University.}

\thankstext{t2}{Supported in part by Grant
R-155-000-106-112 at the National University of Singapore.}

% HISTORY:
\received{\smonth{11} \syear{2009}}
\revised{\smonth{9} \syear{2010}}

% ABSTRACT
%
\begin{abstract}
In this paper we prove the central limit theorem for Hotelling's
$T^2$ statistic when the dimension of the random vectors is
proportional to the sample size.
\end{abstract}

% KEYWORDS
%

\end{frontmatter}

%s1 ###
\section{Introduction and main results}

Since the famous Mar{\v{c}}enko and Pastur
law was found in \cite{MP}, the theory of large sample covariance
matrices has been further developed. Among others, we mention Jonsson
\cite{Jon}, Yin \cite{y1}, Silverstein~\cite{s1}, Watcher \cite
{watch}, Yin, Bai and Krishanaiah \cite{y2}.
Lately, Johnstone \cite{john} discovered the law of the largest
eigenvalue of the Wishart matrix, Bai and Silverstein \cite{b2}
established the
central limit theorems (CLT) of linear spectral statistics and Bai,
Miao and Pan \cite{b1} derived CLT for functionals of the
eigenvalues and eigenvectors. We also refer to \cite{johan,ss,Di}
for CLT on linear statistics of eigenvalues of other
classes of random matrices.

The sample covariance matrix is
defined by
\[
\bolds{\mathcal{S}}=\frac{1}{n}\sum_{j=1}^n({\mathbf s}_j-\bar
{\mathbf s})({\mathbf s}_j-\bar{\mathbf s})^T,
\]
where $\bar{\mathbf s}=n^{-1}\sum_{j=1}^n{\mathbf s}_j$ and ${\mathbf
s}_j=(X_{1j},\ldots,X_{pj})^T$. Here \{$X_{ij}$\}, $i,j=\cdots,$ is a
double array of independent and identically distributed (i.i.d.)
real r.v.'s with $EX_{11}=0$ and $EX_{11}^2=1$.
However, in
the large random matrices theory (RMT) the commonly used sample
covariance matrix is
\[
\bbS=\frac{1}{n}\sum_{j=1}^n{\mathbf s}_j{\mathbf s}_j^T=\frac
{1}{n}\bbX_n\bbX_n^T,
\]
where $\bbX_n=({\mathbf s}_1,\ldots,{\mathbf s}_n)$.

Note that $ \bolds{\mathcal{S}}=\bbS-\bar{\mathbf s}\bar{\mathbf s}^T $ and
thus, by
the rank inequality, there is no difference when one is only
concerned with the limiting empirical spectral distribution (ESD)
of the eigenvalues in large random matrices. Therefore, the
limiting ESD of $\bolds{\mathcal{S}}$ is Mar{\v{c}}enko and Pastur's law
$F_c(x)$ (see \cite{Jon} and \cite{MP}) when
$\lim\frac{p}{n}=c>0$ which has a density function
\[
p_c(x)=
\cases{(2\pi cx)^{-1}\sqrt{(b-x)(x-a)}, &\quad$a\leq x\leq b$,\cr
0, &\quad otherwise,}
\]
and has point mass $1-c^{-1}$ at the origin if $c>1$, where $a=(1-\sqrt
c)^2$ and $b=(1+\sqrt c)^2$.
The Stieljes transform $m(z)$ of $F_c(x)$
satisfies the equation (see \cite{s3})
%
%e1.1 ###
%
\begin{equation}
\label{a67} m(z)=\frac{1}{1-c-czm(z)-z},
\end{equation}
where the Stieljes transform for any function $G(x)$ is defined by
\[
m_G(z)=\int\frac{1}{\lambda-z}\,dG(\lambda), \qquad z\in
{\Bbb{C}}^+\equiv\{z\in{\Bbb C}, v=\Im z>0\}.
\]
Observe that the spectra of $n^{-1}\bbX_n\bbX_n^T$ and
$n^{-1}\bbX_n^T\bbX_n$ are identical except for zero
eigenvalues. This leads to the equality
%
%e1.2 ###
%
\begin{equation}\label{d22}
\underline{m}_n^{\bbS}(z)=-\frac{1-p/n}{z}+\frac{p}{n}m_n^{\bbS}(z)
\end{equation}
and therefore,
%
%e1.3 ###
%
\begin{equation}\label{e1}
z=-\frac{1}{\underline{m}(z)}+\frac{c}{1+\underline{m}(z)},
\end{equation}
where\vspace*{1pt} $m_n^{\bbS}(z)$ and $\underline{m}_n^{\bbS}(z)$ denote,
respectively, the Stieljes transform of the ESD of $n^{-1}\bbX_n\bbX
_n^T$ and $n^{-1}\bbX_n^T\bbX_n$ and, correspondingly, $\underline{m}(z)$
is the limit of $\underline{m}_n^{\bbS}(z)$.

Sample covariance matrices are also of essential importance in multivariate
statistical analysis because many test statistics involve their
eigenvalues and/or eigenvectors. The typical example is $T^2$ statistic
which was
proposed by Hotelling \cite{h}. We refer to \cite{a1} and
\cite{leh} for various uses of the $T^2$ statistic.

The $T^2$ statistic, which is the origin of multivariate linear
hypothesis tests and the associated confidence sets, is defined by
%
%e1.4 ###
%
\begin{equation}
T^2=n(\bar{\mathbf s}-\bolds{\mu}_0)^T\bolds{\mathcal{S}}^{-1}(\bar
{\mathbf
s}-\bolds{\mu}_0),
\end{equation}
whose distribution is invariant if each ${\mathbf s}_j$ is replaced by
$\bolds\Sigma^{1/2}{\mathbf s}_j$ with $\bolds\Sigma$ any
nonsingular $p$ by $p$ matrix when $\bolds{\mu}_0=0$. If
$\{{\mathbf s}_1,\ldots,{\mathbf s}_n\}$ is a sample from the $p$-dimensional
population $N(\bolds{\mu},\bolds\Sigma)$, then
$[T^2/(n-1)][(n-p)/p]$ follows a noncentral $F$
distribution and moreover, the $F$ distribution is central if
$\bolds{\mu}=\bolds{\mu}_0$. When $p$ is fixed, the
limiting distribution of $T^2$ for
$\bolds{\mu}=\bolds{\mu}_0$ is the $\chi^2$-distribution
even if the parent distribution is not normal.

In the recent three or four decades in many research areas,
including signal processing, network security, image processing,
genetics, stock
marketing and other economic problems, people are interested in the case
where $p$ is quite large or proportional to the sample size.
Thus, it will be desirable if one
can obtain the asymptotic distribution of the famous Hotelling
$T^2$ statistic when the dimension of the random vectors is
proportional to the sample size. It is the aim of this work. In
addition, we would like to point out that some discussions about
the two-sample $T^2$ statistic under the assumption that the
underlying r.v.'s are normal were presented in \cite{b5}.

The main results are presented in the following theorems.
\begin{theorem}
\label{th1} Suppose that:

\begin{longlist}[(2)]
\item[(1)] for each $n$ $X_{ij}=X_{ij}^n, i,j=1,2,\ldots,$ are i.i.d. real
r.v.'s with $EX_{11}=\mu, EX_{11}^2=1$
and $EX_{11}^4<\infty$.

\item[(2)]
$p\leq n, c_n=p/n\to c\in(0,1)$ as $n\to\infty$.

Then, when
$\bolds{\mu}_0=(\mu,\ldots,\mu)^T$,
\[
\frac{\sqrt{n}}{\sqrt{2c_n(1-c_n)^{-3}}}\biggl(\frac
{T^2}{n}-c_n(1-c_n)^{-1}\biggr)\stackrel{D}\longrightarrow N(0,1),
\]
where
$F_{c_n}(x)$ denotes $F_{c}(x)$ by substituting $c_n$ for $c$.
\end{longlist}
\end{theorem}
\begin{remark}
When $X_{ij}\sim N(0,1)$, it is well known that $(n-p)T^2/(np)$
follows $F$ distribution with degrees of freedom $p$ and $n-p$,
respectively. As $n\to\infty$ and $p/n \to c$, it follows from strong
law of large numbers and CLT that
\[
\frac{(n-p)T^2/(np)-1}{\sqrt{2/p+2/(n-p)}} \longrightarrow N(0,1).
\]
This is consistent with Theorem \ref{th1}.
\end{remark}
\begin{remark}
Since $\int x^{-1}\,dF_{c}(x)=(1-c)^{-1}$ and $\int
x^{-2}\,dF_{c}(x)=(1-c)^{-3}$ which are derived through differentiating
the following identity [the Stieljes transform $m(z)$ of $F_c(x)$],
\[
\int(x-z)^{-1}p_c(x)\,dx=\frac{-(z+c-1)+(z+c-1)\sqrt{1-4zc(z+c-1)^{-2}}}{2cz},
\]
we actually prove that
\[
\frac{\sqrt{n}}{\sqrt{2c_n\int x^{-2}
\,dF_{c_n}(x)}}\biggl(\frac{T^2}{n}-c_n\int
\frac{dF_{c_n}(x)}{x}\biggr)\stackrel{D}\longrightarrow N(0,1).
\]
\end{remark}

One typical application of Theorem \ref{th1} lies in making inference
on the large-dimensional mean vector of the multivariate model
\[
\mathbf{Z}_j=\Gamma{\mathbf s}_j+\bolds{\mu},\qquad E{\mathbf
s}_j=0,\qquad
j=1,\ldots,n,
\]
where $\Gamma$ is an $m$ by $p$ matrix, $m\leq p$. This model means
that each $\mathbf{Z}_j$
is a linear transformation of some $p$-variate random vector ${\mathbf
s}_j$. It can generate a rich
collection of $\mathbf{Z}_j$ from ${\mathbf s}_j$ with the given
covariance matrix $\bolds\Sigma=\Gamma\Gamma^T$. Most
important, it includes the multivariate normal model.

We will prove Theorem \ref{th1} by establishing Theorem \ref{th2}
which presents asymptotic distributions of random quadratic forms
involving sample means and sample covariance matrices.

For any analytic function $f(\cdot)$, define
\[
f(\bolds{\mathcal{S}})=\bbU^T\operatorname{diag}(f(\lambda_1),\ldots
,f(\lambda_p))\bbU,
\]
where $\bbU^T\operatorname{diag}(\lambda_1,\ldots,\lambda_p)\bbU$ denotes the
spectral decomposition of the matrix $\bolds{\mathcal{S}}$.
\begin{theorem}
\label{th2}
In addition to the assumption $(1)$ of Theorem \ref{th1}, suppose
that $c_n=p/n\to c>0$, $EX_{11}=0$, $g(x)$ is a function with a
continuous first derivative in a neighborhood of $c$ and $f(x)$ is
analytic on an open region containing the interval
%
%e1.5 ###
%
\begin{equation}\label{f4}
\bigl[ I_{(0,1)}(c)\bigl(1-\sqrt{c}\bigr)^2, \bigl(1+\sqrt{c}\bigr)^2\bigr].
\end{equation}
Then,
\[
\biggl(\sqrt{n}\biggl[\frac{\bar{\mathbf s}^Tf(\bolds{\mathcal{S}})\bar
{\mathbf s}}{\|\bar{\mathbf s}\|^2}-\int f(x)\,dF_{c_n}(x)\biggr],\sqrt
{n}\bigl(g(\bar
{\mathbf s}^T\bar{\mathbf s})-g(c_n)\bigr)\biggr)
\stackrel{D}\longrightarrow(X,Y),
\]
where $Y\sim N(0,2c(g'(c))^2)$, which is independent of $X$, a Gaussian
r.v. with $EX=0$ and
%
%e1.6 ###
%
\begin{equation}\label{e3}
\operatorname{Var}(X)=\frac{2}{c}\biggl(\int f^2(x)\,dF_c(x)-\biggl(\int f(x)\,dF_c(x)\biggr)^2\biggr).
\end{equation}
\end{theorem}
\begin{remark}\label{re2}
Let $\mathbf{x}_n=(x_{n1},\ldots
,x_{np})^T\in\mathbb{R}^p$, $\|\mathbf{x}_n\|=1$ where \mbox{$\|\cdot\|$}
denotes the Euclidean norm. Note that, when $\max
_{i}x_{ni}\rightarrow0$
(see \cite{p2}, (1.16), or \cite{s2}),
%
%e1.7 ###
%
\begin{equation}\label{c18}
\sqrt{n}\biggl[\mathbf{x}_n^Tf(\mathbf{S})\mathbf{x}_n-\int
f(x)\,dF_{c_n}(x)\biggr]\stackrel{D}\longrightarrow X.
\end{equation}
This suggests that $\bar{\mathbf s}/\|\bar{\mathbf s}\|$ can be viewed
as a
fixed unit vector $\mathbf{x}_n$ when dealing with
$\bar{\mathbf s}^Tf(\bolds{\mathcal{S}})\bar{\mathbf s}/\|\bar{\mathbf
s}\|^2$
even if $\bar{\mathbf s}$ is not independent of
$\bolds{\mathcal{S}}$.
\end{remark}

\begin{figure}

\includegraphics{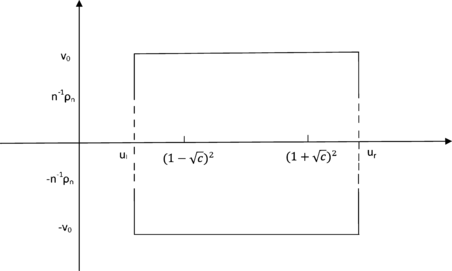}

\caption{Contour $\mathcal{C}$ when $c<1$.}
\label{figure1}
\end{figure}

\begin{figure}[b]

\includegraphics{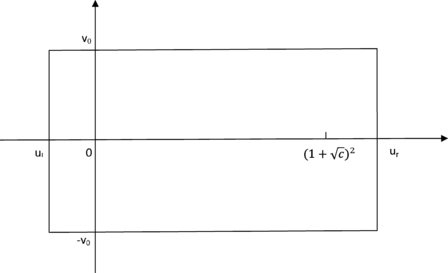}

\caption{Contour $\mathcal{C}$ when $c\geq1$.}
\label{figure2}
\end{figure}

Theorem \ref{th2} relies on Lemma \ref{theo3} below which deals
with the asymptotic joint distribution of
\[
X_n(z)=\sqrt{n}\biggl[\frac{\bar{\mathbf s}^T(\bolds{\mathcal{S}}-z\bbI
)^{-1}\bar{\mathbf s}}{\|\bar{\mathbf s}\|^2}-m_n(z)\biggr],\qquad
Y_n=\sqrt{n}\bigl(g(\bar{\mathbf s}^T\bar{\mathbf s})-g(c_n)\bigr),
\]
where
$m_n(z)=\int(x-z)^{-1}\,dF_{c_n}(x)$. The stochastic process $X_n(z)$
is defined on a contour $\mathcal{C}$, given below. Let $v_0>0$
be arbitrary and set $\mathcal{C}_u=\{u+iv_0, u\in[u_l,u_r]\},$
where $u_l$ is any negative number if the left endpoint of
(\ref{f4}) is zero, otherwise $u_l$ is any positive number smaller
than the left endpoint of (\ref{f4}) and $u_r$ any number larger
than the right endpoint of (\ref{f4}). Then define
\[
\mathcal{C}^{+}=\{u_l+iv\dvtx v\in[0,v_0]\}\cup\mathcal{C}_u\cup
\{u_r+iv\dvtx v\in[0,v_0]\}
\]
and let $\mathcal{C}^{-}$ be the symmetric part of
$\mathcal{C}^{+}$ about the real axis. Then set
\mbox{$\mathcal{C}=\mathcal{C}^{+}\cup\mathcal{C}^{-}$}. See Figures \ref
{figure1} and
\ref{figure2} for a picture of the contour
$\mathcal{C}$ when $c<1$ and $c\geq1$, respectively.

%
%f1 ###
%

%
%f2 ###
%

Let $\bbA^{-1}(z)=(\bbS-zI)^{-1}$. Since it is difficult to control
the spectral norm of $(\bolds{\mathcal{S}}-z\bbI)^{-1}$ or
$\bbA^{-1}(z)$ on the whole contour $\mathcal{C}$, especially\vadjust{\goodbreak} for
$v=0$, we further define $\hat{X}_n(z)$, a truncated version of
$X_n(z)$, as in \cite{b2}. Select a sequence of positive numbers
$\rho_n$ satisfying for some $\beta\in(0,1)$,
%
%e1.8 ###
%
\begin{equation}\label{g21}
\rho_n\downarrow0, \qquad \rho_n\geq n^{-\beta}.
\end{equation}
Let
\[
\mathcal{C}_l=\cases{\{u_l+iv\dvtx v\in[n^{-1}\rho_n,v_0]\}, &\quad
if $u_l>0$,\cr
\{u_l+iv\dvtx v\in[0,v_0]\}, &\quad if $u_l<0$,}
\]
and
\[
\mathcal{C}_r=\{u_r+iv\dvtx v\in[n^{-1}\rho_n,v_0]\}.
\]
Write
$\mathcal{C}_n^+=\mathcal{C}_l\cup\mathcal{C}_u\cup\mathcal{C}_r$.
We can now define the truncated process for $z=u+iv \in\mathcal{C}$ by
%
%e1.9 ###
%
\begin{equation}\label{g20}
\hat{X}_n(z)=\cases{X_n(z),\qquad \mbox{if $z\in
\mathcal{C}_n^+\cup\mathcal{C}_n^-$},\cr
\displaystyle \frac{nv+\rho_n}{2\rho_n}X_n(z_{r1})+\frac{\rho_n-nv}{2\rho
_n}X_n(z_{r2}),\cr
\hspace*{1.5pt}\phantom{X_n(z),}\qquad
\mbox{if $u=u_r,v\in[-n^{-1}\rho_n,n^{-1}\rho_n]$},\cr
\displaystyle \frac{nv+\rho_n}{2\rho_n}X_n(z_{l1})+\frac{\rho_n-nv}{2\rho
_n}X_n(z_{l2}),\cr
\hspace*{1.5pt}\phantom{X_n(z),}\qquad
\mbox{if $u=u_l>0,v\in[-n^{-1}\rho_n,n^{-1}\rho_n]$},}
\end{equation}
where $z_{r1}=u_r+in^{-1}\rho_n, z_{r2}=u_r-in^{-1}\rho_n,
z_{l1}=u_l+in^{-1}\rho_n, z_{l2}=u_l-in^{-1}\rho_n$ and
$\mathcal{C}_n^-$ denotes the symmetric part of $\mathcal{C}_n^+$
about the real axis. A picture of $\mathcal{C}_n^+\cup\mathcal
{C}_n^- $ is the rectangle in Figure \ref{figure1} with the dash line
removed. The
advantage of $\hat{X}_n(z)$ over $X_n(z)$
is that the spectral norm of $\bbA^{-1}(z)$ involved in
$\hat{X}_n(z)$ may be well controlled on the contour $\mathcal{C}$.
Indeed, loosely speaking, all eigenvalues of $\bbS$ are located
inside the interval (\ref{f4}) with a high probability. Therefore,
the spectral norm of $\bbA^{-1}(z)$ corresponding to this case is
bounded on $\mathcal{C}$. If some eigenvalues run outside of the
interval (\ref{f4}) then, at least, we will still have an upper
bound $n\rho_n^{-1}$ for the spectral norm of $\bbA^{-1}(z)$ on
$\mathcal{C}$. But, the probability that some eigenvalues run
outside of the interval (\ref{f4}) is very small, which can offset
$n\rho_n^{-1}$ and even more. This is crucial to establish tightness
of $\hat{X}_n(z)$ on the contour $\mathcal{C}$. On the other hand,
such a truncation does not change the weak limit given in Theorem
\ref{th2} because the truncation has been made only at the intervals
of the length $2\rho_n/n$.

Note that $\hat{X}_n(z)$ may be viewed as a random element in the
metric space $C(\mathcal{C},\mathbb{R}^2)$ of continuous functions
from $\mathcal{C}$ to $\mathbb{R}^2$. We are now in a position to
state Lemma \ref{theo3}.
\begin{lemma}\label{theo3}
Under the assumptions of Theorem \ref{th2}, we have for $z\in
\mathcal{C}$,
\[
(\hat{X}_n(z),Y_n) \stackrel{D}\longrightarrow(X(z),Y),
\]
where $Y\sim N(0,2c(g'(c))^2)$, which is independent of $X(z)$, a
Gaussian stochastic process with mean zero and covariance function
$\operatorname{Cov}(X(z_1),X(z_2))$ equal to
%
%e1.10 ###
%
\begin{equation}\label{e33}
\quad\frac{2}{cz_1z_2[(1+\underline{m}(z_1))(1+\underline
{m}(z_2))-c\underline{m}(z_1)\underline{m}(z_2)]}-\frac{2m(z_1)m(z_2)}{c}.
\end{equation}
\end{lemma}
\begin{remark}\label{re1}
Also, note that $X(z)$ is exactly the weak limit of the stochastic
process $\sqrt{n}(\mathbf{x}_n^T(\bbS-z\bbI)^{-1})\mathbf{x}_n-m_n(z))$ when
$\max_{i}x_{ni}\rightarrow0$, whose covariance function
is
\[
\operatorname{Cov}(X(z_1),X(z_2))=\frac{2(z_2\underline{m}(z_2)-z_1\underline
{m}(z_1))^2}{c^2z_1z_2(z_1-z_2)(\underline{m}(z_1)-\underline{m}(z_2))}
\]
(see \cite{b1} and \cite{p2}).
\end{remark}

We conclude this section by presenting the structure of this work.
In Section \ref{sim}, we present a simulation study to identify when the
asymptotic normality ``kicks in.'' Then we turn to the proof.
To transfer Lemma \ref{theo3} to Theorem \ref{th2} we introduce a
new empirical distribution function
%
%e1.11 ###
%
\begin{equation}\label{f3}F_2^{\bbS}(x)=\sum
_{i=1}^pt_i^2I(\lambda_i\leq
x),
\end{equation}
where
$\bbt=(t_1,\ldots,t_n)^T=\bbU\bar{\mathbf s}/\|\bar{\mathbf s}\|$ and
$\bbU$ is
the eigenvector matrix of $\mathcal{S}$. It turns out that
$F_2^{\bbS}(x)$ and the ESD of $\bbS$ have the same limit, that is,
$F_2^{\bbS}(x)\stackrel{\mathrm{i.p.}}\longrightarrow F_c(x)$. Thus, by
analyticity of $f(x)$,
$\bar{\mathbf s}^Tf(\bolds{\mathcal{S}})\bar{\mathbf s}/\|\bar{\mathbf
s}\|^2$ in
Theorem \ref{th2} is transferred to the Stieljes transform of
$F_2^{\bbS}(x)$,
$\bar{\mathbf s}^T(\bolds{\mathcal{S}}-z\bbI)^{-1}\bar{\mathbf s}/\|\bar
{\mathbf s}\|^2$.
Moreover, note that
%
%e1.12 ###
%
\begin{equation}\label{d21}
\frac{\bar{\mathbf s}^T\bbA^{-1}(z)\bar{\mathbf s}}{1-\bar{\mathbf
s}^T\bbA
^{-1}(z)\bar{\mathbf s}}=\bar{\mathbf s}^T(\bolds{\mathcal{S}}-z\bbI)^{-1}\bar
{\mathbf s}.
\end{equation}
Indeed, this is from the identity (see \cite{s3}, (2.1))
%
%e1.13 ###
%
\begin{equation}\label{g1}
\bbr^T(\bbB+a\bbr\bbr^T)^{-1}=\frac{r^T\bbB^{-1}}{1+a\bbr^T\bbB
^{-1}\bbr},
\end{equation}
where $\bbB$ and $\bbB+a\bbr\bbr^T$ are both invertible,
$\bbr\in\mathbb{R}^p$ and $a\in\mathbb{R}$. The stochastic process
$X_n(z)$ in Lemma \ref{theo3} is then transferred to the
stochastic process $M_n(z)$, where
\[
M_n(z)=\sqrt{n}\biggl(\bar{\mathbf s}^T\bbA^{-1}(z)\bar{\mathbf s}-\frac
{c_nm_n(z)}{1+c_nm_n(z)}\biggr).
\]
The convergence of the stochastic process
$M_n(z)$ is given in Sections \ref{weakcon} and \ref{tightness}. The
proofs of Theorems
\ref{th1} and \ref{th2}, Lemma \ref{theo3} and Remark
\ref{re1} are included in Section \ref{finalproof}. The last section
picks up the
truncation of the underlying r.v.'s and some useful lemmas. At this
point we would like to point out that both this paper and \cite{b2}
deal with Stieljes transform of random variables of interest and use
martingale method to establish CLT. But the random variable of interest
in this paper is a kind of random quadratic forms while \cite{b2} is
concerned with the trace of random matrices.

Throughout this paper, to save notation, $\mathfrak{M}$
may denote different constants on different occasions.

%s2 ###
\section{Simulation study} \label{sim}
In this section, we provide a simulation study to investigate the
performance of normal approximations in Theorem \ref{th1}.
We consider three different populations, the standard normal
distribution, the exponential distribution with parameter 1 and the
Poisson distribution with parameter 1. From each population we generate
5000 samples of order $100\times200$, $200\times400$ and $400\times
800$ matrices, respectively, by routines in \textit{R}. Each $p\times n$
matrix can be regarded as a collection of $n$ observations of
$p$-dimensional vectors ${\mathbf s}$, so we can calculate $T^2$ for each
matrix. Based on 5000 samples, we have 5000 observed $T^2$ which give
us an estimator of the probability
\[
P\biggl(\frac{\sqrt{n}}{\sqrt{2c_n(1-c_n)^{-3}}}\biggl(\frac
{T^2}{n}-c_n(1-c_n)^{-1}\biggr)\leq x\biggr)
\]
by
\[
5000^{-1}\sum I\biggl(\frac{\sqrt{n}}{\sqrt{2c_n(1-c_n)^{-3}}}
\biggl(\frac{T^2}{n}-c_n(1-c_n)^{-1}\biggr)\leq x\biggr).
\]

\begin{figure}

\includegraphics{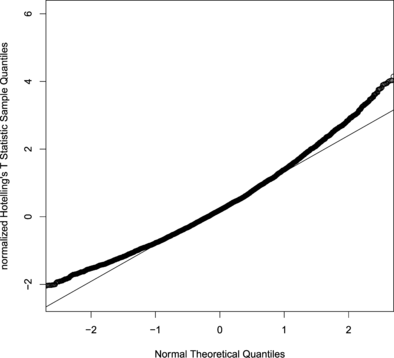}

\caption{Q--Q plot for normal data when $p=100$.}
\label{figure3}\vspace*{3pt}
\end{figure}

\begin{figure}\vspace*{3pt}

\includegraphics{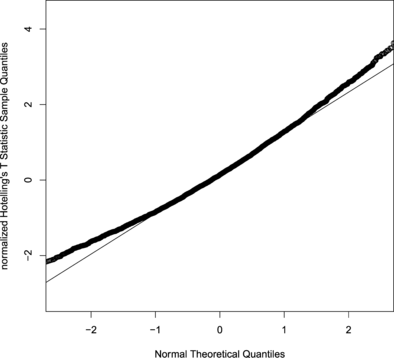}

\caption{Q--Q plot for normal data when $p=200$.}
\label{figure4}
\end{figure}

\begin{figure}

\includegraphics{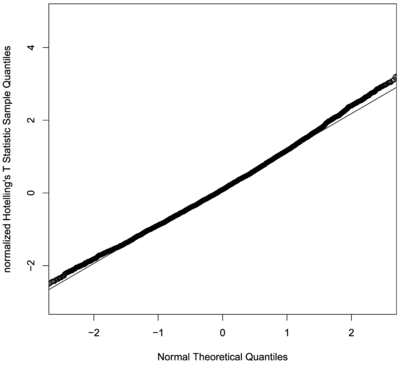}

\caption{Q--Q plot for normal data when $p=400$.}\vspace*{28pt}
\label{figure5}
%
%%f6 ###
%%

\includegraphics{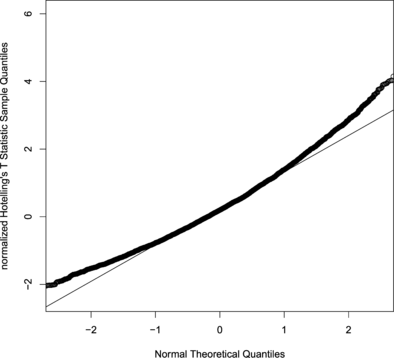}

\caption{Q--Q plot for exponential data when $p=100$.}
\label{figure6}
\end{figure}

\begin{figure}

\includegraphics{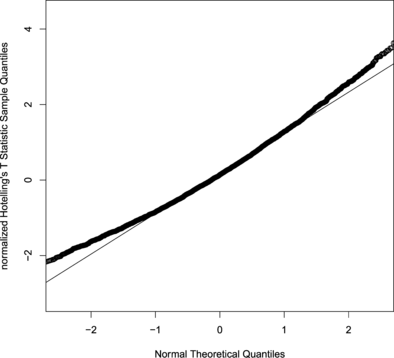}

\caption{Q--Q plot for exponential data when $p=200$.}\vspace*{28pt}
\label{figure7}
%
%%f8 ###
%%

\includegraphics{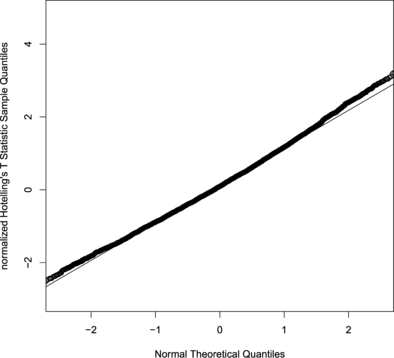}

\caption{Q--Q plot for exponential data when $p=400$.}
\label{figure8}
\end{figure}

\begin{figure}

\includegraphics{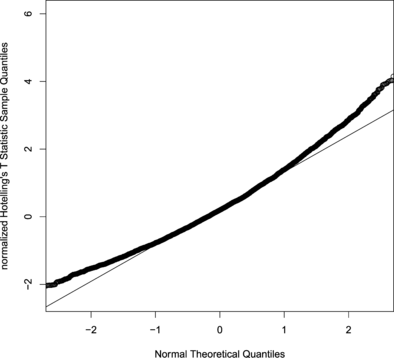}

\caption{Q--Q plot for Poisson data when $p=100$.}\vspace*{28pt}
\label{figure9}
%
%%f10 ###
%%

\includegraphics{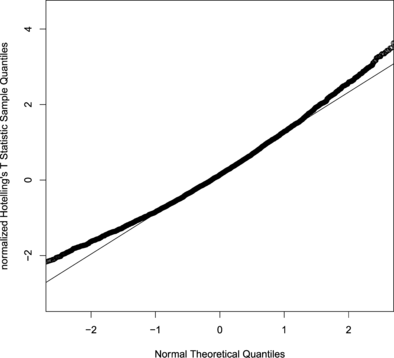}

\caption{Q--Q plot for Poisson data when $p=200$.}
\label{figure10}
\end{figure}

\begin{figure}

\includegraphics{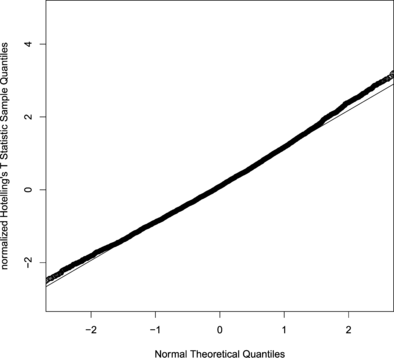}

\caption{Q--Q plot for Poisson data when $p=400$.}
\label{figure11}
\end{figure}

In Figures \ref{figure3}--\ref{figure11}, there are nine curves. In
each figure the horizontal
axis means theoretical quantiles of the standard normal distribution
and the vertical axis indicates sample quantiles of the normalized
Hotelling's $T^2$ statistics. Every curve represents the
quantile-quantile plot for each sampled matrix. From these pictures we
see that the quantiles of $T^2$ get closer to the standard normal one
as the sample size and the dimension increase. Actually, when $p=100$
and $n=200$, normal distributions already ``kick in.''

%
%f3 ###
%
%
%f4 ###
%
%
%f5 ###
%
%
%f7 ###
%
%
%f9 ###
%
%
%f11 ###
%
%

%s3 ###
\section{Weak convergence of the finite-dimensional
distributions}\label{weakcon}

For $z\in\mathcal{C}_n^+$, let $M_n(z)=M_n^{(1)}(z)+M_n^{(2)}(z)$,
where
\[
M_n^{(1)}(z)=\sqrt{n}\bigl(\bar{\mathbf s}^T\bbA^{-1}(z)\bar{\mathbf s}-E\bar
{\mathbf s}^T\bbA^{-1}(z)\bar{\mathbf s}\bigr)
\]
and
\[
M_n^{(2)}(z)=\sqrt{n}\biggl(E\bar{\mathbf s}^T\bbA^{-1}(z)\bar{\mathbf
s}-\frac
{c_nm_n(z)}{1+c_nm_n(z)}\biggr).
\]
In this section the aim is to prove that for any positive integer
$r$ and complex numbers $a_1,\ldots,a_r$,
\[
\sum_{i=1}^ra_iM_n^{(1)}(z_i),\qquad \Im z_i\neq0,
\]
converges in distribution to a Gaussian r.v. and to derive the
asymptotic covariance function. Before proceeding, r.v.'s need to
be truncated. However, we shall postpone the truncation of r.v.'s
until the last section. As a consequence of Lemma \ref{trunlem}, we
assume that the underlying r.v.'s satisfy
%
%e3.1 ###
%
\begin{equation}\label{g8}\quad
|X_{ij}|\leq\varepsilon_n \sqrt{n},\qquad EX_{11}=0,\qquad
E|X_{11}|^2=1,\qquad
E|X_{11}|^4<\infty,
\end{equation}
where $\varepsilon_n$ is a positive sequence which converges to zero
as $n$ goes to infinity.

%s3.1 ###
\subsection{Outline of the argument}

The underlying idea is to write $M_n^{(1)}(z)$ as a sum of
martingale difference sequences and to apply Lemma \ref{lemma3},
CLT for martingale. Define the $\sigma$-field $
{\mathcal F}_j=\sigma({\mathbf s}_1,\ldots,{\mathbf s}_j)$ and let
$E_j(\cdot)=E(\cdot|{\mathcal F}_j)$ and $E_0(\cdot)$ be the
unconditional expectation. We first simplify the martingale
representation of $M_n^{(1)}(z)$ as
$\sum_{j=1}^nY_j(z)+o_p(1)$, where
$Y_j(z)=-2z\underline{m}(z)\* E_{j}(\frac{1}{\sqrt{n}}{\mathbf s}_j^T\bbA
_j^{-1}(z)\bar{\mathbf s}_j)
+z\underline{m}(z)\sqrt{n} E_{j}(\alpha_j(z))$ and $\alpha_j(z)$ and
$\bar{\mathbf s}_j$ are defined in the next subsection. Condition~(ii) in
Lemma \ref{lemma3} is relatively easy to verify. Subsequently,
to identify the asymptotic covariance function of $M_n^{(1)}(z)$,
the following limits in probability need to be determined:
%
%e3.4 ###
%e3.3 ###
%e3.2 ###
%
\begin{eqnarray}
\label{a8}
&\displaystyle \frac{1}{n}\sum_{j=1}^nE_{j-1}[E_{j}({\mathbf s}_j^T\bbA
_j^{-1}(z_1)\bar{\mathbf s}_j)E_j(\bar{\mathbf s}_j^T\bbA
_j^{-1}(z_2){\mathbf s}_j)],&\\
\label{a9}
&\displaystyle \sum_{j=1}^nE_{j-1}[E_{j}({\mathbf s}_j^T\bbA_j^{-1}(z_1)\bar
{\mathbf s}_j)E_{j}(\alpha_j(z_2))],&\\
\label{a10}
&\displaystyle
n\sum_{j=1}^nE_{j-1}[E_{j}(\alpha_j(z_1))E_{j}(\alpha_j(z_2))].&
\end{eqnarray}
As for (\ref{a8}), note that
\[
E_{j-1}[E_{j}({\mathbf s}_j^T\bbA_j^{-1}(z_1)\bar{\mathbf s}_j)E_j(\bar
{\mathbf s}_j^T\bbA_j^{-1}(z_2){\mathbf s}_j)]=E_j(\bar{\mathbf
s}_j^T\bbA
_j^{-1}(z_2))E_{j}(\bbA_j^{-1}(z_1)\bar{\mathbf s}_j)
\]
and $\bar{\mathbf s}_j$ is an average value of all ${\mathbf s}_1,\ldots
,{\mathbf s}_n$
without ${\mathbf s}_j$. Intuitively, the product of two conditional
expectations in the right-hand side of the above formula should be
a multiple of
$\frac{1}{n}\operatorname{tr}[E_{j}(\bbA_j^{-1}(z_1))E_j(\bbA_j^{-1}(z_2))]$. This
turns out to be true. For~(\ref{a10}), a direct calculation
indicates that $E_{j-1}[E_{j}(\alpha_j(z_1))E_{j}(\alpha_j(z_2))]$
involves $\operatorname{tr}[E_j(\bbD_j(z_1))E_j(\bbD_j(z_2))]$ [$\bbD_j(z)$ is
defined in the next subsection]. Then our aim is to transfer it to
$[E_j(\bar{\mathbf s}_j^T\bbA_j^{-1}(z_2))E_{j}(\bbA_j^{-1}(z_1)\bar
{\mathbf s}_j)]^2$
so that the limit of (\ref{a8}) may be used. Essentially, we
expect that (\ref{a8}) and (\ref{a10}) could be reduced to
something like
\[
\frac{1}{n}\sum_{j=1}^nh\biggl(\frac{j-1}{n}\biggr)
\]
for some function $h(x)$. Finally, since the number of ${\mathbf s}_j$
involved in (\ref{a9}) is odd and ${\mathbf s}_j$ is independent of
$\bar{\mathbf s}_j$ we expect that
$\mbox{(\ref{a9})}\stackrel{\mathrm{i.p.}}\rightarrow0$.

%s3.2 ###
\subsection{Notation and estimates}

We first introduce some notation. Let
\begin{eqnarray*}
\bbA_j^{-1}(z)&=&(\bbS-n^{-1}{\mathbf s}_j{\mathbf
s}_j^T-z\bbI)^{-1},\\[-2pt]
\bbA_{ij}^{-1}(z)&=&(\bbS-n^{-1}{\mathbf s}_i{\mathbf
s}_i^T-n^{-1}{\mathbf s}_j{\mathbf s}_j^T-z\bbI)^{-1},
\\[-2pt]
\bar{\mathbf s}_j&=&\bar{\mathbf s}-n^{-1}{\mathbf s}_j,\\[-2pt]
\bbD_j(z)&=&\bbA_j^{-1}(z)\bar{\mathbf s}_j\bar{\mathbf s}_j^T\bbA_j^{-1}(z),
\\[-2pt]
\beta_j(z)&=&\frac{1}{1+({{1}/{n}}){\mathbf s}_j^T \bbA
_j^{-1}(z){\mathbf s}_j},\\[-2pt]
\beta_j^{\rmtr}(z)&=&\frac{1}{1+({1}/{n})\operatorname{tr}\bbA_j^{-1}(z)},\\[-2pt]
b_1(z)&=&\frac{1}{1+({1}/{n})E\operatorname{tr}\bbA_1^{-1}(z)},
\\[-2pt]
\gamma_j(z)&=&\frac{1}{n}{\mathbf s}_j^T \bbA_j^{-1}(z){\mathbf
s}_j-\frac
{1}{n}\operatorname{tr}\bbA_j^{-1}(z),\\[-2pt]
\xi_j(z)&=&\frac{1}{n}{\mathbf s}_j^T \bbA_j^{-1}(z){\mathbf s}_j-\frac
{1}{n}E\operatorname{tr}\bbA_j^{-1}(z),
\\[-2pt]
\alpha_j(z)&=&\frac{1}{n}{\mathbf s}_j^T\bbA_j^{-1}(z)\bar{\mathbf
s}_j\bar
{\mathbf s}_j^T\bbA_j^{-1}(z){\mathbf s}_j
-\frac{1}{n}\bar{\mathbf s}_j^T\bbA_j^{-2}(z)\bar{\mathbf s}_j,
\\[-2pt]
\beta_{ij}(z)&=&\frac{1}{1+({1}/{n}){\mathbf s}_i^T
A_{ij}^{-1}(z){\mathbf s}_i},\\[-2pt]
\beta^{\rmtr}_{ij}(z)&=&\frac{1}{1+({1}/{n})\operatorname{tr}
A_{ij}^{-1}(z)},\\[-2pt]
b_{12}(z)&=&\frac{1}{1+({1}/{n})E\operatorname{tr}A_{12}^{-1}(z)}
\end{eqnarray*}
and
\begin{eqnarray*}
\xi_{ij}(z)&=&\frac{1}{n}{\mathbf s}_i^T
\bbA_{ij}^{-1}(z){\mathbf s}_i-\frac{1}{n}E\operatorname{tr}\bbA_{12}^{-1}(z),\\
\gamma_{ij}(z)&=&\frac{1}{n}{\mathbf s}_i^T
\bbA_{ij}^{-1}(z){\mathbf s}_i-({1}/{n})\operatorname{tr}\bbA_{ij}^{-1}(z).
\end{eqnarray*}

We next list some results to be used later. A direct calculation
indicates that the following equalities are true:
%
%e3.6 ###
%e3.5 ###
%
\begin{eqnarray}
\label{i1}
E({\mathbf s}_1^T\bbA{\mathbf s}_1-\operatorname{tr}\bbA)({\mathbf s}_1^T\bbB{\mathbf
s}_1-\operatorname{tr}\bbB)
&=&
(EX_{11}^4-|EX_{11}^2|^2-2)\sum_{i=1}^pa_{ii}b_{ii}\nonumber\\[-8pt]\\[-8pt]
&&{}+|EX_{11}^2|^2\operatorname{tr}\bbA\bbB^T
+\operatorname{tr}\bbA\bbB;\nonumber
\\
\label{g38}
E[({\mathbf s}_1^T\bbA{\mathbf s}_1-\operatorname{tr}\bbA){\mathbf s}_1^T\bbB\bbr]&=&
EX_{11}^3\sum_{i=1}^pa_{ii}\bbe_i^T\bbB\bbr,
\end{eqnarray}
where $\bbB=(b_{ij})_{p\times p}$ and $\bbA=(a_{ij})_{p\times p}$ are
deterministic complex matrices and $\bbr$ is a deterministic vector.
Here $\bbe_i$ is the vector with the
$i$th element being 1 and zero otherwise. In what follows, to
facilitate the analysis in the subsequent
subsections, we shall assume $v=\Im z> 0$. Note that\vspace*{1pt} $\beta_j(z),
\beta_j^{\rmtr}(z), \beta_{ij}(z), \beta_{ij}^{\rmtr}(z), b_1(z),
b_{12}(z)$ are bounded in absolute value by $|z|/v$ [see \cite{b4},
(3.4)]. From (\ref{g1}) we have
%
%e3.7 ###
%
\begin{eqnarray}\label{i2}
\bbA^{-1}(z)-\bbA_j^{-1}(z)&=&\bbA^{-1}(z)\bigl(\bbA_j(z)-\bbA(z)\bigr)\bbA
_j^{-1}(z)\nonumber\\[-8pt]\\[-8pt]
&=&-\frac{1}{n}\tbbA_j(z)\beta_j(z),\nonumber
\end{eqnarray}
where $\tbbA_j(z)=\bbA_j^{-1}(z){\mathbf s}_j{\mathbf s}_j^T\bbA_j^{-1}(z)$.
From Lemma 2.10 of \cite{b4}, for any matrix $\bbB$,
%
%e3.8 ###
%
\begin{equation}\label{g10}
\bigl|\operatorname{tr}\bigl[\bigl(\bbA^{-1}(z)-\bbA_j^{-1}(z)\bigr)\bbB\bigr]\bigr|
\leq\frac{\|\bbB\|}{v},
\end{equation}
where \mbox{$\|\cdot\|$} denotes the spectral norm of a matrix. Moreover,
Section 4 in \cite{b4} shows that
%
%e3.9 ###
%
\begin{equation}\label{g28}
n^{-k}E|{\operatorname{tr}\bbA_1^{-1}(z)}-E\operatorname{tr}\bbA_1^{-1}(z)|^k=O(n^{-k/2}),\qquad
k\geq2.
\end{equation}
One should also note that (\ref{g28}) is still true when
$\bbA_1^{-1}(z)$ is replaced by $\bbA_{12}^{-1}(z)$.

From now on, we calculate estimates. To simplify the statements,
assume that the spectral norms of nonrandom
$\bbB,\bbB_i,\bbA_i,\bbC$ involved in the equalities
(\ref{a3})--(\ref{a22}) below are all bounded above by a constant.
For $ k\geq2$, it follows from Lemma \ref{lem1}, (\ref{g8}) and
(\ref{g28}) that
%
%e3.10 ###
%
\begin{eqnarray}\label{a3}
n^{-k}E|{\mathbf s}_1^T\bbB{\mathbf s}_1-\operatorname{tr}\bbB|^{k}&=&O(\varepsilon
_n^{2k-4}n^{-1}),\nonumber\\[-8pt]\\[-8pt]
E|\xi_1(z)|^k&=&O(\varepsilon_n^{2k-4}n^{-1})\nonumber
\end{eqnarray}
and that
%
%e3.11 ###
%
\begin{eqnarray}\label{g48}
&&
n^{-k}E|{\mathbf s}_1^T\bbB\bbe_i\bbe_j^T\bbC{\mathbf s}_1|^{k}\nonumber\\
&&\qquad\leq
\mathfrak{M}n^{-k}[E|{\mathbf s}_1^T\bbB\bbe_i\bbe_j^T\bbC{\mathbf
s}_1-\operatorname{tr}(\bbB\bbe_i\bbe_j^T\bbC)|^{k}+E|\bbe_j^T\bbC\bbB\bbe
_i|^k]\\
&&\qquad=O(\varepsilon_n^{2k-4}n^{-2})\nonumber.
\end{eqnarray}

We shall establish the estimates (\ref{g2})--(\ref{g9}) below:
%
%e3.13 ###
%e3.12 ###
%
\begin{eqnarray}\label{g2}
E|{\mathbf s}_1^T\bbB\bar{\mathbf s}_1|^k&=&O\bigl(n^{(
{k}/{2}-2)}\varepsilon
_n^{k-4}\bigr),\qquad
k\geq4,\nonumber\\[-8pt]\\[-8pt]
E|\alpha_1(z)|^k&=&O(n^{-2}\varepsilon_n^{2k-4}),\qquad
k\geq2,\nonumber
\\
\label{a4}
E|{\mathbf s}_{1}^T\bbB{\mathbf
s}_{2}|^k&=&O(n^{k-2}\varepsilon_n^{k-4}),\qquad
k\geq4,
\end{eqnarray}
and for $m\geq0, q\geq1,0\leq r\leq2$,
%
%e3.14 ###
%
\begin{equation}\label{g9}
E\Biggl|\prod_{i=1}^m\frac{1}{n}{\mathbf s}_1^T\bbA_i{\mathbf s}_1\prod
_{j=1}^q\frac{1}{n}({\mathbf s}_1^T\bbB_j{\mathbf s}_1-\operatorname{tr}\bbB
_j)({\mathbf s}_1^T\bbC_l\bar{{\mathbf s}}_1)^r\Biggr|
=O\bigl(n^{-{1}/{2}}\varepsilon_n^{(q-2)\vee0}\bigr).\hspace*{-28pt}
\end{equation}
One should note that (\ref{g2}) and (\ref{a4}) also give the
estimates for $k=2$. For example,
%
%e3.15 ###
%
\begin{equation}\label{g39}
E|{\mathbf s}_1^T\bbB\bar{\mathbf s}_1|^2\leq
(E|{\mathbf s}_{1}^T\bbB\bar{\mathbf s}_{1}|^4)^{1/2}= O(1).
\end{equation}
In addition, from (\ref{a3}) and (\ref{a4}) we also conclude that
%
%e3.16 ###
%
\begin{eqnarray}\label{g41}
E|n^{-1}{\mathbf s}_1^T\bbB{\mathbf s}_1{\mathbf s}_1^T\bbC{\mathbf
s}_2|^4&\leq&
\mathfrak{M}
E|n^{-1}({\mathbf s}_1^T\bbB{\mathbf s}_1-\operatorname{tr}\bbB){\mathbf s}_1^T\bbC
{\mathbf s}_2|^4\nonumber\\
&&{}+\mathfrak{M}E|{\mathbf s}_1^T\bbC{\mathbf s}_2|^4\\
&=&O(n^{5/2}).\nonumber
\end{eqnarray}

Consider (\ref{g2}) first. Note that for $k\geq4$
%
%e3.17 ###
%
\begin{eqnarray}\label{7}
E|\bar{\mathbf s}_1^T\bar{\mathbf s}_1|^k&\leq&
\frac{\mathfrak{M}}{n^{2k}}\Biggl[E\Biggl|\sum_{i=2
}^n{\mathbf s}_i^T{\mathbf s}_i\Biggr|^k+E\biggl|\sum_{i_1\neq i_2, i_1>1,
i_2>1}{\mathbf s}_{i_1}^T{\mathbf s}_{i_2}\biggr|^k\Biggr]\nonumber\\[-8pt]\\[-8pt]
&=&O(1).\nonumber
\end{eqnarray}
Indeed, applying Lemma \ref{lem3} twice gives
\begin{eqnarray*}
E\Biggl|\sum_{i=2}^n{\mathbf s}_i^T{\mathbf s}_i\Biggr|^k&\leq&
\mathfrak{M}E\Biggl|\sum_{i=2}^n\bigl({\mathbf s}_i^T{\mathbf s}_i-E({\mathbf
s}_i^T{\mathbf s}_i)\bigr)\Biggr|^k+\mathfrak{M}\Biggl|\sum_{i=2}^nE({\mathbf
s}_i^T{\mathbf s}_i)\Biggr|^k\\
&\leq& \mathfrak{M}\Biggl(\sum_{i=2}^n
E\bigl({\mathbf s}_i^T{\mathbf s}_i-E({\mathbf s}_i^T{\mathbf s}_i)\bigr)^2
\Biggr)^{k/2}\\
&&{}+\mathfrak{M}\sum_{i=2}^nE|{\mathbf s}_i^T{\mathbf
s}_i-E({\mathbf s}_i^T{\mathbf s}_i)|^k+\mathfrak{M}n^{2k}\\
&\leq& \mathfrak{M}n^{k}+
\mathfrak{M}n\Biggl[\Biggl(\sum_{m=1}^pE(X_{m2}^2-1)^2\Biggr)^{k/2}+\sum
_{m=1}^pE|X_{m2}^2-1|^k\Biggr]\\
&&{}+\mathfrak{M}n^{2k}\\
&\leq& \mathfrak{M}n^{2k},
\end{eqnarray*}
while using Lemma \ref{lem3} three times we obtain
\begin{eqnarray*}
E\biggl|\sum_{i_1\neq i_2, i_1>1,i_2>1}{\mathbf s}_{i_1}^T{\mathbf s}_{i_2}\biggr|^k
&\leq& n^{k-1}\sum_{i_1>1}E\biggl|\sum_{i_2>1 ,i_1\neq
i_2}{\mathbf s}_{i_1}^T{\mathbf s}_{i_2}\biggr|^k\leq
n^{k}E\Biggl|\sum_{i=3}^n{\mathbf s}_{2}^T{\mathbf s}_{i}\Biggr|^k \\
&\leq&
\mathfrak{M}n^{k}E\Biggl|\sum_{i=3}^nE[({\mathbf s}_{2}^T{\mathbf
s}_{i})^2|\mathcal{G}_{i-1}]\Biggr|^{k/2}+\mathfrak{M}n^{k}\sum
_{i=3}^nE|{\mathbf s}_{2}^T{\mathbf s}_{i}|^k\\
&\leq&
\mathfrak{M}\bigl[n^{({3}/{2})k}E|{\mathbf s}_{2}^T{\mathbf
s}_{2}-E{\mathbf s}_{2}^T{\mathbf
s}_{2}|^{k/2}+n^{2k}+n^{k+1}E|{\mathbf s}_{2}^T{\mathbf
s}_{3}|^k\bigr]\\
&=&O(n^{2k}),
\end{eqnarray*}
where $\mathcal{G}_{i}=\sigma({\mathbf s}_2,\ldots,{\mathbf s}_i)$. It follows
from (\ref{7}) that for $k\geq4$
%
%e3.18 ###
%
\begin{equation}\label{a22}
E|\bar{\mathbf s}_1^T\bbB\bar{\mathbf s}_1|^{k}
=E\|\bar{\mathbf s}_1^T\bbB\bar{\mathbf s}_1\|^{k}\leq
E(\|\bar{\mathbf s}_1^T\|\|\bbB\|\|\bar{\mathbf s}_1\|)^{k}\leq
\mathfrak{M}E|\bar{\mathbf s}_1^T\bar{\mathbf
s}_1|^{k}\leq\mathfrak{M},\hspace*{-28pt}
\end{equation}
where \mbox{$\|\cdot\|$} denotes the spectral norm of a matrix. This,
together with Lemma \ref{lem1}, ensures that for $k\geq4$
\begin{eqnarray*}
E|{\mathbf s}_1^T\bbB\bar{\mathbf s}_1|^k
&=&E|{\mathbf s}_1^T\bbB\bar{\mathbf s}_1\bar{\mathbf s}_1^T\bbB
^*{\mathbf s}_1|^{{k}/{2}}\\
&\leq&
\mathfrak{M}E|{\mathbf s}_1^T\bbB\bar{\mathbf s}_1\bar{\mathbf
s}_1^T\bbB
^*{\mathbf s}_1-\bar{\mathbf s}_1^T\bbB^*\bbB\bar{\mathbf s}_1|^{{k}/{2}}
+\mathfrak{M}E|\bar{\mathbf s}_1^T\bbB^*\bbB\bar{\mathbf s}_1|^{
{k}/{2}}\\
&\leq&[\mathfrak{M}n^{{k}/{2}-2}\varepsilon_n^{k-4}+\mathfrak
{M}]E|\bar{\mathbf s}_1^T\bbB^*\bbB\bar{\mathbf s}_1|^{
{k}/{2}}+\mathfrak{M}\\
&\leq& \mathfrak{M}n^{{k}/{2}-2}\varepsilon_n^{k-4},
\end{eqnarray*}
which gives the first estimate in (\ref{g2}) as well as the order of
$E|\alpha_1(z)|^k$.

Second, consider (\ref{a4}). Let $\mathbf{y}=(y_1,\ldots,y_p)^T=\bbB
{\mathbf s}_2$ and then, by Lemma \ref{lem3} and (\ref{a3}), for $k\geq4$,
%
%e3.19 ###
%
\begin{eqnarray}\label{g16}
E|{\mathbf s}_{1}^T\mathbf{y}|^k&\leq&
\mathfrak{M}E\Biggl(\sum_{m=1}^p|y_m|^2\Biggr)^{k/2}+M\sum
_{m=1}^pE|X_{m1}|^kE|y_m|^k\nonumber\\
&\leq&\mathfrak{M}E|\mathbf{y}^*\mathbf{y}|^{k/2}+\mathfrak{M}n^{
{k}/{2}-2}\varepsilon_n^{k-4}E|\mathbf{y}^*\mathbf{y}|^{k/2}\nonumber\\
&\leq&
\mathfrak{M}(1+n^{{k}/{2}-2}\varepsilon_n^{k-4})E|{\mathbf
s}_2^T\bbB^*\bbB{\mathbf s}_2-\operatorname{tr}(\bbB^*\bbB)|^{k/2}\\
&&{}+\mathfrak
{M}n^{k/2}+\mathfrak{M}n^{k-2}\varepsilon_n^{k-4}
\nonumber\\
&=&O(n^{k-2}\varepsilon_n^{k-4}),\nonumber
\end{eqnarray}
where we also use the fact that for $k\geq4$
\[
\sum_{m}|y_m|^k\leq\biggl(\sum_{m}|y_m|^2\biggr)^{k/2}.
\]

As for (\ref{g9}), if $m=0$ and $r=0$, then (\ref{g9}) directly
follows from (\ref{a3}) and the H\"{o}lder inequality. If $m\geq
1$ and $r=0$, then by induction on $m$ we have
\begin{eqnarray*}
&&E\Biggl|\prod_{i=1}^m\frac{1}{n}{\mathbf s}_1^T\bbA_i{\mathbf s}_1\prod
_{j=1}^q\frac{1}{n}({\mathbf s}_1^T\bbB_j{\mathbf s}_1-\operatorname{tr}\bbB_j)\Biggr|
\\
&&\qquad\leq
E\Biggl|\prod_{i=1}^{m-1}\frac{1}{n}{\mathbf s}_1^T\bbA_i{\mathbf s}_1\frac
{1}{n}({\mathbf s}_1^T\bbA_m{\mathbf s}_1-\operatorname{tr}\bbA_m)\prod
_{j=1}^q\frac{1}{n}({\mathbf s}_1^T\bbB_j{\mathbf s}_1-\operatorname{tr}\bbB_j)\Biggr|\\
&&\qquad\quad{}+\mathfrak{M}E\Biggl|\prod_{i=1}^{m-1}\frac{1}{n}\operatorname{tr}\bbA
_i\prod_{j=1}^q\frac{1}{n}({\mathbf s}_1^T\bbB_j{\mathbf s}_1-\operatorname{tr}\bbB
_j)\Biggr|\\
&&\qquad=O\bigl(n^{-{1}/{2}}\varepsilon_n^{(q-2)\vee0}\bigr).
\end{eqnarray*}
Repeating the argument above gives
\[
E\Biggl|\prod_{i=1}^m\frac{1}{n}{\mathbf s}_1^T\bbA_i{\mathbf s}_1\prod
_{j=1}^q\frac{1}{n}({\mathbf s}_1^T\bbB_j{\mathbf s}_1-\operatorname{tr}\bbB
_j)\Biggr|^2=O\bigl(n^{-1}\varepsilon_n^{(2q-4)\vee
0}\bigr)
\]
[$m=0$ by (\ref{a3}) and $m\geq1$ by induction]. Thus, for the
case $m\geq1$ and $2\geq r\geq1$, by (\ref{g2}) we obtain
\begin{eqnarray*}
&&E\Biggl|\prod_{i=1}^m\frac{1}{n}{\mathbf s}_1^T\bbA_i{\mathbf s}_1\prod
_{j=1}^q\frac{1}{n}({\mathbf s}_1^T\bbB_j{\mathbf s}_1-\operatorname{tr}\bbB
_j)({\mathbf s}_1^T\bbC_1\bar{{\mathbf s}}_1)^r\Biggr|
\\
&&\qquad\leq\Biggl(E\Biggl|\prod_{i=1}^m\frac{1}{n}{\mathbf s}_1^T\bbA
_i{\mathbf s}_1\prod_{j=1}^q\frac{1}{n}({\mathbf s}_1^T\bbB_j{\mathbf
s}_1-\operatorname{tr}\bbB_j)\Biggr|^2E|{\mathbf s}_1^T\bbC_1\bar{{\mathbf s}}_1|^{2r}\Biggr)^{1/2}
\\
&&\qquad=O\bigl(n^{-{1}/{2}}\varepsilon_n^{(q-2)\vee0}\bigr).
\end{eqnarray*}
When $m=0$ and $2\geq r\geq1$, (\ref{g9}) can be obtained
similarly. Thus, we have proved (\ref{g9}).

%s3.3 ###
\subsection{The simplification of $M_n^{(1)}(z)$} \label{simplif}

To develop CLT for $M_n^{(1)}(z)$, we write it
as a sum of martingale difference sequences. When simplifying such
a martingale representation, a well-known trick is to use the fact
that
%
%e3.20 ###
%
\begin{equation}\label{h1}
E_j[h(\operatorname{tr}\bbA_j^{-1}(z))]=E_{j-1}[h(\operatorname{tr}\bbA_j^{-1}(z))],
\end{equation}
where $h(x)$ is some function. For example, when
$h(x)=1/(1+n^{-1}x)$, (\ref{h1}) becomes
$E_j(\beta_j^{\rmtr})=E_{j-1}(\beta_j^{\rmtr})$.

Notice that
$E_j(\bar{\mathbf s}_j^T\bbA_j^{-1}(z)\bar{\mathbf s}_j)=E_{j-1}(\bar
{\mathbf s}_j^T\bbA_j^{-1}(z)\bar{\mathbf s}_j)$.
We then write
%
%e3.21 ###
%
\begin{eqnarray}\label{a1}
M_n^{(1)}(z)&=&\sqrt{n}\sum_{j=1}^n[E_{j}(\bar{\mathbf s}^T\bbA
^{-1}(z)\bar{\mathbf s})-E_{j-1}(\bar{\mathbf s}^T\bbA^{-1}(z)\bar
{\mathbf s})]\nonumber\\
&=&\sqrt{n}\sum_{j=1}^n\bigl[
E_{j}\bigl(\bar{\mathbf s}^T\bbA^{-1}(z)\bar{\mathbf s}-\bar{\mathbf
s}_j^T\bbA
_j^{-1}(z)\bar{\mathbf s}_j\bigr)\nonumber\\[-8pt]\\[-8pt]
&&\hspace*{33.9pt}{}
-E_{j-1}\bigl(\bar{\mathbf s}^T\bbA^{-1}(z)\bar{\mathbf s}-\bar{\mathbf
s}_j^T\bbA
_j^{-1}(z)\bar{\mathbf s}_j\bigr)\bigr]\nonumber\\
&=&\sqrt{n}\sum_{j=1}^n
[(E_{j}-E_{j-1})(a_{n1}+a_{n2}+a_{n3})],\nonumber
\end{eqnarray}
where
\begin{eqnarray*}
a_{n1}&=&(\bar{\mathbf s}-\bar{\mathbf s}_j)^T\bbA^{-1}(z)\bar{\mathbf
s},\qquad
a_{n2}=\bar{\mathbf s}_j^T\bigl(\bbA^{-1}(z)-\bbA_j^{-1}(z)\bigr)\bar{\mathbf
s},\\
a_{n3}&=&\bar{\mathbf s}_j^T\bbA_j^{-1}(z)(\bar{\mathbf s}-\bar{\mathbf s}_j).
\end{eqnarray*}
The above sum involving $a_{n1}$ and $a_{n2}$ will be further
simplified below.

First, splitting $\bbA^{-1}(z)$ into the sum of
$\bbA^{-1}(z)-\bbA_j^{-1}(z)$ and $\bbA_j^{-1}(z)$ and splitting
$\bar{\mathbf s}$ into the sum of $\bar{\mathbf s}_j$ and ${\mathbf
s}_j/n$, by
(\ref{i2}) we then have
%
%e3.22 ###
%
\begin{equation}
a_{n1}=a_{n1}^{(1)}+a_{n1}^{(2)}+a_{n1}^{(3)}+a_{n1}^{(4)},
\end{equation}
where
\[
a_{n1}^{(1)}=-\frac{1}{n^3}({\mathbf s}_j^T \bbA_j^{-1}(z){\mathbf
s}_j)^2\beta_j(z),\qquad a_{n1}^{(2)}=-
\frac{1}{n^2}{\mathbf s}_j^T \tbbA_j(z)\bar{{\mathbf s}}_j\beta_j(z)
\]
and
\[
a_{n1}^{(3)}=\frac{1}{n^2}{\mathbf s}_j^T \bbA_j^{-1}(z){\mathbf
s}_j,\qquad
a_{n1}^{(4)}=\frac{1}{n}{\mathbf s}_j^T \bbA_j^{-1}(z)\bar{\mathbf s}_j.
\]
Using (\ref{h1}) and
%
%e3.23 ###
%
\begin{equation}
\label{a2}
\beta_j(z)=\beta_j^{\rmtr}(z)-\beta_j(z)\beta_j^{\rmtr}(z)\gamma_j(z),
\end{equation}
we have
\begin{eqnarray*}
&&(E_{j}-E_{j-1})\bigl(a_{n1}^{(1)}\bigr)\\
&&\qquad=(E_{j}-E_{j-1})\biggl[\frac{1}{n^3}({\mathbf s}_j^T
\bbA_j^{-1}(z){\mathbf s}_j)^2\beta_j^{\rmtr}(z)\biggr]-\varsigma_n
\\
&&\qquad=(E_{j}-E_{j-1})\biggl[\frac{1}{n}\gamma_j^2(z)\beta
_j^{\rmtr}(z)\biggr]\\
&&\qquad\quad{}+(E_{j}-E_{j-1})\biggl[\frac{2}{n^2}\gamma_j(z)\beta
_j^{\rmtr}(z)\operatorname{tr}\bbA_j^{-1}(z)\biggr]-\varsigma_n,
\end{eqnarray*}
where $\varsigma_n=(E_{j}-E_{j-1})\frac{1}{n^3}({\mathbf s}_j^T
\bbA_j^{-1}(z){\mathbf s}_j)^2\beta_j(z)\beta_j^{\rmtr}(z)\gamma_j(z)$. This,
together with (\ref{g9}), shows that
\begin{eqnarray*}
&&E\Biggl|\sqrt{n}\sum_{j=1}^n(E_{j}-E_{j-1})\bigl(a_{n1}^{(1)}\bigr)\Biggr|^2\\
&&\qquad=n\sum
_{j=1}^nE\bigl|(E_{j}-E_{j-1})\bigl(a_{n1}^{(1)}\bigr)\bigr|^2
\\
&&\qquad\leq
\mathfrak{M}E|\gamma_1(z)|^4+\mathfrak{}E|\gamma_1(z)|^2+\mathfrak
{M}E\biggl|\gamma_1(z)\frac{1}{n^2}({\mathbf s}_1^T
\bbA_1^{-1}(z){\mathbf s}_1)^2\biggr|^2 \\
&&\qquad=O(n^{-{1}/{2}}),
\end{eqnarray*}
which gives
\[
\sqrt{n}\sum_{j=1}^n(E_{j}-E_{j-1})\bigl(a_{n1}^{(1)}\bigr)\stackrel
{\mathrm{i.p.}}\longrightarrow0.
\]
By (\ref{a3}) it is a simple matter to verify that
\[
\sqrt{n}\sum_{j=1}^n(E_{j}-E_{j-1})\bigl(a_{n1}^{(3)}\bigr)\stackrel
{\mathrm{i.p.}}\longrightarrow
0.
\]
Appealing to (\ref{g9}) we have
\[
E\Biggl|\sum_{j=1}^n(E_{j}-E_{j-1})\gamma_j(z)\frac{1}{\sqrt
{n}}{\mathbf s}_j^T\bbA_j^{-1}(z)\bar{\mathbf s}_j\beta_j^{\rmtr}(z)\Biggr|^2=
O(n^{-1/2})
\]
and
\[
E\Biggl|\sqrt{n}\sum_{j=1}^n(E_{j}-E_{j-1})\frac{1}{n^2}{\mathbf s}_j^T
\tbbA_j(z)\bar{{\mathbf s}}_j\beta_j(z)\gamma_j(z)\beta
_j^{\rmtr}(z)\Biggr|^2=O(n^{-1/2}),
\]
which, together with (\ref{a2}), leads to
\begin{eqnarray*}
&&
\sqrt{n}\sum_{j=1}^n(E_{j}-E_{j-1})\bigl(a_{n1}^{(2)}\bigr)\\
&&\qquad=-\sum
_{j=1}^nE_{j}\biggl[\bigl(1-\beta_j^{\rmtr}(z)\bigr)\frac{1}{\sqrt{n}}{\mathbf s}_j^T\bbA
_j^{-1}(z)\bar{\mathbf s}_j\biggr]+o_p(1).
\end{eqnarray*}
This ensures that
%
%e3.24 ###
%
\begin{eqnarray}\label{a5}
&&\sqrt{n}\sum_{j=1}^n(E_{j}-E_{j-1})(a_{n1})\nonumber
\\
&&\qquad=\sum
_{j=1}^nE_{j}\biggl(\beta_j^{\rmtr}(z)\frac{1}{\sqrt{n}}{\mathbf s}_j^T\bbA
_j^{-1}(z)\bar{\mathbf s}_j\biggr)+o_p(1)\\
&&\qquad=-z\underline{m}(z)\sum_{j=1}^nE_{j}\biggl(\frac{1}{\sqrt
{n}}{\mathbf s}_j^T\bbA_j^{-1}(z)\bar{\mathbf
s}_j\biggr)+o_p(1),\nonumber
\end{eqnarray}
because, by (2.17) in \cite{b2}, (\ref{g28}) and (\ref{g2}),
%
%e3.25 ###
%
\begin{equation}\label{g40}
E\bigl|\bigl(\beta_j^{\rmtr}(z)+z\underline{m}(z)\bigr){\mathbf s}_j^T\bbA_j^{-1}(z)\bar
{\mathbf s}_j\bigr|^2=o(1).
\end{equation}

Second, splitting $\bar{\mathbf s}$ into the sum of $\bar{\mathbf s}_j$ and
${\mathbf s}_j/n$ further gives
\[
a_{n2}=-\frac{1}{n^2}\bar{\mathbf s}_j^T\tbbA_j(z){\mathbf s}_j\beta
_j(z)-\frac{1}{n}\bar{\mathbf s}_j^T\tbbA_j(z)\bar{\mathbf s}_j\beta_j(z)
\]
and thus, as in treating $a_{n1}^{(2)}$, we have
\begin{eqnarray*}
&&\sqrt{n}\sum_{j=1}^n(E_{j}-E_{j-1})(a_{n2})\\
&&\qquad=-\sum
_{j=1}^n(E_{j}-E_{j-1})\biggl[\bigl(1-\beta_j^{\rmtr}(z)\bigr)\frac{1}{\sqrt{n}}\bar
{\mathbf s}_j^T\bbA_j^{-1}(z){\mathbf s}_j\biggr]
\\
&&\qquad\quad{}
-\frac{1}{\sqrt{n}}\sum
_{j=1}^n(E_{j}-E_{j-1})[\bar{\mathbf s}_j^T\tbbA_j(z)\bar{\mathbf
s}_j\beta_j^{\rmtr}(z)]+o_p(1)
\\
&&\qquad=-\bigl(1+z\underline{m}(z)\bigr)\sum_{j=1}^nE_{j}\biggl(\frac{1}{\sqrt
{n}}\bar{\mathbf s}_j^T\bbA_j^{-1}(z){\mathbf s}_j\biggr)\\
&&\qquad\quad{}+z\underline{m}(z)\sum_{j=1}^n\sqrt{n}E_{j}(\alpha_j(z))+o_p(1),
\end{eqnarray*}
where in the last step we also use the estimate
\begin{eqnarray*}
E\bigl|\bigl(\beta_j^{\rmtr}(z)+z\underline{m}(z)\bigr)\alpha_j(z)\bigr|^2&=&E\bigl[E\bigl(\bigl|\bigl(\beta
_j^{\rmtr}(z)+z\underline{m}(z)\bigr)\alpha_j(z)\bigr|^2|\sigma({\mathbf s}_i,i\neq
j)\bigr)\bigr]
\\
&=&E\bigl[|\beta_j^{\rmtr}(z)+z\underline{m}(z)|^2E\bigl(|\alpha_j(z)|^2|\sigma
({\mathbf s}_i,i\neq
j)\bigr)\bigr]\\
&=&o(n^{-2}),
\end{eqnarray*}
which is from (2.17) in \cite{b2}, (\ref{g28}) and (\ref{g2}).

Recalling
$Y_j(z)=-2z\underline{m}(z)E_{j}(\frac{1}{\sqrt{n}}{\mathbf s}_j^T\bbA
_j^{-1}(z)\bar{\mathbf s}_j)
+z\underline{m}(z)\sqrt{n}E_{j}(\alpha_j(z))$, so far we have
proved
\[
M_n^{(1)}(z)=\sum_{j=1}^nY_j(z)+o_p(1).
\]
Consequently, for finite dimension convergence of $M_n^{(1)}(z)$,
we need consider only the sum
%
%e3.26 ###
%
\begin{equation}\label{a7}
\sum_{i=1}^ra_i\sum_{j=1}^nY_j(z_i)=\sum
_{j=1}^n\sum_{i=1}^ra_iY_j(z_i).
\end{equation}

Next we verify condition (ii) of Lemma \ref{lemma3}. Recalling
$\bbD_j(z)=\bbA_j^{-1}(z)\bar{\mathbf s}_j\bar{\mathbf s}_j^T\times\break\bbA
_j^{-1}(z)$, write
\[
\alpha_j(z)=\alpha_j^{(1)}(z)+\alpha_j^{(2)}(z)+\alpha_j^{(3)}(z),
\]
where
\begin{eqnarray*}
\alpha_j^{(3)}(z)&=&\frac{1}{n}\sum_{h\neq
l}\bbe_h^T\bbD_j(z)\bbe_lX_{hj}X_{lj},
\\
\alpha_j^{(2)}(z)&=&\frac{1}{n}\sum_{h=1}^p\bbe_h^T\bbD
_j(z)\bbe_h[X_{hj}^2I(|X_{hj}|\leq
\log n)-EX_{hj}^2I(|X_{hj}|\leq\log n)]
\end{eqnarray*}
and
\[
\alpha_j^{(1)}(z)=\frac{1}{n}\sum_{h=1}^p\bbe_h^T\bbD
_j(z)\bbe_h[X_{hj}^2I(|X_{hj}|>
\log n)-EX_{hj}^2I(|X_{hj}|> \log n)].
\]
Lemma \ref{lem4} and (\ref{a22}) show that
$E|\alpha_j^{(3)}(z)|^4= O(n^{-4})$. Lemma \ref{lem3} and
(\ref{a22}) give $E|\alpha_j^{(2)}(z)|^4= O(n^{-4}(\log n)^4)$
because
%
%e3.27 ###
%
\begin{eqnarray}\label{c4}
\sum_{h=1}^n|\bbe_h^T\bbD_j(z)\bbe_h|^k&\leq&
\Biggl|\sum_{h=1}^n\bar{\mathbf s}_j^T\bbA_j^{-1}(\bar{z})\bbe_h\bbe
_h^T\bbA_j^{-1}(z)\bar{\mathbf s}_j\Biggr|^k\nonumber\\[-8pt]\\[-8pt]
&=&(\bar{\mathbf s}_j^T\bbA
_j^{-1}(\bar
{z})\bbA_j^{-1}(z)\bar{\mathbf s}_j)^k,\nonumber
\end{eqnarray}
where $k=2$ or $4$ and $\bbA_j^{-1}(\bar{z})$ denotes the complex
conjugate of $\bbA_j^{-1}(z)$. We conclude from (\ref{c4}) and
$EX_{11}^4I(|X_{11}|> \log n)\rightarrow0$ that
$E|\alpha_j^{(1)}(z)|^2=o(n^{-2})$. Therefore, we obtain
\begin{eqnarray*}
&&\sum_{j=1}^nE\Biggl|\sum_{i=1}^ra_iY_j(z_i)
\Biggr|^2I\Biggl(\Biggl|\sum_{i=1}^ra_iY_j(z_i)\Biggr|\geq\varepsilon\Biggr)\\
&&\qquad\leq
4\sum_{j=1}^n\sum_{h=1}^4E\Biggl|\sum
_{i=1}^ra_iY_j^{(h)}(z_i)\Biggr|^2I\Biggl(\Biggl|\sum
_{i=1}^ra_iY_j^{(h)}(z_i)\Biggr|\geq\varepsilon/4\Biggr)
\\
&&\qquad\leq
\frac{\mathfrak{M}}{\varepsilon^2}\sum_{j=1}^n\sum
_{h=2}^4E\Biggl|\sum_{i=1}^ra_iY_j^{(h)}(z_i)\Biggr|^4+4\sum
_{j=1}^nE\Biggl|\sum_{i=1}^ra_iY_j^{(1)}(z_i)
\Biggr|^2\rightarrow
0,
\end{eqnarray*}
where
$Y_j^{(h)}(z)=z\underline{m}(z)\sqrt{n}E_{j}(\alpha_j^{(h)}(z)),h=1,2,3$
and
$Y_j^{(4)}(z)=-2z\underline{m}(z)\times E_{j}(\frac{1}{\sqrt
{n}}{\mathbf s}_j^T\bbA_j^{-1}(z)\bar{\mathbf s}_j)$.
Here we also use $E|Y_j^{(4)}(z)|^4=O(n^{-2})$ by (\ref{g2}). Thus,
the condition (ii) of Lemma \ref{lemma3} is satisfied. Hence, the next
task is to
find, for
$z_1,z_2 \in\mathbb{C}\setminus\mathbb{R}$, the limit in
probability of
%
%e3.28 ###
%
\begin{equation}
\label{i3}
\sum_{j=1}^nE_{j-1}(Y_j(z_1)Y_j(z_2)).
\end{equation}
To this end, it is enough to find the limits in probability for
(\ref{a8}), (\ref{a9}) and (\ref{a10}).

The limits of (\ref{a8})--(\ref{a10}) and finally (\ref{i3}) will be determined in the
subsequent subsections.

%s3.4 ###
\subsection{\texorpdfstring{The limit of (\protect\ref{a8})}{The limit of (3.2)}}

Our aim is to prove that
%
%e3.29 ###
%
\begin{eqnarray}
\label{c1}\mbox{(\ref{a8})}&=&\frac{z_1z_2\underline{m}(z_1)\underline
{m}(z_2)}{n}\sum_{j=1}^n\frac{j-1}{n^2}\operatorname{tr}(E_j(\bbA
_{j}^{-1}(z_2))E_j(\bbA_{j}^{-1}(z_1)))\nonumber\\[-8pt]\\[-8pt]
&&{}+o_p(1).\nonumber
\end{eqnarray}
The strategy is to first replace $\bar{\mathbf s}_j$ by
$\frac{1}{n}\sum_{i\neq j}^n{\mathbf s}_i$, then replace the
resulting quadratic forms in terms of ${\mathbf s}_i$ by its
corresponding trace and $\beta_{ij}(z_2)$ by its corresponding
limit.

To this end, introduce $\underline{\bbA}^{-1}_j(z)$ and
$\underline{\bar{\mathbf s}}_j$ like $\bbA^{-1}_j(z)$ and $\bar{\mathbf s}_j$,
respectively, but $\underline{\bbA}^{-1}_j(z)$ and
$\underline{\bar{\mathbf s}}_j$ are now defined by
${\mathbf s}_{1},\ldots,{\mathbf s}_{j-1},\underline{{\mathbf
s}}_{j+1},\ldots
,\underline{{\mathbf s}}_n$
instead of ${\mathbf s}_{1},\ldots,{\mathbf s}_{j-1},{\mathbf
s}_{j+1},\ldots
,{\mathbf s}_n$.
Here $\{\underline{{\mathbf s}}_{j+1},\ldots,\underline{{\mathbf s}}_n\}
$ are
i.i.d. copies of ${\mathbf s}_1$ and independent of
$\{{\mathbf s}_j,j=1,\ldots,n\}$. Therefore, (\ref{a8}) is equal to
\[
\frac{1}{n}\sum_{j=1}^n\operatorname{tr}[E_{j}(\bbA_j^{-1}(z_1)\bar{\mathbf
s}_j)E_{j}(\bar{\mathbf s}_j^T\bbA_j^{-1}(z_2))]
=\frac{1}{n}\sum_{j=1}^nE_j[\bar{\mathbf s}_j^T\bbA
_j^{-1}(z_2)\underline{\bbA}_j^{-1}(z_1)\underline{\bar{\mathbf s}}_j].
\]
Applying $\bar{\mathbf s}_j=\frac{1}{n}\sum_{i\neq j}^n{\mathbf s}_i$ and
(\ref{g1}) further gives
%
%e3.30 ###
%
\begin{equation}\label{c5}\qquad
E_j[\bar{\mathbf s}_j^T\bbA_j^{-1}(z_2)\underline{\bbA
}_j^{-1}(z_1)\underline{\bar{\mathbf s}}_j]
=\frac{1}{n}\sum_{i\neq
j}^nE_j[\beta_{ij}(z_2){\mathbf s}_i^T\bbA_{ij}^{-1}(z_2)\underline{\bbA
}_j^{-1}(z_1)\underline{\bar{\mathbf s}}_j].
\end{equation}

The next aim is to replace $\beta_{ij}(z_2)$ in the equality above
by $\beta_{ij}^{\rmtr}(z_2)$. To this end, consider the case $i>j$
first. By (\ref{g9})
%
%e3.31 ###
%
\begin{equation}\label{g3}
E\bigl|E_j\bigl[\bigl(\beta_{ij}(z_2)-\beta_{ij}^{\rmtr}(z_2)\bigr){\mathbf s}_i^T\bbA
_{ij}^{-1}(z_2)\underline{\bbA}_j^{-1}(z_1)\underline{\bar{\mathbf s}}_j\bigr]\bigr|
=O(n^{-1/2}).
\end{equation}
Second, when $i<j$, break $\underline{\bbA}_j^{-1}(z_1)$ into the
sum of $\underline{\bbA}_{\uij}^{-1}(z_1)$ and
$\underline{\bbA}_{j}^{-1}(z_1)-\underline{\bbA}_{\uij}^{-1}(z_1)$,
$\underline{\bar{\mathbf s}}_j$ into the sum of $\underline{\bar
{\mathbf s}}_{\uij}$
and $\underline{\bar{\mathbf s}}_j-\underline{\bar{\mathbf s}}_{\uij}$, where
$\underline{\bbA}_{\uij}(z_1)=\underline{\bbA}_j(z_1)-n^{-1}{\mathbf
s}_i{\mathbf s}_i^T$
and $\underline{\bar{\mathbf s}}_{\uij}=\underline{\bar{\mathbf
s}}_j-{\mathbf s}_i/n$.
Then, when $i<j$, with notation
\[
\underline\beta_{\uij}(z)=\frac{1}{1+({1}/{n}){\mathbf s}_i^T\underline
{\bbA}_{\uij}^{-1}(z){\mathbf s}_i},
\]
we have
%
%e3.32 ###
%
\begin{equation}\label{g4}
E_j\bigl[\bigl(\beta_{ij}(z_2)-\beta_{ij}^{\rmtr}(z_2)\bigr){\mathbf s}_i^T\bbA
_{ij}^{-1}(z_2)\underline{\bbA}_j^{-1}(z_1)\underline{\bar{\mathbf s}}_j\bigr]
=c_{n1}+c_{n2}+c_{n3}+c_{n4},\hspace*{-28pt}
\end{equation}
where
\begin{eqnarray*}
c_{n1}&=&E_j\bigl[\bigl(\beta_{ij}(z_2)-\beta_{ij}^{\rmtr}(z_2)\bigr){\mathbf s}_i^T\bbA
_{ij}^{-1}(z_2)\underline{\bbA}_{\uij}^{-1}(z_1)\underline{\bar{\mathbf s}}_{\uij}\bigr],
\\
c_{n2}&=&\frac{1}{n}E_j\bigl[\bigl(\beta_{ij}(z_2)-\beta_{ij}^{\rmtr}(z_2)\bigr){\mathbf
s}_i^T\bbA_{ij}^{-1}(z_2)\underline{\bbA}_{\uij}^{-1}(z_1){\mathbf s}_i\bigr],
\\
c_{n3}&=&-\frac{1}{n}E_j\bigl[\bigl(\beta_{ij}(z_2)-\beta_{ij}^{\rmtr}(z_2)\bigr){\mathbf
s}_i^T\bbA_{ij}^{-1}(z_2)\underline{\bbA}_{\uij}^{-1}(z_1){\mathbf
s}_i{\mathbf s}_i^T\underline{\bbA}_{\uij}^{-1}(z_1)\underline\beta
_{\uij}(z_1)\underline{\bar{\mathbf s}}_{\uij}\bigr]
\end{eqnarray*}
and
\[
c_{n4}=-\frac{1}{n^2}E_j\bigl[\bigl(\beta_{ij}(z_2)-\beta_{ij}^{\rmtr}(z_2)\bigr){\mathbf
s}_i^T\bbA_{ij}^{-1}(z_2)\underline{\bbA}_{\uij}^{-1}(z_1){\mathbf
s}_i{\mathbf s}_i^T\underline{\bbA}_{\uij}^{-1}(z_1)\underline\beta
_{\uij}(z_1){\mathbf s}_{i}\bigr].
\]
It follows from (\ref{g9}) that $E|c_{nj}|\leq\mathfrak
{M}n^{-1/2},j=1,2,3,4$.
Thus, $\beta_{ij}(z_2)$ in (\ref{c5}) can be replaced by
$\beta_{ij}^{\rmtr}(z_2)$, as expected.

In what follows we use the notation $o_{L_1}(1)$ to denote
convergence to zero in~$L_1$. Moreover,\vspace*{1pt} note that
$E_j[\beta_{ij}^{\rmtr}(z_2){\mathbf s}_i^T\bbA_{ij}^{-1}(z_2)\underline
{\bbA}_j^{-1}(z_1)\underline{\bar{\mathbf s}}_j]=0$
when $i>j$. This, together with (\ref{g3}) and (\ref{g4}), implies
that
%
%e3.33 ###
%
\begin{eqnarray}
&&
E_j[\bar{\mathbf s}_j^T\bbA_{j}^{-1}(z_2)\underline{\bbA
}_j^{-1}(z_1)\underline{\bar{\mathbf s}}_j]\nonumber\\
&&\qquad=\frac{1}{n}\sum
_{i\neq j}^nE_j[\beta_{ij}^{\rmtr}(z_2){\mathbf s}_i^T\bbA_{ij}^{-1}(z_2)\underline
{\bbA}_j^{-1}(z_1)\underline{\bar{\mathbf
s}}_j]+o_{L_1}(1)\nonumber\\[-8pt]\\[-8pt]
&&\qquad=\frac{1}{n}\sum_{i<j}E_j[\beta_{ij}^{\rmtr}(z_2){\mathbf s}_i^T\bbA
_{ij}^{-1}(z_2)\underline{\bbA}_j^{-1}(z_1)\underline
{\bar{\mathbf s}}_j]+o_{L_1}(1)\nonumber\\
&&\qquad=d_{n1}+d_{n2}+d_{n3}+o_{L_1}(1),\nonumber
\end{eqnarray}
where
\begin{eqnarray*}
d_{n1}&=&\frac{1}{n^2}\sum_{i<j}E_j[\beta_{ij}^{\rmtr}(z_2){\mathbf
s}_i^T\bbA_{ij}^{-1}(z_2)\underline{\bbA}_{\uij}^{-1}(z_1){\mathbf
s}_i\underline{\beta}_{\uij}(z_1)],
\\
d_{n2}&=&\frac{1}{n}\sum_{i<j}E_j[\beta_{ij}^{\rmtr}(z_2){\mathbf
s}_i^T\bbA_{ij}^{-1}(z_2)\underline{\bbA}_{ij}^{-1}(z_1)\underline
{\bar{\mathbf s}}_{ij}]
\end{eqnarray*}
and
\[
d_{n3}=-\frac{1}{n^2}\sum_{i<j}E_j[\beta_{ij}^{\rmtr}(z_2){\mathbf
s}_i^T\bbA_{ij}^{-1}(z_2)\underline{\bbA}_{ij}^{-1}(z_1){\mathbf
s}_i{\mathbf s}_i^T\underline{\bbA}_{ij}^{-1}(z_1)\underline{\beta
}_{ij}(z_1)\underline{\bar{\mathbf s}}_{ij}].
\]
Here, in the last step, we apply
$\underline{\bar{\mathbf s}}_j={\mathbf s}_i/n+\underline{\bar{\mathbf
s}}_{ij}$ first,
then use (\ref{g1}) and finally split
$\underline{\bbA}_{j}^{-1}(z_1)$ into two parts as before.

We claim that the terms $d_{n2}$ and $d_{n3}$ are both negligible.
To see this, we first prove the following estimate:
%
%e3.34 ###
%
\begin{equation}\label{g6}
E\biggl|\frac{1}{n}\sum_{i<j}{\mathbf s}_i^T\bbA
_{ij}^{-1}(z_2)\underline{\bbA}_{ij}^{-1}(z_1)\underline{\bar{\mathbf
s}}_{ij}\biggr|^2=o(1).
\end{equation}
Indeed, the left-hand side of (\ref{g6}) may be expanded as
%
%e3.35 ###
%
\begin{equation}\label{g7}
\frac{1}{n^2}\sum_{i_1<j,
i_2<j}E({\mathbf s}_{i_1}^T\bbA_{i_1j}^{-1}(z_2)\underline{\bbA
}_{i_1j}^{-1}(z_1)\underline{\bar{\mathbf s}}_{i_1j}
{\mathbf s}_{i_2}^T\bbA_{i_2j}^{-1}(\bar{z}_2)\underline{\bbA
}_{i_2j}^{-1}(\bar{z}_1)\underline{\bar{\mathbf s}}_{i_2j}).
\end{equation}
From (\ref{g2}), the term corresponding to $i_1=i_2$ in (\ref{g7})
is bounded by
\[
\frac{1}{n^2}\sum_{i_1<j}E|{\mathbf s}_{i_1}^T\bbA
_{i_1j}^{-1}(z_2)\underline{\bbA}_{i_1j}^{-1}(z_1)\underline{\bar
{\mathbf s}}_{i_1j}|^2=O\biggl(\frac{1}{n}\biggr).
\]
To treat the case $i_1\neq i_2$, we need to further split
$\bbA_{i_1j}^{-1}(z_2)$ as the sum of $\bbA_{i_1i_2j}^{-1}(z_2)$
and $\bbA_{i_1j}^{-1}(z_2)-\bbA_{i_1i_2j}^{-1}(z_2)$, where
$\bbA_{i_1i_2j}(z_2)=\bbA_{i_1j}(z_2)-n^{-1}{\mathbf s}_{i_2}{\mathbf
s}_{i_2}^T$.
Moreover, both $\underline{\bbA}_{i_1j}^{-1}(z_1)$ and
$\underline{\bar{\mathbf s}}_{i_1j}$ are also needed to be similarly
split. To simplify notation, define
\begin{eqnarray*}
\beta_{i_1i_2j}(z)&=&\frac{1}{1+({1}/{n}){\mathbf s}_{i_2}^T\bbA
_{i_1i_2j}^{-1}(z){\mathbf s}_{i_2}},\\
\underline{\beta}_{i_1i_2j}(z)&=&\frac{1}{1+({1}/{n}){\mathbf
s}_{i_2}^T\underline{\bbA}_{i_1i_2j}^{-1}(z){\mathbf s}_{i_2}}
\end{eqnarray*}
and
\begin{eqnarray*}
\underline{\bbA}_{i_1i_2j}(z)&=&\underline{\bbA}_{i_1j}(z)-{\mathbf
s}_{i_2}{\mathbf s}_{i_2}^T,\qquad
\underline{\bar{\mathbf s}}_{i_1i_2j}=\underline{\bar{\mathbf
s}}_{i_1j}-\frac{{\mathbf s}_{i_2}}{n},\\
\zeta_{i_2j}&=&{\mathbf s}_{i_2}^T\bbA_{i_2j}^{-1}(\bar{z}_2)\underline
{\bbA}_{i_2j}^{-1}(\bar{z}_1)\underline{\bar{\mathbf s}}_{i_2j}.
\end{eqnarray*}
By (\ref{g1}), (\ref{g2}), (\ref{a4}) and (\ref{g41}) we have
\begin{eqnarray*}
\hspace*{-4pt}&&\frac{1}{n}|E({\mathbf s}_{i_1}^T\bbA
_{i_1j}^{-1}(z_2)\underline{\bbA}_{i_1j}^{-1}(z_1){\mathbf s}_{i_2}
\zeta_{i_2j})|\\
\hspace*{-4pt}&&\qquad=\frac{1}{n}|E({\mathbf s}_{i_1}^T\bbA
_{i_1j}^{-1}(z_2)\underline{\bbA}_{i_1i_2j}^{-1}(z_1)\underline
{\beta}_{i_1i_2j}(z_1){\mathbf s}_{i_2}
\zeta_{i_2j})|\\
\hspace*{-4pt}&&\qquad\leq
\frac{\mathfrak{M}}{n}|E({\mathbf s}_{i_1}^T\bbA
_{i_1i_2j}^{-1}(z_2)\underline{\bbA}_{i_1i_2j}^{-1}(z_1)\underline
{\beta}_{i_1i_2j}(z_1){\mathbf s}_{i_2}
\zeta_{i_2j})|\\
\hspace*{-4pt}&&\qquad\quad{}+\frac{\mathfrak{M}}{n^2}\bigl|E\bigl({\mathbf s}_{i_1}^T\bbA
_{i_1i_2j}^{-1}(z_2){\mathbf s}_{i_2}{\mathbf s}_{i_2}^T\bbA
_{i_1i_2j}^{-1}(z_2)\beta_{i_1i_2j}(z_2)\\
\hspace*{-4pt}&&\qquad\quad\hspace*{78.4pt}{}\times\underline{\bbA
}_{i_1i_2j}^{-1}(z_1)\underline{\beta}_{i_1i_2j}(z_1){\mathbf s}_{i_2}
\zeta_{i_2j}\bigr)\bigr|\\
\hspace*{-4pt}&&\qquad\leq\frac{\mathfrak{M}}{n}
(E|{\mathbf s}_{i_1}^T\bbA_{i_1i_2j}^{-1}(z_2)\underline{\bbA
}_{i_1i_2j}^{-1}(z_1){\mathbf s}_{i_2}|^2E|\zeta_{i_2j}|^2)^{1/2}\\
\hspace*{-4pt}&&\qquad\quad{}+\frac{\mathfrak{M}}{n^2}(E|{\mathbf s}_{i_1}^T\bbA
_{i_1i_2j}^{-1}(z_2){\mathbf s}_{i_2}{\mathbf s}_{i_2}^T\bbA
_{i_1i_2j}^{-1}(z_2)\underline{\bbA}_{i_1i_2j}^{-1}(z_1){\mathbf
s}_{i_2}|^2E|\zeta_{i_2j}|^2)^{1/2}\\
\hspace*{-4pt}&&\qquad=O(n^{-3/8});
\\
\hspace*{-4pt}&&\frac{1}{n}|E({\mathbf s}_{i_1}^T\bbA
_{i_1i_2j}^{-1}(z_2){\mathbf s}_{i_2}{\mathbf s}_{i_2}^T\beta
_{i_1i_2j}(z_2)\bbA_{i_1i_2j}^{-1}(z_2)\underline{\bbA
}_{i_1i_2j}^{-1}(z_1)\underline{\bar{\mathbf s}}_{i_1i_2j}
\zeta_{i_2j})|\\
\hspace*{-4pt}&&\qquad\leq
\frac{\mathfrak{M}}{n}(E|{\mathbf s}_{i_1}^T\bbA
_{i_1i_2j}^{-1}(z_2){\mathbf s}_{i_2}|^4E|{\mathbf s}_{i_2}^T\bbA
_{i_1i_2j}^{-1}(z_2)\underline{\bbA
}_{i_1i_2j}^{-1}(z_1)\underline{\bar{\mathbf
s}}_{i_1i_2j}|^4)^{1/4}(E|\zeta_{i_2j}|^2)^{1/2}\\
\hspace*{-4pt}&&\qquad=O(n^{-1/2});
\\
\hspace*{-4pt}&&\frac{1}{n}|E({\mathbf s}_{i_1}^T\bbA
_{i_1i_2j}^{-1}(z_2)\underline{\bbA}_{i_1i_2j}^{-1}(z_1){\mathbf
s}_{i_2}{\mathbf s}_{i_2}^T\underline{\beta}_{i_1i_2j}(z_1)\underline
{\bbA}_{i_1i_2j}^{-1}(z_1)\underline{\bar{\mathbf s}}_{i_1i_2j}
\zeta_{i_2j})|\\
\hspace*{-4pt}&&\qquad\leq
\frac{\mathfrak{M}}{n}(E|{\mathbf s}_{i_1}^T\bbA
_{i_1i_2j}^{-1}(z_2)\underline{\bbA}_{i_1i_2j}^{-1}(z_1){\mathbf
s}_{i_2}|^4E|{\mathbf s}_{i_2}^T\underline{\bbA
}_{i_1i_2j}^{-1}(z_1)\underline{\bar{\mathbf
s}}_{i_1i_2j}|^4)^{1/4}(E|\zeta_{i_2j}|^2)^{1/2}\\
\hspace*{-4pt}&&\qquad=O(n^{-1/2});
\\
\hspace*{-4pt}&&\frac{1}{n^2}|E({\mathbf s}_{i_1}^T\bbA
_{i_1i_2j}^{-1}(z_2){\mathbf s}_{i_2}{\mathbf s}_{i_2}^T\beta
_{i_1i_2j}(z_2)\bbA_{i_1i_2j}^{-1}(z_2)\\
\hspace*{-4pt}&&\hspace*{1.1pt}\quad{}\times
\underline{\bbA}_{i_1i_2j}^{-1}(z_1){\mathbf s}_{i_2}{\mathbf
s}_{i_2}^T\underline{\bbA}_{i_1i_2j}^{-1}(z_1)\underline{\beta
}_{i_1i_2j}(z_1)\underline{\bar{\mathbf s}}_{i_1i_2j}
\zeta_{i_2j})| \\
\hspace*{-4pt}&&\qquad\leq
\frac{\mathfrak{M}}{n^2}(E|{\mathbf s}_{i_1}^T\bbA
_{i_1i_2j}^{-1}(z_2){\mathbf s}_{i_2}{\mathbf s}_{i_2}^T\bbA
_{i_1i_2j}^{-1}(z_2)\underline{\bbA}_{i_1i_2j}^{-1}(z_1){\mathbf
s}_{i_2}|^4
E|{\mathbf s}_{i_2}^T\underline{\bbA
}_{i_1i_2j}^{-1}(z_1)\underline{\bar{\mathbf s}}_{i_1i_2j}|^4)^{1/4}\\
\hspace*{-4pt}&&\qquad\quad{}\times(E|\zeta_{i_2j}|^2)^{1/2} =O(n^{-3/8}).
\end{eqnarray*}
The above four estimates, together with
the fact that
\[
E\bigl({\mathbf s}_{i_1}^T\bbA_{i_1i_2j}^{-1}(z_2)\underline{\bbA
}_{i_1i_2j}^{-1}(z_1)
\times\underline{\bar{\mathbf s}}_{i_1i_2j} \zeta_{i_2j}
\bigr)=0,\qquad i_1\neq
i_2,
\]
imply that all terms in (\ref{g7}) corresponding to $i_1\neq i_2$
are bounded in absolute value by $\mathfrak{M}n^{-3/8}$, which ensures
(\ref{g6}).

Consider the term $d_{n2}$ now. In view of (\ref{g28}) and
(\ref{g9}) we may substitute $b_{12}(z_2)$ for
$\beta_{ij}^{\rmtr}(z_2)$ in the term $d_{n2}$ first and then apply
(\ref{g6}) to conclude that $E|d_{n2}|=o(1)$. As for the term
$d_{n3}$, it follows from (\ref{g28}) and (\ref{g9}) that
$\beta_{ij}^{\rmtr}(z_2),\underline{\beta}_{\uij}(z_1)$ and
${\mathbf s}_i^T\bbA_{ij}^{-1}(z_2)\underline{\bbA
}_{\uij}^{-1}(z_1){\mathbf s}_i$
can be replaced by $b_{12}(z_2),\underline{b}_{12}(z_1)$ and
$\frac{1}{n}\operatorname{tr}\bbA_{ij}^{-1}(z_2)\underline{\bbA}_{\uij}^{-1}(z_1)$,
respectively, where
\[
\underline{b}_{12}(z)=\frac{1}{1+({1}/{n})E\operatorname{tr}\underline{\bbA
}_{12}^{-1}(z)}
\]
[note: $\underline{b}_{12}(z)=b_{12}(z)$]. Moreover, by an
inequality similar to (\ref{g10}) we have
\begin{eqnarray*}
&&\biggl|E_j\biggl[{\mathbf s}_i^T\underline{\bbA}_{\uij}^{-1}(z_1)\underline{\bar
{\mathbf s}}_{\uij}\frac{1}{n}\bigl(\operatorname{tr}(\bbA_{ij}^{-1}(z_2)\underline{\bbA
}_{\uij}^{-1}(z_1))
-\operatorname{tr}(\bbA_{j}^{-1}(z_2)\underline{\bbA}_{j}^{-1}(z_1))\bigr)\biggr]\biggr|\\
&&\qquad\leq
\mathfrak{M}\frac{E_j|{\mathbf s}_i^T\underline{\bbA
}_{\uij}^{-1}(z_1)\underline{\bar{\mathbf s}}_{\uij}|}{n}.
\end{eqnarray*}
Therefore, from (\ref{g2}) we obtain
\[
d_{n3}=-\frac{b_{12}(z_2)\underline{b}_{12}(z_1)}{n^2}E_j\biggl[\operatorname{tr}(\bbA
_{j}^{-1}(z_2)\underline{\bbA}_{j}^{-1}(z_1))\sum_{i<j}{\mathbf
s}_i^T\underline{\bbA}_{\uij}^{-1}(z_1)\underline{\bar{\mathbf
s}}_{\uij}\biggr]+o_{L_1}(1).
\]
As in (\ref{g6}) we may prove that (even simpler)
%
%e3.36 ###
%
\begin{equation}\label{g30}
E\biggl|\frac{1}{n}\sum_{i<j}{\mathbf s}_i^T\underline{\bbA
}_{\uij}^{-1}(z_1)\underline{\bar{\mathbf s}}_{\uij}\biggr|^2=o(1),
\end{equation}
which then implies that $E|d_{n3}|=o(1)$.

As for $d_{n1}$, we conclude from (\ref{g28}), (\ref{g9}) and
(\ref{g10}) that
\begin{eqnarray*}
d_{n1}&=&\frac{b_{12}(z_2)b_{12}(z_1)}{n^2}\sum
_{i<j}\operatorname{tr}E_j[\bbA_{ij}^{-1}(z_2)\underline{\bbA
}_{\uij}^{-1}(z_1)]+o_{L_1}(1)\\
&=&\frac{b_{12}(z_2)b_{12}(z_1)}{n^2}(j-1)\operatorname{tr}[E_j(\bbA
_{j}^{-1}(z_2))E_j(\underline{\bbA}_{j}^{-1}(z_1))]+o_{L_1}(1).
\end{eqnarray*}
Summarizing the above, we have thus proved that
%
%e3.37 ###
%
\begin{eqnarray}\label{g13}\quad
&&E_j(\bar{\mathbf s}_j^T\bbA_j^{-1}(z_2))E_j(\bbA_j^{-1}(z_1)\bar
{\mathbf s}_j)
\nonumber\\
&&\qquad=\frac{j-1}{n^2}b_{12}(z_2)b_{12}(z_1)\operatorname{tr}[E_j(\bbA
_{j}^{-1}(z_2))E_j(\underline{\bbA
}_{j}^{-1}(z_1))]+o_{L_1}(1)\\
&&\qquad=\frac{j-1}{n^2}z_1z_2\underline{m}(z_1)\underline
{m}(z_2)\operatorname{tr}[E_j(\bbA_{j}^{-1}(z_2))E_j(\underline{\bbA
}_{j}^{-1}(z_1))]+o_{L_1}(1),\nonumber\vadjust{\goodbreak}
\end{eqnarray}
using the fact that, by (2.17) in \cite{b2} and (\ref{g10}),
%
%e3.38 ###
%
\begin{equation}\label{g47} b_{12}(z)\rightarrow
-z\underline{m}(z).
\end{equation}
This implies (\ref{c1}).

%s3.5 ###
\subsection{\texorpdfstring{The limit of (\protect\ref{a9})}{The limit of (3.3)}}

Our goal is to show that
%
%e3.39 ###
%
\begin{equation}
\label{c2} \mbox{(\ref{a9})}\stackrel{\mathrm{i.p.}}\longrightarrow0.
\end{equation}

In view of (\ref{g38}) we have
%
%e3.40 ###
%
\begin{equation}\label{g43}
(\ref{a9})=\frac{EX_{11}^3}{n}\sum_{j=1}^n\sum
_{i=1}^p[E_{j}(\bbD_j(z_2))]_{ii}
[E_{j}(\bbe_i^T\bbA_j^{-1}(z_1)\bar{\mathbf s}_j)].
\end{equation}
We first prove that $\bbe_i^T\bbA_j^{-1}(z_1)\bar{\mathbf s}_j$ above may
be replaced by $E(\bbe_i^T\bbA_j^{-1}(z_1)\bar{\mathbf s}_j)$. Using
martingale decompositions as in (\ref{a1}) and the fact that
$\bbe_i^T\bbA_{j}^{-1}(z)\bar{{\mathbf s}}_{j}=\bar{\mathbf s}_j^T\bbA
_j^{-1}(z)\bbe_i$,
we obtain that
%
%e3.41 ###
%
\begin{eqnarray}\label{c6}
&&\bar{\mathbf s}_j^T\bbA_j^{-1}(z_2)\bbe_iE_j[\theta_{ij}(z_1)]
\nonumber\\
&&\qquad=[\bar{\mathbf s}_j^T\bbA_j^{-1}(z_2)\bbe_i-E(\bar{\mathbf s}_j^T\bbA
_j^{-1}(z_2)\bbe_i)]E_j[\theta_{ij}(z_1)]\nonumber\\[-8pt]\\[-8pt]
&&\qquad\quad{}+E(\bar{\mathbf s}_j^T\bbA_j^{-1}(z_2)\bbe_i)E_j[\theta_{ij}(z_1)]
\nonumber\\
&&\qquad=\theta_{ij}(z_2) \times E_j[\theta_{ij}(z_1)]
+E(\bar{\mathbf s}_j^T\bbA_j^{-1}(z_2)\bbe_i)E_j[\theta
_{ij}(z_1)],
\nonumber
\end{eqnarray}
where
\[
\theta_{ij}(z)=\bbe_i^T\bbA_j^{-1}(z)\bar{{\mathbf s}}_j-E(\bbe
_i^T\bbA_j^{-1}(z)\bar{{\mathbf s}}_j)=\sum_{m\neq
j}^n(E_{m}-E_{m-1})(\theta_{ijm}(z))
\]
and
\begin{eqnarray*}
\theta_{ijm}(z)&=&\bbe_i^T\bbA_{j}^{-1}(z)\bar{{\mathbf s}}_{j}-\bbe
_i^T\bbA_{jm}^{-1}(z)\bar{{\mathbf s}}_{jm}\\
&=&\biggl[-\frac{1}{n^2}\bbe_i^T\bbA_{jm}^{-1}(z_1){\mathbf s}_{m}{\mathbf
s}_{m}^T\bbA_{jm}^{-1}(z)\beta_{mj}(z){\mathbf s}_{m}
\\
&&\hspace*{3.1pt}{}-\frac{1}{n}\bbe_i^T\bbA
_{jm}^{-1}(z){\mathbf s}_{m}{\mathbf s}_{m}^T\bbA_{jm}^{-1}(z)\beta
_{mj}(z)\bar{{\mathbf s}}_{jm}
+\frac{1}{n}\bbe_i^T\bbA_{jm}^{-1}(z){\mathbf s}_{m}\biggr].
\end{eqnarray*}

As in (\ref{g16}), one can verify that
%
%e3.42 ###
%
\begin{eqnarray}\label{g42}
E|n^{-1}\bbe_i^T\bbA_{jm}^{-1}(z){\mathbf s}_{m}|^k&=&O(n^{-k}),\qquad k=2
\mbox{ or } 4,\nonumber\\[-8pt]\\[-8pt]
E|n^{-1}\bbe_i^T\bbA_{jm}^{-1}(z){\mathbf
s}_{m}|^8&=&O(n^{-6}).\nonumber
\end{eqnarray}
Thus, for $k=2$ or $4$, via (\ref{g2}),
\[
E\biggl|\frac{1}{n}\bbe_i^T\bbA_{jm}^{-1}(z){\mathbf s}_{m}{\mathbf
s}_{m}^T\bbA
_{jm}^{-1}(z)\bar{{\mathbf s}}_{jm}\biggr|^k
=O(n^{-2}\varepsilon_n^{k-2})\vadjust{\goodbreak}
\]
and, via (\ref{a3}),
\[
E\biggl|\frac{1}{n^2}\bbe_i^T\bbA_{jm}^{-1}(z_1){\mathbf s}_{m}{\mathbf
s}_{m}^T\bbA_{jm}^{-1}(z){\mathbf s}_{m}\biggr|^k=O\bigl(n^{-2-{3(k-2)}/{4}}\bigr).
\]
These yield that $E|\theta_{ijm}(z)|^2=O(n^{-2})$, $E|\theta
_{ijm}(z)|^4=O(n^{-2}\varepsilon_n)$ and then
%
%e3.43 ###
%
\begin{eqnarray}
\label{g46}
E|\theta_{ij}(z)|^2&=&O(n^{-1}),\nonumber\\[-8pt]\\[-8pt]
E|\theta_{ij}(z)|^4&=&O(n^{-1}\varepsilon_n).\nonumber
\end{eqnarray}

Therefore,
%
%e3.44 ###
%
\begin{eqnarray}
\label{g17}
&& \Biggl[E\sum_{i=1}^p|[E_{j}(\bbD
_j(z_2))]_{ii}E_{j}(\theta_{ij}(z_1))|\Biggr]^2\nonumber\\
&&\qquad\leq
\sum_{i=1}^pE|\bbe_i^T\bbA_j^{-1}(z_2)\bar{\mathbf s}_j|^2\sum
_{i=1}^pE|\bar{\mathbf s}_j^T\bbA_j^{-1}(z_2)\bbe_iE_{j}(\theta
_{ij}(z_1))|^2
\nonumber\\
&&\qquad\leq \mathfrak{M} \sum_{i=1}^p[E|\theta_{ij}(z_2)|^4
E|\theta_{ij}(z_1)|^4]^{1/2}\\
&&\qquad\quad{} + \mathfrak{M}
\sum_{i=1}^p|E(\bar{\mathbf s}_j^T\bbA_j^{-1}(z_2)\bbe
_i)|^2E|\theta_{ij}(z_1)|^2
\nonumber\\
&&\qquad=O(\varepsilon_n).\nonumber
\end{eqnarray}
Here, by (\ref{a22})
\begin{eqnarray*}
\sum_{i=1}^p|E(\bar{\mathbf s}_j^T\bbA_j^{-1}(z_2)\bbe
_i)|^2E|\theta_{ij}(z_2)|^2&\leq&
\frac{\mathfrak{M}}{n}\sum_{i=1}^p|E(\bar{\mathbf s}_j^T\bbA
_j^{-1}(z_2)\bbe_i)|^2
\\
&\leq&
\frac{\mathfrak{M}}{n}E(\bar{\mathbf s}_j^T\bbA_j^{-1}(z_2)\bbA
_j^{-1}(\bar{z}_2)\bar{\mathbf s}_j)\\
&\leq&
\frac{\mathfrak{M}}{n}.
\end{eqnarray*}
Thus, $\bbe_i^T\bbA_j^{-1}(z_1)\bar{\mathbf s}_j$ involved in (\ref
{g43}) may be
replaced by $E(\bbe_i^T\bbA_j^{-1}(z_1)\bar{\mathbf s}_j)$, as
expected.

In addition, by (\ref{a22}) and (\ref{a58})
%
%e3.45 ###
%
\begin{eqnarray}\label{g18}
&&E\sum_{i=1}^p|[E_{j}(\bbD_j(z_2))]_{ii}
E(\bbe_i^T\bbA_j^{-1}(z_1)\bar{\mathbf s}_j)|\nonumber\\
&&\qquad\leq E\sum_{i=1}^p[E_{j}(\bbA_j^{-1}(\bar{z}_2)\bar{\mathbf s}_j\bar
{\mathbf s}_j^T\bbA_j^{-1}(z_2))]_{ii}
|E(\bbe_i^T\bbA_j^{-1}(z_1)\bar{\mathbf s}_j)|\\
&&\qquad \leq
{\max_{i}}|E(\bbe_i^T\bbA_1^{-1}(z_1)\bar{\mathbf s}_1)|E(\bar
{\mathbf s}_j^T\bbA_j^{-1}(z_2)\bbA_j^{-1}(\bar{z}_2)\bar{\mathbf
s}_j)\rightarrow0.\nonumber
\end{eqnarray}

It follows from (\ref{g17}) and (\ref{g18}) that
%
%e3.46 ###
%
\begin{equation}
\label{g19}E\sum_{i=1}^p|[E_{j}(\bbD_j(z_2))]_{ii}
E_j(\bbe_i^T\bbA_j^{-1}(z_1)\bar{\mathbf s}_j)|\rightarrow0,
\end{equation}
which then ensures (\ref{c2}).

%s3.6 ###
\subsection{\texorpdfstring{The limit of (\protect\ref{a10})}{The limit of (3.4)}}

The goal in this section is to prove that
%
%e3.47 ###
%
\begin{eqnarray}\label{c3}
\mbox{(\ref{a10})}&=&\frac{2z_1^2z_2^2\underline{m}^2(z_1)\underline
{m}^2(z_2)}{n}\nonumber\\[-8pt]\\[-8pt]
&&{}\times\sum_{j=1}^n\frac{(j-1)^2}{n^4}[\operatorname{tr}
(E_j(\bbA_{j}^{-1}(z_2))E_j(\bbA_{j}^{-1}(z_1)))]^2+o_p(1).\nonumber
\end{eqnarray}

First, (\ref{i1}) shows that (\ref{a10}) is equal to
%
%e3.48 ###
%
\begin{eqnarray}\label{a53}
&&\frac{E|X_{11}|^4-3}{n}\sum_{j=1}^n\sum
_{i=1}^pE_j(\bbD_j(z_1))_{ii}E_j(\bbD_j(z_2))_{ii}\nonumber\\[-8pt]\\[-8pt]
&&\qquad{}+\frac{2}{n}\sum_{j=1}^n
\operatorname{tr}[E_j(\bbD_j(z_1))E_j(\bbD_j(z_2))].\nonumber
\end{eqnarray}
To prove (\ref{c3}), the strategy is to substitute
$E_j(\bar{\mathbf s}_j^T\bbA_j^{-1}(z))$ for each
$\bar{\mathbf s}_j^T\bbA_j^{-1}(z)$ involved in $E_j(\bbD_j(z))$
by a martingale method.
As we shall see, the above first term converges to zero in
probability and the second term has a close connection with
(\ref{a8}).

Consider the second term of (\ref{a53}) first. Write
%
%e3.49 ###
%
\begin{eqnarray}\label{g12}
&&\operatorname{tr}[E_j(\bbD_j(z_1))E_j(\bbD_j(z_2))]\nonumber\\
&&\qquad=E_j[\bar{\mathbf s}_j^T\bbA_j^{-1}(z_1)\underline{\bbA
}_j^{-1}(z_2)\underline{\bar{\mathbf s}}_j\underline{\bar{\mathbf
s}}_j^T\underline{\bbA}_j^{-1}(z_2)\bbA_j^{-1}(z_1)\bar{\mathbf
s}_j]\\
&&\qquad=E_j[\bar{\mathbf s}_j^T\bbA_j^{-1}(z_1)\underline{\bbA
}_j^{-1}(z_2)\underline{\bar{\mathbf s}}_j\underline{\bar{\mathbf
s}}_j^T\underline{\bbA}_j^{-1}(z_2)E_j(\bbA_j^{-1}(z_1)\bar
{\mathbf s}_j)]
+f_{n},\nonumber
\end{eqnarray}
where
\[
f_n=E_j\bigl[\bar{\mathbf s}_j^T\bbA_j^{-1}(z_1)\underline{\bbA
}_j^{-1}(z_2)\underline{\bar{\mathbf s}}_j\underline{\bar{\mathbf
s}}_j^T\underline{\bbA}_j^{-1}(z_2)\bigl(\bbA_j^{-1}(z_1)\bar{\mathbf
s}_j-E_j(\bbA_j^{-1}(z_1)\bar{\mathbf s}_j)\bigr)\bigr].
\]

We claim that
%
%e3.50 ###
%
\begin{equation}\label{g14}
E|f_n|=o(1).\vadjust{\goodbreak}
\end{equation}
To see this, let
$\underline{E}_{\uij}=E(\cdot|{\mathbf s}_1,\ldots,{\mathbf
s}_i,\underline
{{\mathbf s}}_{j+1},
\ldots,\underline{{\mathbf s}}_n)$. Then, recalling the definitions of
$\underline{\bbA}_j^{-1}(z)$ and $\underline{\bar{\mathbf s}}_j$ as
before, we obtain a martingale decomposition
\begin{eqnarray*}
&&\underline{\bar{\mathbf s}}_j^T\underline{\bbA}_j^{-1}(z_2)
\bigl(\bbA_j^{-1}(z_1)\bar{\mathbf s}_j-E_{jj}(\bbA_j^{-1}(z_1)\bar{\mathbf
s}_j)\bigr)\\
&&\qquad=\sum_{i=j+1}^n\bigl(\underline{E}_{\uij}[\underline
{\bar{\mathbf s}}_j^T\underline{\bbA}_j^{-1}(z_2)\bbA_j^{-1}(z_1)\bar
{\mathbf s}_j]-
\underline{E}_{(i-1)j}[\underline{\bar{\mathbf s}}_j^T\underline
{\bbA}_j^{-1}(z_2)\bbA_j^{-1}(z_1)\bar{\mathbf s}_j]\bigr)\\
&&\qquad=\sum_{i=j+1}^n\bigl(\underline{E}_{\uij}-\underline
{E}_{(i-1)j}\bigr)[\underline{\bar{\mathbf s}}_j^T\underline{\bbA
}_j^{-1}(z_2)\bbA_j^{-1}(z_1)\bar{\mathbf s}_j-\underline{\bar{\mathbf
s}}_j^T\underline{\bbA}_j^{-1}(z_2)\bbA_{ij}^{-1}(z_1)\bar{\mathbf
s}_{ij}]\\
&&\qquad=f_{n1}+f_{n2},
\end{eqnarray*}
where
\[
f_{n1}=\frac{1}{n}\sum_{i=j+1}^n\bigl(\underline
{E}_{\uij}-\underline{E}_{(i-1)j}\bigr)[\underline{\bar{\mathbf
s}}_j^T\underline{\bbA}_j^{-1}(z_2)\bbA_{ij}^{-1}(z_1){\mathbf s}_i\beta
_{ij}(z_1)]
\]
and
\[
f_{n2}=-\frac{1}{n}\sum_{i=j+1}^n\bigl(\underline
{E}_{\uij}-\underline{E}_{(i-1)j}\bigr)[\underline{\bar{\mathbf
s}}_j^T\underline{\bbA}_j^{-1}(z_2)\bbA_{ij}^{-1}(z_1){\mathbf
s}_i{\mathbf s}_i^T\bbA_{ij}^{-1}(z_1)\bar{{\mathbf s}}_{ij}\beta_{ij}(z_1)].
\]
Note that $\underline{\bar{\mathbf s}}_j$ is independent of ${\mathbf
s}_i$ for
$i>j$. Then applying (\ref{g2}) yields
\begin{eqnarray*}
E|f_{n1}|^2\leq
\frac{\mathfrak{M}}{n^2}\sum_{i=j+1}^nE|\underline{\bar
{\mathbf s}}_j^T\underline{\bbA}_j^{-1}(z_2)\bbA_{ij}^{-1}(z_1){\mathbf
s}_i|^2
=O\biggl(\frac{1}{n}\biggr)
\end{eqnarray*}
and
\begin{eqnarray*}
E|f_{n2}|^2&\leq&
\frac{\mathfrak{M}}{n^2}\sum_{i=j+1}^nE|\underline{\bar
{\mathbf s}}_j^T\underline{\bbA}_j^{-1}(z_2)\bbA_{ij}^{-1}(z_1){\mathbf
s}_i{\mathbf s}_i^T\bbA_{ij}^{-1}(z_1)\bar{{\mathbf s}}_{ij}|^2\\
&\leq&\frac{\mathfrak{M}}{n^2}\sum_{i=j+1}^n\bigl(E|\underline
{\bar{\mathbf s}}_j^T\underline{\bbA}_j^{-1}(z_2)\bbA
_{ij}^{-1}(z_1){\mathbf s}_i|^4E|{\mathbf s}_i^T\bbA_{ij}^{-1}(z_1)\bar
{{\mathbf s}}_{ij}|^4\bigr)^{1/2}\\
&=&O\biggl(\frac{1}{n}\biggr),
\end{eqnarray*}
which ensures that
\[
E\bigl|\underline{\bar{\mathbf s}}_j^T\underline{\bbA}_j^{-1}(z_2)
\bigl(\bbA_j^{-1}(z_1)\bar{\mathbf s}_j-E_{jj}(\bbA_j^{-1}(z_1)\bar{\mathbf
s}_j)\bigr)\bigr|^2=O\biggl(\frac{1}{n}\biggr).\vadjust{\goodbreak}
\]
So (\ref{g14}) follows from the above estimate and
\[
E|\bar{\mathbf s}_j^T\bbA_j^{-1}(z_1)\underline{\bbA
}_j^{-1}(z_2)\underline{\bar{\mathbf s}}_j|^2=O(1),
\]
which may be obtained immediately by checking the argument of
(\ref{a22}).

As in (\ref{g14}) we may also prove that
%
%e3.51 ###
%
\begin{eqnarray}\label{g15}\quad
&&E\bigl|E_j\bigl[\bar{\mathbf s}_j^T\bbA_j^{-1}(z_1)\underline{\bbA
}_j^{-1}(z_2)\nonumber\\[-8pt]\\[-8pt]
&&\qquad\hspace*{3.1pt}{}\times\underline{\bar{\mathbf s}}_j\bigl(\underline{\bar{\mathbf
s}}_j^T\underline{\bbA}_j^{-1}(z_2)
-
E_j(\underline{\bar{\mathbf s}}_j^T\underline{\bbA}_j^{-1}(z_2))
\bigr)E_j(\bbA_j^{-1}(z_1)\bar{\mathbf s}_j)\bigr]\bigr|
=o(1).\nonumber
\end{eqnarray}
Therefore, combining (\ref{g12})--(\ref{g15}) with
(\ref{g13}) we have
%
%e3.52 ###
%
\begin{eqnarray}\label{g31}\quad
&&\operatorname{tr}[E_j(\bbD_j(z_1))E_j(\bbD_j(z_2))]\nonumber\\
&&\qquad= E_j[\bar{\mathbf s}_j^T\bbA_j^{-1}(z_1)\underline{\bbA
}_j^{-1}(z_2)\underline{\bar{\mathbf s}}_jE_j(\underline{\bar{\mathbf
s}}_j^T\underline{\bbA}_j^{-1}(z_2))E_j(\bbA_j^{-1}(z_1)\bar{\mathbf
s}_j)]\nonumber\\
&&\qquad\quad{}+o_{L_1}(1)\nonumber\\
&&\qquad=E_j(\bar{\mathbf s}_j^T\bbA_j^{-1}(z_1))E_j(\bbA_j^{-1}(z_2)\bar
{\mathbf s}_j)E_j(\bar{\mathbf s}_j^T\bbA_j^{-1}(z_2))E_j(\bbA
_j^{-1}(z_1)\bar
{\mathbf s}_j)\\
&&\qquad\quad{} +o_{L_1}(1)\nonumber\\
&&\qquad=\frac{(j-1)^2}{n^4}z_1^2z_2^2\underline{m}^2(z_1)\underline
{m}^2(z_2)[\operatorname{tr}(E_j(\bbA_{j}^{-1}(z_2))E_j\underline{\bbA
}_{j}^{-1}(z_1))]^2\nonumber\\
&&\qquad\quad{}+o_{L_1}(1).\nonumber
\end{eqnarray}

We now turn to the first term in (\ref{a53}) and claim that
%
%e3.53 ###
%
\begin{equation}\label{a60}
\frac{1}{n}\sum_{j=1}^n\sum_{i=1}^pE_j(\bbD
_j(z_1))_{ii}E_j(\bbD_j(z_2))_{ii}\stackrel{\mathrm{i.p.}}\longrightarrow0.
\end{equation}
Indeed, it follows from (\ref{c6}) that
\begin{eqnarray*}
&&E\Biggl|\sum_{i=1}^pE_j(\bbD_j(z_2))_{ii}E_j(\theta
_{ij}(z_1)\bar{\mathbf s}_j^T\bbA_{j}^{-1}(z_1)\bbe_i)\Biggr|
\\
&&\qquad\leq\sum_{i=1}^pE|E_j(\bbD_j(z_2))_{ii}E_j(\theta_{ij}(z_1))^2|\\
&&\qquad\quad{}+\sum_{i=1}^pE|E_j(\bbD_j(z_2))_{ii}E_j(\theta
_{ij}(z_1))E(\bar{\mathbf s}_j^T\bbA_{j}^{-1}(z_1)\bbe_i)|.
\end{eqnarray*}
The second term above is not greater than
\[
{\max_{i}}|E(\bar{\mathbf s}_j^T\bbA_{j}^{-1}(z_1)\bbe_i)|\sum
_{i=1}^pE|E_j(\bbD_j(z_2))_{ii}E_j(\theta_{ij}(z_1))|,
\]
which converges to zero by (\ref{g17}) and (\ref{a58}). Moreover,
by (\ref{a22}) and (\ref{g46})
\begin{eqnarray*}
&&\Biggl(\sum_{i=1}^pE|E_j(\bbD_j(z_2))_{ii}E_j(\theta
_{ij}(z_1))^2|\Biggr)^2\\
&&\qquad\leq
\sum_{i=1}^p
E|(\bbD_j(z_2))_{ii}|^2\sum_{i=1}^pE|\theta_{ij}(z_1)|^4
\\
&&\qquad\leq
E\Biggl(\sum_{i=1}^p\bar{{\mathbf s}}_j^T\bbA_j^{-1}(\bar{z}_1)\bbe
_i\bbe_i^T\bbA_j^{-1}(z_1)\bar{{\mathbf s}}_j\Biggr)^2\sum
_{i=1}^pE|\theta_{ij}(z_1)|^4\\
&&\qquad=O(\varepsilon_n).
\end{eqnarray*}
In addition, it follows from Lemma \ref{lem2} and (\ref{g19}) that
\begin{eqnarray*}
&&
E\Biggl|\sum_{i=1}^pE_j(\bbD_j(z_2))_{ii}E(\bbe_i^T\bbA
_j^{-1}(z_1)\bar{{\mathbf s}}_j)E_j(\bar{\mathbf s}_j^T\bbA
_{j}^{-1}(z_1)\bbe_i)\Biggr|
\\
&&\qquad
\leq{\max_{i}}|E(\bbe_i^T\bbA_1^{-1}(z_1)\bar{{\mathbf s}}_1)|\sum
_{i=1}^pE|E_j(\bbD_j(z_2))_{ii}E_j(\bar{\mathbf s}_j^T\bbA
_{j}^{-1}(z_1)\bbe_i)|
\rightarrow0.
\end{eqnarray*}
Consequently, the proof of (\ref{a60}) is complete. Thus,
(\ref{c3}) follows from (\ref{g31}), (\ref{a60}) and (\ref{a53}).

%s3.7 ###
\subsection{\texorpdfstring{The limit of (\protect\ref{i3})}{The limit of (3.28)}}

Note that (see \cite{b2}, (2.18))
%
%e3.54 ###
%
\begin{eqnarray}\label{d13}\qquad
&&\operatorname{tr}(E_j(\bbA_{j}^{-1}(z_2))E_j(\bbA_{j}^{-1}(z_1))
)\biggl[1-\frac{(j-1)p}{n^2}\frac{\underline{m}_n(z_1)\underline
{m}_n(z_2)}{(1+\underline{m}_n(z_1))(1+\underline
{m}_n(z_2))}\biggr]\nonumber\\[-8pt]\\[-8pt]
&&\qquad=\frac{p}{z_1z_2(1+\underline{m}_n(z_1))(1+\underline{m}_n(z_2))}+l_n,\nonumber
\end{eqnarray}
where $E|l_n|\leq\mathfrak{M}\sqrt{n}$ and $\underline{m}_n(z)$ is defined
like $m_n(z)$, but corresponding to $\underline{m}(z)$. Obviously,
$\underline{m}_n(z)\rightarrow\underline{m}(z)$. This implies that
\begin{eqnarray*}
&&\frac{(j-1)z_1z_2\underline{m}(z_1)\underline{m}(z_2)}{n^2}\operatorname{tr}
(E_j(\bbA_{j}^{-1}(z_2))E_j(\bbA_{j}^{-1}(z_1)))
\\
&&\qquad=\frac{z_1z_2(1+\underline{m}(z_1))(1+\underline{m}(z_2))}{p}\operatorname{tr}
(E_j(\bbA_{j}^{-1}(z_2))E_j(\bbA_{j}^{-1}(z_1)))\\
&&\qquad\quad{}-1+o_{L_1}(1),
\end{eqnarray*}
which, together with (\ref{g13}) and (\ref{g31}), leads to
\begin{eqnarray*}
&&
4\operatorname{tr}[E_{j}(\bbA_j^{-1}(z_1)\bar{\mathbf s}_j)E_{j}(\bar{\mathbf
s}_j^T\bbA_j^{-1}(z_2))]+2\operatorname{tr}(E_j(\bbD_j(z_1))E_j(\bbD
_j(z_2)))
\\
&&\qquad=\frac{4(j-1)z_1z_2\underline{m}(z_1)\underline{m}(z_2)}{n^2}\operatorname{tr}
(E_j(\bbA_{j}^{-1}(z_2))E_j(\underline{\bbA}_{j}^{-1}(z_1)))
\\
&&\qquad\quad{}
+\frac{2(j-1)^2z_1^2z_2^2\underline{m}^2(z_1)\underline
{m}^2(z_2)}{n^4}[\operatorname{tr}(E_j(\bbA_{j}^{-1}(z_2))E_j(\underline
{\bbA}_{j}^{-1}(z_1)))]^2\\
&&\qquad\quad{}+o_{L_1}(1)
\\
&&\qquad
=-2+2z_1^2z_2^2\bigl(1+\underline{m}(z_1)\bigr)^2\bigl(1+\underline{m}(z_2)\bigr)^2\frac
{[\operatorname{tr}(E_j(\bbA_{j}^{-1}(z_2))E_j(\underline{\bbA
}_{j}^{-1}(z_1)))]^2}{p^2}\\
&&\qquad\quad{}+o_{L_1}(1).
\end{eqnarray*}

Further, we conclude from (\ref{d13}) that
\begin{eqnarray*}\label{d20}
&&\frac{1}{np^2}\sum_{j=1}^n[\operatorname{tr}(E_j(\bbA
_{j}^{-1}(z_1))E_j(\bbA_{j}^{-1}(z_2)))]^2\nonumber\\
&&\qquad=\frac{1}{z_1^2z_2^2(1+\underline{m}(z_1))^2(1+\underline
{m}(z_2))^2}\\
&&\qquad\quad{}\times\frac{1}{n}\sum_{j=1}^n\frac{1}{(1-
{(j-1)p}/{n^2}({\underline{m}(z_1)\underline
{m}(z_2)}/({(1+\underline{m}(z_1))(1+\underline{m}(z_2))})))^2}\\
&&\qquad\quad{} + o_p(1)
\nonumber\\
&&\qquad\stackrel{\mathrm{i.p.}}\longrightarrow
\frac{1}{z_1^2z_2^2(1+\underline{m}(z_1))^2(1+\underline
{m}(z_2))^2}\\
&&\qquad\quad\hspace*{7.7pt}{}\times\int^1_0\frac{dx}{(1-x({c\underline
{m}(z_1)\underline{m}(z_2)}/({(1+\underline{m}(z_1))(1+\underline
{m}(z_2))})))^2}\nonumber\\
&&\qquad=
\frac{1}{z_1^2z_2^2(1+\underline{m}(z_1))(1+\underline
{m}(z_2))[(1+\underline{m}(z_1))(1+\underline{m}(z_2))-c\underline
{m}(z_1)\underline{m}(z_2)]}.
\end{eqnarray*}
It follows that
%
%e3.55 ###
%
\begin{eqnarray}\label{f8}\quad
\mbox{(\ref{i3})}&=&z_1z_2\underline{m}(z_1)\underline{m}(z_2)\frac
{1}{n}\sum_{j=1}^n\bigl[4\operatorname{tr}[E_{j}(\bbA_j^{-1}(z_1)\bar{\mathbf
s}_j)E_{j}(\bar{\mathbf s}_j^T\bbA_j^{-1}(z_2))]
\nonumber\\
&&\hspace*{121.2pt}{}+2\operatorname{tr}[E_j(\bbD_j(z_1))E_j(\bbD
_j(z_2))]\bigr]\nonumber\\[-8pt]\\[-8pt]
&&{}+o_p(1)
\nonumber\\
&\stackrel{\mathrm{i.p.}}\longrightarrow&
\frac{2cz_1z_2\underline{m}^2(z_1)\underline
{m}^2(z_2)}{(1+\underline{m}(z_1))(1+\underline{m}(z_2))-c\underline
{m}(z_1)\underline{m}(z_2)}.\nonumber
\end{eqnarray}

%s4 ###
\section{\texorpdfstring{Tightness of ${\hat{M}}_n^{(1)}(z)$ and convergence
of $M_n^{(2)}(z)$}{Tightness of ${M}_n^{(1)}(z)$ and convergence
of $M_n^{(2)}(z)$}}\label{tightness}

First, we proceed to prove the tightness of $\hat{M}_n^{(1)}(z)$ for
$z\in\mathcal{C}$, which is a truncated version of $M_n(z)$ as in
(\ref{g20}). By (\ref{g2}) we have
\[
E\Biggl|\sum_{i=1}^ma_i\sum_{j=1}^nY_j(z_i)\Biggr|^2=\sum
_{j=1}^nE\Biggl|\sum_{i=1}^ma_iY_j(z_i)\Biggr|^2\leq
\mathfrak{M},\qquad v_0=\Im z_i,
\]
which ensures that condition (i) of Theorem 12.3 in \cite{bili2}
is satisfied, as pointed out in \cite{b2}. Here $Y_j(z)$ is
defined in (\ref{a7}). Condition (ii) of Theorem 12.3 in \cite{bili2}
will be
verified if the following holds:
%
%e4.1 ###
%
\begin{equation}\label{c7}
E\frac{|M_n^{(1)}(z_1)-M_n^{(1)}(z_2)|^2}{|z_1-z_2|^2}\leq\mathfrak{M}\qquad
\mbox{for } z_1,z_2\in\mathcal{C}_n^+\cup\mathcal{C}_n^-.
\end{equation}
In the sequel, since $\mathcal{C}_n^+$ and $\mathcal{C}_n^-$ are
symmetric, we shall prove the above inequality on $\mathcal{C}_n^+$
only. Throughout this section, all bounds including $O(\cdot)$ and
$o(\cdot)$ expressions hold uniformly for $z\in\mathcal{C}_n^+$.

In view of our truncation steps, (1.9a) and (1.9b) in \cite{b2}
apply to our case as well, that is, for any
$\eta_1>(1+\sqrt{c})^2$, $0<\eta_2<I(0,1)(c)(1-\sqrt{c})^2$ and
any positive $l$
%
%e4.2 ###
%
\begin{equation}\label{g25}
P(\|\bbS\|\geq\eta_1)=o(n^{-l}),\qquad
P\bigl(\lambda_{\min}(\bbS)\leq\eta_2\bigr)=o(n^{-l}).
\end{equation}
Note that when either $z\in\mathcal{C}_u$ or
$z\in\mathcal{C}_l$ and $u_l<0$, $\|\bbA_j^{-1}(z)\|$ is bounded in $n$.
But this is not the case for $z\in\mathcal{C}_r$ or
$z\in\mathcal{C}_l$ and $u_l>0$. In general, for $z\in\mathcal{C}_n^+$,
we have
%
%e4.3 ###
%
\begin{equation}\label{g22}
\|\bbA_j^{-1}(z)\|\leq M+ v^{-1}I\bigl(\|\bbA_j\|\geq h_r \mbox{ or }
\lambda_{\min}(\bbA_j)\leq h_l \bigr).
\end{equation}
Here, $\bbA_j=\bbS-{\mathbf s}_j{\mathbf s}_j^T$,
$h_r\in((1+\sqrt{c})^2,u_r)$ and $h_l\in(u_l,(1-\sqrt{c})^2)$.

Note that $\bbA^{-1}(z_1)-\bbA^{-1}(z_2)=(z_2-z_1)\bbA^{-1}(z_1)\bbA
^{-1}(z_2)$. As in Section \ref{simplif}, we then write
%
%e4.4 ###
%
\begin{eqnarray}\label{a62}
&&\frac{M_n^{(1)}(z_1)-M_n^{(1)}(z_2)}{z_1-z_2}\nonumber\\
&&\qquad=-\sqrt{n}\sum_{j=1}^n(E_j-E_{j-1})[\bar{\mathbf s}^T\bbA
^{-1}(z_1)\bbA^{-1}(z_2)\bar{\mathbf s}\\
&&\qquad\quad\hspace*{84.7pt}\hspace*{14.4pt}{}-\bar{\mathbf s}_j^T\bbA
_j^{-1}(z_1)\bbA_j^{-1}(z_2)\bar{\mathbf s}_j].\nonumber
\end{eqnarray}
Moreover, expanding the above difference we get
\[
\bar{\mathbf s}^T\bbA^{-1}(z_1)\bbA^{-1}(z_2)\bar{\mathbf s}-\bar
{\mathbf s}_j^T\bbA_j^{-1}(z_1)\bbA_j^{-1}(z_2)\bar{\mathbf s}_j
= q_{n1}+q_{n2}+q_{n3},
\]
where
\begin{eqnarray*}
q_{n1}&=&(\bar{\mathbf s}^T-\bar{\mathbf s}_j^T)\bbA^{-1}(z_1)\bbA
^{-1}(z_2)\bar{\mathbf s},\\
q_{n2}&=&\bar{\mathbf s}_j^T\bigl(\bbA^{-1}(z_1)\bbA^{-1}(z_2)-\bbA
_j^{-1}(z_1)\bbA_j^{-1}(z_2)\bigr)\bar{\mathbf s}
\end{eqnarray*}
and
\[
q_{n3}=\bar{\mathbf s}_j^T\bbA_j^{-1}(z_1)\bbA_j^{-1}(z_2)(\bar{\mathbf
s}-\bar{\mathbf s}_j).
\]

It follows from (\ref{g21}), (\ref{g2}), (\ref{g22}) and (\ref{g25})
that
\begin{eqnarray*}
E\Biggl|\sqrt{n}\sum_{j=1}^n(E_j-E_{j-1})q_{n3}\Biggr|^2&\leq&
\frac{1}{n}\sum_{j=1}^nE|\bar{\mathbf s}_j^T\bbA
_j^{-1}(z_1)\bbA_j^{-1}(z_2){\mathbf s}_j|^2
\\
&\leq&\mathfrak{M}+ \mathfrak{M}n^{8}\rho_n^{-4}P\bigl(\|\bbA_1\|\geq
h_r \mbox{ or }
\lambda_{\min}(\bbA_1)\leq h_l\bigr)\\
&\leq&\mathfrak{M},
\end{eqnarray*}
where we use, on the event $(\|\bbA_j\|\geq h_r \mbox{ or }
\lambda_{\min}(\bbA_j)\leq h_l)$, by (\ref{g8}),
%
%e4.5 ###
%
\begin{eqnarray}\label{g26}
|\bar{\mathbf s}_j^T\bbA_j^{-1}(z_1)\bbA_j^{-1}(z_2){\mathbf s}_j|&\leq&
\|\bar{\mathbf s}_j\|\|{\mathbf s}_j\|\|\bbA_j^{-1}(z_1)\bbA
_j^{-1}(z_2)\|\nonumber\\[-8pt]\\[-8pt]
&\leq&\mathfrak{M}
v^{-2}n^{2}\leq\mathfrak{M}n^{4}\rho_n^{-2}.\nonumber
\end{eqnarray}

For $q_{n2}$, expanding its difference term by term we
have
\[
q_{n2}=q_{n2}^{(1)}+\cdots+q_{n2}^{(6)},
\]
where
\begin{eqnarray*}
q_{n2}^{(1)}&=&\frac{1}{n^2}\bar{\mathbf s}_j^T\beta_j(z_1)\beta
_j(z_2)\tbbA_j(z_1)\tbbA_j(z_2)\bar{\mathbf s}_j,\\
q_{n2}^{(2)}&=&-\frac{1}{n}\bar{\mathbf s}_j^T\beta_j(z_1)\tbbA
_j(z_1)\bbA_j^{-1}(z_2)\bar{\mathbf s}_j,
\\
q_{n2}^{(3)}&=&-\frac{1}{n}\bar{\mathbf s}_j^T\beta_j(z_2)\bbA
_j^{-1}(z_1)\tbbA_j(z_2)\bar{\mathbf s}_j,\\
q_{n2}^{(4)}&=&\frac{1}{n^3}\bar{\mathbf s}_j^T\beta_j(z_1)\beta
_j(z_2)\tbbA_j(z_1)\tbbA_j(z_2){\mathbf s}_j
\end{eqnarray*}
and
\begin{eqnarray*}
q_{n2}^{(5)}&=&-\frac{1}{n^2}\bar{\mathbf s}_j^T\beta_j(z_1)\tbbA
_j(z_1)\bbA_j^{-1}(z_2){\mathbf s}_j,\\
q_{n2}^{(6)}&=&-\frac{1}{n^2}\bar{\mathbf s}_j^T\beta_j(z_2)\bbA
_j^{-1}(z_1)\tbbA_j(z_2){\mathbf s}_j.
\end{eqnarray*}
We conclude from (\ref{g9}), (\ref{g25}), (\ref{g22}) and
(\ref{g26}) that
\[
E\Biggl|\sqrt{n}\sum_{j=1}^n(E_j-E_{j-1})q_{n2}^{(6)}\Biggr|^2\leq
\mathfrak{M}+\mathfrak{M}v^{-8}n^{8}P\bigl(\|\bbS\|\geq h_r \mbox{ or }
\lambda_{\min}(\bbA_1)\leq h_l\bigr)\leq\mathfrak{M},
\]
where we use, on the event $(\|\bbS\|\geq h_r \mbox{ or }
\lambda_{\min}(\bbA_1)\leq h_l)$,
%
%e4.6 ###
%
\begin{equation}\label{g23}
|\beta_j(z)|=|1-n^{-1}{\mathbf s}_j^T\bbA^{-1}(z){\mathbf s}_j|\leq
1+n^{-1}v^{-1}\|{\mathbf s}_j\|^2\leq\mathfrak{M}v^{-1}n
\end{equation}
by (\ref{i2}). Similar argument shows that
\[
E\Biggl|\sqrt{n}\sum
_{j=1}^n(E_j-E_{j-1})q_{n2}^{(6)}\Biggr|^2=O(1),\qquad j=2,\ldots,5.
\]
Moreover, write $q_{n1}=q_{n1}^{(1)}+q_{n1}^{(2)}+q_{n1}^{(3)}$,
where
\begin{eqnarray*}
q_{n1}^{(1)}&=&\frac{1}{n^2}\beta_j(z_1)\beta_j(z_2){\mathbf s}_j^T\bbA
_j^{-1}(z_1)\bbA_j^{-1}(z_2){\mathbf s}_j,\\
q_{n1}^{(2)}&=&\frac{1}{n}\beta_j(z_1){\mathbf s}_j^T\bbA_j^{-1}(z_1)\bbA
_j^{-1}(z_2)\bar{\mathbf s}_{j}
\end{eqnarray*}
and
\[
q_{n1}^{(3)}=-\frac{1}{n^2}\beta_j(z_1)\beta_j(z_2){\mathbf s}_j^T\bbA
_j^{-1}(z_1)\tbbA_j(z_2)\bar{\mathbf s}_{j}.
\]
The argument for $q_{n2}^{(6)}$ also works for
$q_{n1}^{(j)},j=1,2,3$, and thus,
\[
E\Biggl|\sqrt{n}\sum_{j=1}^n(E_j-E_{j-1})q_{n1}\Biggr|^2\leq\mathfrak{M}.
\]
The proof of (\ref{c7}) is complete.

Next, consider $M_n^{(2)}(z)$. By
$\bar{\mathbf s}=n^{-1}\sum_{i=1}^n{\mathbf s}_i$, (\ref{g1}) and an
equality similar to (\ref{a11}) we obtain
\begin{eqnarray*}
\sqrt{n}E(\bar{\mathbf s}^T\bbA^{-1}(z)\bar{\mathbf s})&=&\frac{1}{\sqrt
{n}}\sum_{i=1}^nE(\beta_i(z){\mathbf s}_i^T\bbA_i^{-1}(z)\bar
{\mathbf s})\\
&=&\frac{1}{\sqrt{n}}\sum_{i=1}^nE(\beta_i(z){\mathbf s}_i^T\bbA
_i^{-1}(z)\bar{\mathbf s}_i)\\
&&{}+\frac{1}{n^{3/2}}\sum
_{i=1}^nE(\beta_i(z){\mathbf s}_i^T\bbA_i^{-1}(z){\mathbf s}_i)\\
&=&\frac{b_1(z)}{\sqrt{n}}E(\operatorname{tr}\bbA_1^{-1}(z))+b_1(z)t_{n1}+b_1(z)t_{n2},
\end{eqnarray*}
where
\begin{eqnarray*}
t_{n1}&=&-\sqrt{n}E(\beta_1(z)\xi_1(z){\mathbf s}_1^T\bbA_1^{-1}(z)\bar
{\mathbf s}_1),\\
t_{n2}&=&-\frac{1}{\sqrt{n}}E(\beta_1(z)\xi_1(z){\mathbf s}_1^T\bbA
_1^{-1}(z){\mathbf s}_1).
\end{eqnarray*}
Again, using an equality similar to (\ref{a11}) further gives
\[
t_{n1}= t_{n1}^{(1)}+ t_{n1}^{(2)},\qquad t_{n2}= t_{n2}^{(1)}+
t_{n2}^{(2)},
\]
where
\begin{eqnarray*}
t_{n1}^{(1)}&=&-\sqrt{n}b_1(z)E(\xi_1(z){\mathbf s}_1^T\bbA_1^{-1}(z)\bar
{\mathbf s}_1),\\
t_{n1}^{(2)}&=&\sqrt{n}b_1(z)E(\beta_1(z)\xi
_1^2(z){\mathbf s}_1^T\bbA_1^{-1}(z)\bar{\mathbf s}_1)
\end{eqnarray*}
and
\begin{eqnarray*}
t_{n2}^{(1)}&=&-\frac{b_1(z)}{\sqrt{n}}E(\xi_1(z){\mathbf s}_1^T\bbA
_1^{-1}(z){\mathbf s}_1),\\
t_{n2}^{(2)}&=&\frac{b_1(z)}{\sqrt{n}}E(\beta_1(z)\xi_1^2(z){\mathbf
s}_1^T\bbA_1^{-1}(z){\mathbf s}_1).
\end{eqnarray*}

Note that $|b_1(z)|\leq\mathfrak{M}$ for $z\in\mathcal{C}_n$ (see
\cite{b2}, three lines below (3.6)). It follows from (\ref{g28}),
(\ref{a3}),
(\ref{g2}), (\ref{g25}), (\ref{g22}) and (\ref{g23}) that
\[
\bigl|t_{n1}^{(2)}\bigr|\leq\mathfrak{M}\varepsilon_n+ \mathfrak{M}n^{10}\rho
_n^{-4}P\bigl(\|\bbS\|\geq
h_r \mbox{ or } \lambda_{\min}(\bbA_1)\leq h_l\bigr)\leq\mathfrak
{M}\varepsilon_n,
\]
because
$|\beta_i(z)\xi_i^2(z){\mathbf s}_i^T\bbA_i^{-1}(z)\bar{\mathbf
s}_i|\leq
n^{5}v^{-4}$ on the event $(\|\bbS\|\geq h_r \mbox{ or }
\lambda_{\min}(\bbA_1)\leq h_l)$. This argument clearly applies
to $t_{n2}^{(2)}$ as well and so $|t_{n2}^{(2)}|\leq
\mathfrak{M}\varepsilon_n$. Notice that $\frac{1}{n}E[\operatorname{tr}
(\bbA_1^{-1}(z))]=E[\bbe_m^T\bbA_1^{-1}(z)\bbe_m]$. This and
(\ref{g38}) show that
\begin{eqnarray*}
\bigl|t_{n1}^{(1)}\bigr|&=&\Biggl|-\frac{b_1(z)EX_{11}^3}{\sqrt{n}}\sum
_{m=1}^pE(\bbe_m^T\bbA_1^{-1}(z)\bbe_m\bbe_m^T\bbA_1^{-1}(z)\bar
{\mathbf s}_1)\Biggr|\\
&=&\Biggl|-\frac{b_1(z)EX_{11}^3({1}/{n})E\operatorname{tr} (\bbA_1^{-1}(z))}{\sqrt
{n}}\sum_{m=1}^pE(\bbe_m^T\bbA_1^{-1}(z)\bar{\mathbf s}_1)\Biggr|+o(1)\\
&\leq&\mathfrak{M}\biggl|b_1(z)EX_{11}^3\frac{1}{n}E\operatorname{tr} (\bbA
_1^{-1}(z))\biggr|\max_{m}\sqrt{n}|E(\bbe_m^T\bbA_1^{-1}(z)\bar
{\mathbf s}_1)|+o(1)\\
&=&o(1),
\end{eqnarray*}
where we make use of the facts that by (\ref{a58}), (\ref{g25}) and
(\ref{g22}),
\[
\max_{m}\sqrt{n}|E(\bbe_m^T\bbA_1^{-1}(z)\bar{\mathbf
s}_1)|=o(1),\qquad
\frac{1}{n}E\operatorname{tr} (\bbA_1^{-1}(z))=O(1)
\]
and that by (\ref{g29}), (\ref{g46}), (\ref{g25}) and (\ref{g22}),
\begin{eqnarray*}
&&E\bigl|\bigl(\bbe_m^T\bbA_1^{-1}(z)\bbe_m-E(\bbe_m^T\bbA_1^{-1}(z)\bbe
_m)\bigr)\bigl(\bbe_m^T\bbA_1^{-1}(z)\bar{\mathbf s}_1-E(\bbe_m^T\bbA
_1^{-1}(z)\bar{\mathbf s}_1)\bigr)\bigr|\\
&&\qquad\leq\mathfrak{M}n^{-1}+\mathfrak{M}v^{-2}nP\bigl(\|\bbA_1\|\geq h_r
\mbox{ or }
\lambda_{\min}(\bbA_1)\leq h_l\bigr)=O(n^{-1}).
\end{eqnarray*}
Note that
$n^{-1}E(\xi_i(z){\mathbf s}_i^T\bbA_i^{-1}(z){\mathbf s}_i)=E\gamma
_i^2(z)+n^{-2}E(\operatorname{tr}\bbA_i^{-1}(z)-E\operatorname{tr}\bbA_i^{-1}(z))^2$
and then applying (\ref{g28}), (\ref{a3}), (\ref{g25}) and
(\ref{g22}) gives $t_{n2}^{(1)}=O(n^{-1/2})$.

Summarizing the above we obtain
\[
\sqrt{n}E(\bar{\mathbf s}^T\bbA^{-1}(z)\bar{\mathbf s})=\frac
{b_1(z)}{n^{1/2}}E(\operatorname{tr}\bbA_1^{-1}(z))+o(1).
\]
Moreover, it is proven in \cite{b2}, Section 4, that
$n(E\operatorname{tr}\bbA^{-1}(z)/n-c_nm_n(z))$ is bounded for $z\in\mathcal{C}_n$.
In addition, by (\ref{g10}), (\ref{g25}) and (\ref{g22}) we have
\[
\sqrt{n}\biggl|\frac{E\operatorname{tr}\bbA_1^{-1}(z)}{n}-\frac{E\operatorname{tr}\bbA
^{-1}(z)}{n}\biggr|\leq\frac{\mathfrak{M}}{\sqrt{n}}.
\]
It follows that $n(E\operatorname{tr}\bbA_1^{-1}(z)/n-c_nm_n(z))$ is bounded. This,
together with the boundedness of $b_1(z)$, shows that
\[
\sup_{z\in\mathcal{C}_n}\sqrt{n}\biggl(E\bar{\mathbf s}^T\bbA
^{-1}(z)\bar{\mathbf s}-\frac{c_nm_n(z)}{1+c_nm_n(z)}\biggr)\rightarrow
0.
\]

%s5 ###
\section{\texorpdfstring{Proofs of Lemma \protect\ref{theo3}, Theorems \protect
\ref{th1} and \protect\ref{th2}}{Proofs of Lemma 1, Theorems 1 and 2}}
\label{finalproof}

\mbox{}

\begin{pf*}{Proof of Lemma \ref{theo3}} To finish Lemma
\ref{theo3}, $\bar{\mathbf s}^T\bar{\mathbf s}-c_n$ needs to be written
as a sum
of martingale difference sequence so that we can get a CLT for $\bar
{\mathbf s}^T\bar{\mathbf s}-c_n$ and, more importantly, obtain
the asymptotic covariance between $\bar{\mathbf s}^T\bar{\mathbf
s}-c_n$ and
$\bar{\mathbf s}^T\bbA^{-1}(z)\bar{\mathbf s}$.

Thus, write
%
%e5.1 ###
%
\begin{eqnarray}\label{f2}
\sqrt{n}(\bar{\mathbf s}^T\bar{\mathbf s}-c_n)
&=&\sqrt{n}\sum_{j=1}^n(E_j-E_{j-1})(\bar{\mathbf s}^T\bar{\mathbf
s})\nonumber\\
&=&\sqrt{n}\sum_{j=1}^n(E_j-E_{j-1})(\bar{\mathbf s}^T\bar{\mathbf
s}-\bar{\mathbf s}_j^T\bar{\mathbf s}_j)
\nonumber\\[-8pt]\\[-8pt]
&=&\sqrt{n}\sum_{j=1}^n(E_j-E_{j-1})\biggl(2\frac{\bar{\mathbf s}_j^T{\mathbf
s}_j}{n}+\frac{{\mathbf s}_j^T{\mathbf s}_j}{n^2}\biggr)\nonumber\\
&=&\frac{2}{\sqrt{n}}\sum_{j=1}^nE_j(\bar{\mathbf s}_j^T{\mathbf
s}_j)+o_p(1),\nonumber
\end{eqnarray}
because
\[
E\Biggl|\sqrt{n}\sum_{j=1}^n(E_j-E_{j-1})\biggl(\frac{{\mathbf s}_j^T{\mathbf
s}_j}{n^2}\biggr)\Biggr|^2=\frac{1}{n^3}\sum_{j=1}^nE|E_j({\mathbf s}_j^T{\mathbf
s}_j)-p|^2=O\biggl(\frac{1}{n}\biggr).
\]

From (\ref{g2}) we have
\begin{eqnarray*}
\sum_{j=1}^nE\biggl|\frac{1}{\sqrt{n}}E_j(\bar{\mathbf s}_j^T{\mathbf
s}_j)\biggr|^2I\biggl(\frac{1}{\sqrt{n}}E_j(\bar{\mathbf s}_j^T{\mathbf s}_j)\geq
\varepsilon\biggr)
&\leq&\frac{1}{\varepsilon^2}\sum_{j=1}^nE\biggl|\frac{1}{\sqrt
{n}}E_j(\bar{\mathbf s}_j^T{\mathbf s}_j)\biggr|^4\\
&=&O(n^{-1}),
\end{eqnarray*}
which implies condition (ii) of Lemma \ref{lemma3}. Look at
condition (i) of Lemma \ref{lemma3} next. It is easily seen that
\[
E_{j-1}[E_j(\bar{\mathbf s}_j^T{\mathbf s}_j)]^2=E_{j}(\bar{\mathbf
s}_j^T)E_{j}(\bar{\mathbf s}_j)=
\frac{1}{n^2}\sum_{k_1<j,k_2<j}{\mathbf s}_{k_1}^T{\mathbf s}_{k_2}.
\]
Furthermore, for the term corresponding to $k_1=k_2$, we have
\[
E\biggl|\frac{1}{n^2}\sum_{k_1<j}[{\mathbf s}_{k_1}^T{\mathbf
s}_{k_1}-E({\mathbf s}_{k_1}^T{\mathbf s}_{k_1})]\biggr|^2=
\frac{1}{n^4}\sum_{k_1<j}E|{\mathbf s}_{k_1}^T{\mathbf
s}_{k_1}-E({\mathbf s}_{k_1}^T{\mathbf s}_{k_1})|^2=O\biggl(\frac{1}{n^2}\biggr).
\]
On the other hand, when $k_1\neq k_2$,
\begin{eqnarray*}
E\biggl|\frac{1}{n^2}\sum_{k_1\neq
k_2}{\mathbf s}_{k_1}^T{\mathbf s}_{k_2}\biggr|^2&=&\frac{1}{n^4}\sum
_{k_1\neq
k_2,h_1\neq
h_2}E[{\mathbf s}_{k_1}^T{\mathbf s}_{k_2}{\mathbf s}_{h_1}^T{\mathbf
s}_{h_2}]\\
&=&\frac
{2}{n^4}\sum_{k_1\neq
k_2}E({\mathbf s}_{k_1}^T{\mathbf s}_{k_2})^2=O\biggl(\frac{1}{n}\biggr).
\end{eqnarray*}
It follows that
%
%e5.2 ###
%
\begin{equation}\label{f7}
\frac{4}{n}\sum_{j=1}^nE_{j-1}[E_j(\bar{\mathbf s}_j^T{\mathbf
s}_j)]^2=\frac{4}{n}\sum_{j=1}^n\frac
{c(j-1)}{n}+o_p(1)\stackrel{\mathrm{i.p.}}\longrightarrow
4c\int^1_0x\,dx=2c.\hspace*{-28pt}
\end{equation}
Therefore, by Lemma \ref{lemma3}
%
%e5.3 ###
%
\begin{equation}\label{f5}
\sqrt{n}(\bar{\mathbf s}^T\bar{\mathbf s}-c_n)\stackrel
{D}\longrightarrow N(0,2c).
\end{equation}

We conclude from Sections \ref{sim} and \ref{weakcon} that $\hat
{M}_n(z)$ converges
weakly to a Gaussian process on $\mathcal{C}$. Moreover,
$m_n(z)\rightarrow m(z)$ uniformly on $\mathcal{C}$ by (4.2) in
\cite{b2} and (\ref{d22}). These, together with (\ref{d21}),
(\ref{f5}), (\ref{a7}) and (\ref{f2}), give, for any constants $a_1$
and~$a_2$,
%
%e5.4 ###
%
\begin{eqnarray}%\label{f9}
\label{e4}
&&
a_1X_n(z)+a_2\sqrt{n}\bigl(g(\|\bar{\mathbf s}\|^2)-g(c_n)\bigr)\nonumber\\
&&\qquad=\tilde{a}_1(z)\sqrt{n}\biggl[\bar{\mathbf s}^T\bbA^{-1}(z)\bar
{\mathbf s}-\frac{c_nm_n(z)}{1+c_nm_n(z)}\biggr]\nonumber\\[-8pt]\\[-8pt]
&&\qquad\quad{}+\tilde{a}_2(z) \sqrt{n}(\|\bar{\mathbf s}\|^2-c_n)+o_p(1)\nonumber\\
&&\qquad=\sum_{j=1}^nl_j(z) +o_p(1)\nonumber,
\end{eqnarray}
where $\tilde{a}_1(z)=a_1(1+cm(z))^2/c,
\tilde{a}_2(z)=a_2g'(c_n)-a_1m(z)/c$ and
\[
l_j(z)=\tilde{a}_1(z)Y_j(z)+\tilde{a}_2(z)\frac{2}{\sqrt
{n}}E_j(\bar{\mathbf s}_j^T{\mathbf s}_j).
\]
Here, the first $o_p(1)$ denotes convergence in probability to zero
in the $C$ space and in the first step we use the fact that $
g(x)=g(c_n)+g'(a)(x-c_n)+o(|x-c_n|)$
as $x\rightarrow c_n$.
Thus, tightness of $\hat{X}_n(z)$ is from that of $\hat{M}_n(z)$.

Since $b_{1}(z)\rightarrow1/(1+cm(z))$ and
$b_{1}(z)\rightarrow-z\underline{m}(z) $ by (2.17) in \cite{b2}, we
have
%
%e5.5 ###
%
\begin{equation}
\label{g33}1/\bigl(1+cm(z)\bigr)=-z\underline{m}(z).
\end{equation}
Moreover, we assume for the moment that
%
%e5.6 ###
%
\begin{equation}
\label{g32}
\sum_{j=1}^nE_{j-1}\biggl[Y_j(z)\frac{2}{\sqrt{n}}E_j(\bar{\mathbf
s}_j^T{\mathbf s}_j)\biggr]\stackrel{\mathrm{i.p.}}\longrightarrow
\frac{2cm(z)}{(1+cm(z))^2}.
\end{equation}
It follows from (\ref{f8}), (\ref{f7}), (\ref{g32}) and (\ref{g33})
that
\begin{eqnarray*}
&&\sum_{j=1}^nE_{j-1}[l_j(z_1)l_j(z_2)]\\
&&\qquad=\tilde{a}_1(z_1)\tilde
{a}_1(z_2)\frac{2cz_1z_2\underline{m}^2(z_1)\underline
{m}^2(z_2)}{(1+\underline{m}(z_1))(1+\underline{m}(z_2))-c\underline
{m}(z_1)\underline{m}(z_2)}
\\
&&\qquad\quad{}
+2c\tilde{a}_2(z_1)\tilde{a}_2(z_2)+\tilde{a}_1(z_1)\tilde
{a}_2(z_2)\frac{2cm(z_1)}{(1+cm(z_1))^2}\\
&&\qquad\quad{}+\tilde{a}_1(z_2)\tilde
{a}_2(z_1)\frac{2cm(z_2)}{(1+cm(z_2))^2}
+o_p(1)
\\
&&\qquad=a_1^2\times\mbox{(\ref{e33})}+a_2^2 \times2c(g'(c))^2+o_p(1).
\end{eqnarray*}
Thus, Lemma \ref{theo3} follows from the above argument, Lemma
\ref{lemma3} and Cram\'er--Wold's device.

Now consider (\ref{g32}). Write
\begin{eqnarray*}
&&E_{j-1}[E_j({\mathbf s}_j^T\bbA_j^{-1}(z)\bar{\mathbf s}_j)E_j(\bar
{\mathbf s}_j^T{\mathbf s}_j)]\\
&&\qquad=E_j(\bar{\mathbf s}_j^T)E_j(\bbA
_j^{-1}(z)\bar{\mathbf s}_j)
=\frac{1}{n}\sum_{i<j}E_j({\mathbf s}_i^T\bbA_{ij}^{-1}(z)\bar
{\mathbf s}_j\beta_{ij}(z))
\\
&&\qquad=
\frac{1}{n^2}\sum_{i<j}E_j({\mathbf s}_i^T\bbA_{ij}^{-1}(z){\mathbf
s}_i\beta_{ij}(z))
+\frac{1}{n}\sum_{i<j}E_j({\mathbf s}_i^T\bbA
_{ij}^{-1}(z)\bar{\mathbf s}_{ij}\beta_{ij}(z)),
\end{eqnarray*}
where we use $\bar{\mathbf s}_j=1/n\sum_{i\neq j}{\mathbf s}_i$ in the
second step and $\bar{\mathbf s}_j=\bar{\mathbf s}_{ij}+{\mathbf
s}_i/n$ in the last
step. By (\ref{g28}), (\ref{g2}) and (\ref{a3})
\[
E\biggl|\frac{1}{n}\sum_{i<j}E_j\bigl({\mathbf s}_i^T\bbA_{ij}^{-1}(z)\bar
{\mathbf s}_{ij}(\beta_{ij}(z))-b_{12}(z)\bigr)\biggr|=O\biggl(\frac{1}{\sqrt{n}}\biggr),
\]
which, together with (\ref{g30}), yields
\[
E\biggl|\frac{1}{n}\sum_{i<j}E_j({\mathbf s}_i^T\bbA_{ij}^{-1}(z)\bar
{\mathbf s}_{ij}\beta_{ij}(z))\biggr|=o(1).
\]
On the other hand, appealing to (\ref{g10}), (\ref{g28}) and
(\ref{a3}) ensures that
\[
\frac{1}{n^2}\sum_{i<j}E_j({\mathbf s}_i^T\bbA_{ij}^{-1}(z){\mathbf
s}_i\beta_{ij}(z))=\frac{j-1}{n}\frac{n^{-1}E\operatorname{tr}\bbA
^{-1}(z)}{1+n^{-1}E\operatorname{tr}\bbA^{-1}(z)}+o_{L_1}(1).
\]
Therefore, we obtain
%
%e5.7 ###
%
\begin{eqnarray}\label{d24}
&&\frac{1}{n}\sum_{j=1}^nE_{j-1}[E_j({\mathbf s}_j^T\bbA
_j^{-1}(z)\bar{\mathbf s}_j)E_j(\bar{\mathbf s}_j^T{\mathbf
s}_j)]\nonumber\\
&&\qquad=\frac{n^{-1}E\operatorname{tr}\bbA^{-1}(z)}{1+n^{-1}E\operatorname{tr}\bbA^{-1}(z)}\frac
{1}{n}\sum_{j=1}^n\frac{j-1}{n}+o_{L_1}(1)
\\
&&\qquad\stackrel{\mathrm{i.p.}}\longrightarrow\frac{cm(z)}{2(1+cm(z))}.\nonumber
\end{eqnarray}

Next, by the Markov inequality and the Doob inequality
\begin{eqnarray*}
P\biggl(\max_{
i,j}\frac{1}{n}\biggl|\sum_{k<j}v_{ik}\biggr|\geq\varepsilon\biggr) &\leq&
\frac{\sum_{i=1}^nE(\max_{j}({1}/{n})|{\sum
_{k<j}v_{1k}}|)^4}{\varepsilon^4}
\\
&\leq&
\frac{\mathfrak{M}nE(({1}/{n})|{\sum
_{k<j}v_{ik}}|)^4}{\varepsilon^4}
\leq\frac{\mathfrak{M}}{n},
\end{eqnarray*}
which implies
\[
\max_{i,j}\biggl|\frac{1}{n}\sum_{k<j}v_{ik}\biggr|\stackrel
{\mathrm{i.p.}}\longrightarrow
0.
\]
This and (\ref{g38}) ensure that
%
%e5.8 ###
%
\begin{eqnarray}\label{d23}
&&
\sum_{j=1}^nE_{j-1}[E_j\alpha_j(z)E_j(\bar{\mathbf s}_j^T{\mathbf
s}_j)]\nonumber\\
&&\qquad=\frac{EX_{11}^3}{n}\sum_{j=1}^n\sum
_{i=1}^p[E_{j}\bbD_j(z_2)]_{ii}
[E_{j}(\bbe_i^T\bar{\mathbf s}_j)]\nonumber\\
&&\qquad\leq\max_{
i,j}\biggl|\frac{1}{n}\sum_{k<j}v_{ik}\biggr|\frac{\mathfrak{M}}{n}\sum
_{j=1}^n\sum_{i=1}^p[E_{j}(\bbA_j^{-1}(z_2)\bar{\mathbf s}_j\bar
{\mathbf s}_j^T\bbA_j^{-1}(\bar{z}_2))]_{ii}\\
&&\qquad\leq\max_{
i,j}\biggl|\frac{1}{n}\sum_{k<j}v_{ik}\biggr|\frac{\mathfrak{M}}{n}\sum
_{j=1}^nE_{j}(\bar{\mathbf s}_j^T\bbA_j^{-1}(\bar{z}_2)\bbA
_j^{-1}(z_2)\bar{\mathbf s}_j)
\nonumber\\
&&\qquad\stackrel{\mathrm{i.p.}}\longrightarrow0,\nonumber
\end{eqnarray}
because\vspace*{1pt}
(\ref{a22}) implies that
$n^{-1}\sum_{j=1}^nE_{j}(\bar{\mathbf s}_j^T\bbA_j^{-1}(\bar
{z}_2)\bbA_j^{-1}(z_2)\bar{\mathbf s}_j)$
is uniformly integrable. Based on (\ref{d23}) and (\ref{d24}) we
have (\ref{g32}).
\end{pf*}
\begin{pf*}{Proof of Remark \ref{re1}} By (\ref{e1}) we get
%
%e5.9 ###
%
\begin{equation}\label{g34}
\frac{\underline{m}(z_1)-\underline{m}(z_2)}{(z_1-z_2)}=
\frac{\underline{m}(z_1)\underline{m}(z_2)(1+\underline
{m}(z_1))(1+\underline{m}(z_2))}{(1+\underline{m}(z_1))(1+\underline
{m}(z_2))-c\underline{m}(z_1)\underline{m}(z_2)}.
\end{equation}
Then
\begin{eqnarray*}
&&\frac{2}{cz_1z_2[(1+\underline{m}(z_1))(1+\underline
{m}(z_2))-c\underline{m}(z_1)\underline{m}(z_2)]}
\\
&&\qquad=\frac{2(\underline{m}(z_1)-\underline
{m}(z_2))}{cz_1z_2(z_1-z_2)\underline{m}(z_1)\underline
{m}(z_2)(1+\underline{m}(z_1))(1+\underline{m}(z_2))}
\\
&&\qquad
=\frac{2(\underline{m}(z_1)-\underline
{m}(z_2))}{z_1z_2(z_1-z_2)(1+\underline{m}(z_1))^2(1+\underline{m}(z_2))^2}
\\
&&\qquad\quad{}+\frac{2(\underline{m}(z_1)-\underline
{m}(z_2))}{cz_1z_2(z_1-z_2)(1+\underline{m}(z_1))(1+\underline{m}(z_2))}\\
&&\qquad\quad\hspace*{10.8pt}{}\times\biggl[\frac{1}{\underline{m}(z_1)\underline{m}(z_2)}-\frac
{c}{(1+\underline{m}(z_1))(1+\underline{m}(z_2))}\biggr]
\\
&&\qquad=\frac{2(\underline{m}(z_1)-\underline
{m}(z_2))}{z_1z_2(z_1-z_2)(1+\underline{m}(z_1))^2(1+\underline{m}(z_2))^2}\\
&&\qquad\quad{}+\frac{2}{cz_1z_2(1+\underline{m}(z_1))(1+\underline{m}(z_2))}
\\
&&\qquad=\frac{2(\underline{m}(z_1)-\underline
{m}(z_2))}{z_1z_2(z_1-z_2)(1+\underline{m}(z_1))^2(1+\underline{m}(z_2))^2}\\
&&\qquad\quad{}+\frac{2m(z_1)m(z_2)}{c},
\end{eqnarray*}
where in the first step and the third step we use (\ref{g34}) and in
the last step we use (\ref{g33}). On the other hand, via (\ref{e1})
one can verify that
\[
\frac{2(\underline{m}(z_1)-\underline
{m}(z_2))}{z_1z_2(z_1-z_2)(1+\underline{m}(z_1))^2(1+\underline{m}(z_2))^2}
=\frac{2(z_2\underline{m}(z_2)-z_1\underline
{m}(z_1))^2}{c^2z_1z_2(z_1-z_2)(\underline{m}(z_1)-\underline{m}(z_2))},
\]
which is exactly the covariance function in Lemma 2 of \cite{b1}.
Therefore, Remark \ref{re1} holds.
\end{pf*}
\begin{pf*}{Proof of Theorem \ref{th2}} The idea from Lemma
\ref{theo3} to Theorem \ref{th2} is similar to that in \cite{b2}.
First, by the Cauchy formula we have
\[
\int f(x)\,dG(x)=-\frac{1}{2\pi i}\oint f(z)m_G(z)\,d(z),
\]
where the contour contains the support of $G(x)$ on which $f(x)$ is
analytic. Then, with probability one, we have
\[
\int f(x)\,dG_n(x)=-\frac{1}{2\pi i}\oint f(z)X_n(z)\,d(z)
\]
for all $n$ large, where the complex integral is over $\mathcal{C}$
and
\[
G_n(x)=\sqrt{n}\bigl(F_2^{\bbS}(x)-F_{c_n}(x)\bigr).
\]
Further,
\[
\biggl|\int
f(z)\bigl(X_n(z)-\hat{X}_n(z)\bigr)\,dz\biggr|\leq\frac{\mathfrak{M}\rho_n}{\sqrt
{n}(u_r-\lambda_{\max}(\bolds{
\mathcal{S}}))}+\frac{\mathfrak{M}\rho_n}{\sqrt{n}(\lambda_{\min
}(\bolds{
\mathcal{S}})-u_l)} \stackrel{\mathrm{a.s.}}\longrightarrow0,
\]
where, with probability one, $\lambda_{\max}
(\bolds{\mathcal{S}})\rightarrow(1+\sqrt{c})^2$ by \cite{j} and
$\lambda_{\min}(\bolds{\mathcal{S}})\rightarrow(1-\sqrt{c})^2$ by
\cite{z}.
Second, note
that for any constants $a_1$ and $a_2$
\[
(\hat{X}_n(z),Y_n)\rightarrow a_1\oint f(z)\hat{X}_n(z)\,dz+ a_2Y_n
\]
is a continuous mapping. Therefore, the right-hand side above converges
in distribution by Lemma \ref{theo3}. Moreover, Remark \ref{re1}
shows that (\ref{e3}) follows from (1.12) and (1.15) in \cite{b1}.
\end{pf*}
\begin{pf*}{Proof of Theorem \ref{th1}}
By taking $f(x)=x^{-1}$ and $g(x)=x$ in Theorem \ref{th2} and noting
that $c_n\to c$
as $n\to\infty$, we can complete the proof.
\end{pf*}

\vspace*{-10pt}

\begin{appendix}\label{app}
%s6 ###
\section*{Appendix}

%s6.1 ###
\subsection{Some lemmas} We collect some results needed to prove Lemma
\ref{theo3}.

\begin{lemma}[(Burkholder \cite{burk})]\label{lem3}
Let $\{Y_i\}$ be a complex martingale difference sequence with
respect to the increasing $\sigma$-field $\{\mathcal{F}_i\}$. Then
for $k\geq2$
\[
E\biggl|\sum_{i}Y_i\biggr|^k\leq
\mathfrak{M}_kE\biggl(\sum_{i}E(|Y_i|^2|\mathcal
{F}_{i-1})\biggr)^{k/2}+\mathfrak{M}_kE\biggl(\sum
_{i}|Y_i|^k\biggr).
\]
\end{lemma}
\begin{lemma}[(Theorem 35.12 of Billingsley \cite{bili})]\label{lemma3}
Suppose for each n, $Y_{n,1},Y_{n,2},\break\ldots, Y_{n,r_n}$ is a real
martingale difference sequence with respect to the increasing
$\sigma$-field $\{{\mathcal{F}}_{n,j}\}$ having second moments. If as
$n\rightarrow\infty$
\begin{eqnarray*}
\mbox{\textup{(i)}\hspace*{6.5pt}}\quad
\sum_{j=1}^{r_n}E(Y^2_{n,j}|{\mathcal{F}}_{n,j-1})&\stackrel
{\mathit{i.p.}}\rightarrow&\sigma^2,\\
\mbox{\textup{(ii)}}\quad
\sum_{j=1}^{r_n}E\bigl(Y^2_{n,j}I_{(|Y_{n,j}|\geq\varepsilon
)}\bigr)&\rightarrow&
0,
\end{eqnarray*}
where $\sigma^2$ is a positive constant and $\varepsilon$ is an
arbitrary positive number, then
\[
\sum_{j=1}^{r_n}Y_{n,j}\stackrel{D}\rightarrow N(0,
\sigma^2).
\]
\end{lemma}
\begin{lemma}[(\cite{b4}, Lemma 2.7)]\label{lem1}
Let $\bbY=(Y_1,\ldots,Y_p)^T$, where $Y_i$'s are i.i.d. real r.v.'s
with mean $0$ and variance $1$. Let $\bbB=(b_{ij})_{p\times p}$, a
deterministic complex matrix. Then for any $k\geq2$, we have
\[
E|Y^T\bbB Y-\operatorname{tr}\bbB|^k\leq
\mathfrak{M}_k(EY_1^4\operatorname{tr}\bbB\bbB^*)^{k/2}+\mathfrak
{M}_kE(Y_1)^{2k}\operatorname{tr}(\bbB\bbB^*)^{k/2},
\]
where $\bbB^*$ denotes the complex conjugate transpose of $\bbB$.
\end{lemma}
\begin{lemma}\label{lem4}
Let $\bbC=(c_{ij})_{p\times p}$ be a deterministic complex matrix
with $c_{jj}=0$ and $\bbY=(Y_1,\ldots,Y_p)^T$, defined in Lemma
\ref{lem1}. Then for any $k\geq2$,
%
%e6.1 ###
%
\setcounter{equation}{0}
\begin{equation}
E|Y^T\bbC Y|^k\leq\mathfrak{M}_k(E|Y_1|^{k})^2(\operatorname{tr}\bbC\bbC^*)^{k/2}.
\end{equation}
\end{lemma}

Lemma \ref{lem4} directly follows from the argument of Lemma A.1
in \cite{b4}.
\begin{lemma}\label{lem2}
Under the assumptions of Theorem \ref{theo3}, as
$n\rightarrow\infty$,
%
%e6.2 ###
%
\begin{equation}
\label{a58}\max_{i}\sqrt{n}|E(\bbe_i^T\bbA_1^{-1}(z)\bar
{\mathbf s}_1)|\rightarrow0.
\end{equation}
\end{lemma}
\begin{pf}
We first prove that for $i\neq j$, $
\sup_{i,j}\sqrt{n}|E(\bbe_j^T\bbA_1^{-1}(z)\bbe
_i)|\rightarrow
0. $ To this end, write
\[
\bbA_1(z)+z\bbI=\frac{1}{n}\sum_{m=2}^n{\mathbf s}_m{\mathbf s}_m^T.
\]
Multiplying by $\bbA_1^{-1}(z)$ from the right on both sides of
the above equality gives
\[
\bbI+z\bbA_1^{-1}(z)=\frac{1}{n}\sum_{m=2}^n{\mathbf s}_m{\mathbf
s}_m^T\bbA_{m1}^{-1}(z)\beta_{m1}(z).
\]
Using
%
%e6.3 ###
%
\begin{equation}
\label{a11}
\beta_{m1}(z)=b_{12}(z)-\beta_{m1}(z)b_{12}(z)\xi_{m1}(z)
\end{equation}
we obtain
%
%e6.4 ###
%
\begin{eqnarray}\label{a12}
\bbI+z\bbA_1^{-1}(z)&=&\frac{b_{12}(z)}{n}\sum_{m=2}^n{\mathbf
s}_m{\mathbf s}_m^T\bbA_{m1}^{-1}(z)\nonumber\\[-8pt]\\[-8pt]
&&{}-\frac{b_{12}(z)}{n}\sum_{m=2
}^n{\mathbf s}_m{\mathbf
s}_m^T\bbA_{m1}^{-1}(z)\beta_{m1}(z)\xi_{m1}(z).\nonumber
\end{eqnarray}
It follows that for $i\neq j$
%
%e6.5 ###
%
\begin{eqnarray}\label{a57}
&&  z\sqrt{n}E(\bbe_j^T\bbA_1^{-1}(z)\bbe_i)\nonumber\\
&&\quad=\frac{b_{12}(z)}{\sqrt{n}}\Biggl(\sum_{m=2}^nE(\bbe
_j^T\bbA_{m1}^{-1}(z)\bbe_i)-
\sum_{m=2}^nE(\bbe_j^T{\mathbf s}_m{\mathbf s}_m^T\bbA
_{m1}^{-1}(z)\beta_{m1}(z)\xi_{m1}(z)\bbe_i)\Biggr)\hspace*{-32pt}\\
&&\quad=b_{12}(z)\sqrt{n}\bigl(E(\bbe_j^T\bbA_{21}^{-1}(z)\bbe_i)-
E(\bbe_j^T{\mathbf s}_2{\mathbf s}_2^T\bbA_{21}^{-1}(z)\beta_{21}(z)\xi
_{21}(z)\bbe_i)\bigr).\nonumber
\end{eqnarray}
As in
(\ref{a3}), by Lemma \ref{lem1} and (\ref{g28}),
%
%e6.6 ###
%
\begin{equation}\label{g49} E|\xi_{21}(z)|^k=O(\varepsilon
_n^{2k-4}n^{-1}),\qquad k\geq2.
\end{equation}
Here and in what follows (in this lemma) $O(\varepsilon
_n^{2k-4}n^{-1})$ and other
bounds are independent of $i$ and $j$.

We conclude from (\ref{g48}) that
\begin{eqnarray*}
&&b_{12}(z)\sqrt{n}E(\bbe_j^T\bbA_{21}^{-1}(z)\bbe_i)\\
&&\qquad=b_{12}(z)\sqrt{n}\biggl[E(\bbe_j^T\bbA_1^{-1}(z)\bbe_i)+
E\biggl(\bbe_j^T\bbA_{21}^{-1}(z)\frac{{\mathbf s}_2{\mathbf s}_2^T}{n}\bbA
_{21}^{-1}(z)\bbe_i\beta_{21}(z)\biggr)\biggr]\\
&&\qquad=b_{12}(z)\sqrt{n}E(\bbe_j^T\bbA_1^{-1}(z)\bbe_i)+O(n^{-1/2}).
\end{eqnarray*}
For the second term in (\ref{a57}), first, by a martingale method
similar to (\ref{a1}) and (\ref{g48}) we have, for $\bbe_l=\bbe_i$
or $\bbe_j$,
%
%e6.7 ###
%
\begin{eqnarray}\label{g29}
&&E|\bbe_l^T\bbA_{21}^{-1}(z_1)\bbe_j-E(\bbe_l^T\bbA
_{21}^{-1}(z_1)\bbe_j)|^2\nonumber\\
&&\qquad=E\Biggl|\sum_{m=3}^n(E_m-E_{m-1})\bigl[\bbe_l^T\bigl(\bbA
_{21}^{-1}(z_1)-\bbA_{m21}^{-1}(z_1)\bigr)\bbe_j\bigr]\Biggr|^2\\
&&\qquad\leq\frac{M}{n^2}\sum_{m=3}^nE|{\mathbf s}_m^T\bbA
_{m21}^{-1}(z_1)\bbe_j\bbe^T_l\bbA_{m21}^{-1}(z_1){\mathbf s}_m|^2=O(n^{-1}).\nonumber
\end{eqnarray}
This and (\ref{g28}) ensure that
\begin{eqnarray*}
&&\biggl|\frac{1}{n}E\bigl[\bbe_j^T\bbA_{21}^{-1}(z)\bbe_i\bigl(\operatorname{tr}\bbA
_{21}^{-1}(z)-E\operatorname{tr}\bbA_{21}^{-1}(z)\bigr)\bigr]\biggr|
\\[2pt]
&&\qquad=\biggl|\frac{1}{n}E\bigl[\bigl(\bbe_j^T\bbA_{21}^{-1}(z)\bbe_i-E\bbe_j^T\bbA
_{21}^{-1}(z)\bbe_i\bigr)\bigl(\operatorname{tr}\bbA_{21}^{-1}(z)-E\operatorname{tr}\bbA_{21}^{-1}(z)\bigr)\bigr]\biggr|
\\[2pt]
&&\qquad\leq\frac{\mathfrak{M}}{n}\bigl(E|\bbe_j^T\bbA_{21}^{-1}(z)\bbe
_i-E\bbe_j^T\bbA_{21}^{-1}(z)\bbe_i|^2E|\operatorname{tr}\bbA_{21}^{-1}(z)-E\operatorname{tr}\bbA
_{21}^{-1}(z)|^2\bigr)^{1/2}\\[2pt]
&&\qquad\leq\frac{\mathfrak{M}}{n}.
\end{eqnarray*}
Second, appealing to (\ref{i1}) gives
\begin{eqnarray*}
&&
E(\bbe_j^T{\mathbf s}_2{\mathbf s}_2^T\bbA_{21}^{-1}(z)\bbe_i\gamma
_{21}(z))\\[2pt]
&&\qquad=E\bigl(\bigl({\mathbf s}_2^T\bbA_{21}^{-1}(z)\bbe_i\bbe_j^T{\mathbf
s}_2-\bbe
_j^T\bbA_{21}^{-1}(z)\bbe_i\bigr)\gamma_{21}(z)\bigr)
\\[2pt]
&&\qquad=\frac{EX_{11}^4-3}{n}E(\bbe_j^T\bbA_{21}^{-1}(z)\bbe_i\bbe
_j^T\bbA_{21}^{-1}(z)\bbe_j)+\frac{2}{n}E(\bbe_j^T\bbA
_{21}^{-2}(z)\bbe_i).
\end{eqnarray*}
It follows that
\begin{eqnarray*}
&&\sqrt{n}E(\bbe_j^T{\mathbf s}_2{\mathbf s}_2^T\bbA_{21}^{-1}(z)\bbe
_i\xi_{21}(z))
\\[2pt]
&&\qquad=\sqrt{n}E(\bbe_j^T{\mathbf s}_2{\mathbf s}_2^T\bbA
_{21}^{-1}(z)\bbe
_i\gamma_{21}(z))\\[2pt]
&&\qquad\quad{}+
\sqrt{n}E\biggl[\bbe_j^T\bbA_{21}^{-1}(z)\bbe_i\frac{1}{n}\bigl(\operatorname{tr}\bbA
_{21}^{-1}(z)-E\operatorname{tr}\bbA_{21}^{-1}(z)\bigr)\biggr]
\\[2pt]
&&\qquad=O(n^{-1/2}).
\end{eqnarray*}
On the other hand, in view of (\ref{g48}) and (\ref{g49}) we
obtain
\[
\sqrt{n}E(\bbe_j^T{\mathbf s}_2{\mathbf s}_2^T\bbA_2^{-1}(z)\bbe_i\beta
_{21}(z)\xi_{21}^2(z))=O(\varepsilon_n).
\]
Therefore, by (\ref{a11}) we find
\begin{eqnarray*}
&&\sqrt{n}E(\bbe_j^T{\mathbf s}_2{\mathbf s}_2^T\bbA_{21}^{-1}(z)\beta
_{21}(z)\xi_{21}(z)\bbe_i)
\\[2pt]
&&\qquad=\sqrt{n}b_{12}(z)[E(\bbe_j^T{\mathbf s}_2{\mathbf s}_2^T\bbA
_{21}^{-1}(z)\bbe_i\xi_{21}(z))
-E(\bbe_j^T{\mathbf s}_2{\mathbf s}_2^T\bbA_2^{-1}(z)\bbe_i\beta
_{21}(z)\xi
_{21}^2(z))]
\\[2pt]
&&\qquad=O(\varepsilon_n).
\end{eqnarray*}

Therefore, combining the above argument with (\ref{g47}), we have
%
%e6.8 ###
%
\begin{equation}\label{g44}
\sup_{i\neq
j}\bigl|\sqrt{n}E(\bbe_j^T\bbA_1^{-1}(z)\bbe_i)\bigr|\rightarrow0.
\end{equation}

Next, applying (\ref{a11}) two times gives
\begin{eqnarray*}
&&E(\bbe_i^T\bbA_1^{-1}(z_1)\bar{\mathbf s}_1)\\[2pt]
&&\qquad=\frac{1}{n}\sum_{m=2
}^nE(\bbe_i^T\bbA_{m1}^{-1}(z_1){\mathbf s}_m\beta_{m1}(z_1))
\\
&&\qquad=\frac{b_{12}^2(z_1)(n-1)}{n}[-E(\bbe_i^T\bbA
_{21}^{-1}(z_1){\mathbf s}_2\xi_{21}(z_1))\\
&&\qquad\quad\hspace*{69pt}{}+E(\bbe_i^T\bbA_{21}^{-1}(z_1){\mathbf s}_2\beta_{21}(z_1)\xi
_{21}^2(z_1))].
\end{eqnarray*}
Obviously, we conclude from (\ref{g49}), (\ref{g48}) and
H\"{o}lder's inequality that
\[
\biggl|\frac{n-1}{n}E(\bbe_i^T\bbA_{21}^{-1}(z_1){\mathbf s}_2\beta
_{21}(z_1)\xi_{21}^2(z_1))\biggr|
=O(n^{-1/2}\varepsilon_n),
\]
while (\ref{g38}), (\ref{g10}) and (\ref{g44}) yield
\begin{eqnarray*}
&&\max_{i}\biggl|\frac{n-1}{n}E(\bbe_i^T\bbA_{21}^{-1}(z_1){\mathbf s}_2\xi
_{21}(z_1))\biggr|\\
&&\qquad=\max_{i}\Biggl|\frac{EX_{11}^3(n-1)}{n^2}\sum
_{j=1}^pE[\bbe_i^T\bbA_{21}^{-1}(z_1)\bbe_j(\bbA
_{21}^{-1}(z_1))_{jj}]\Biggr|\\
&&\qquad\leq\frac{|EX_{11}^3|}{n}\max_{i}\sum_{j\neq
i}^p|E[\bbe_i^T\bbA_1^{-1}(z_1)\bbe_j(\bbA
_1^{-1}(z_1))_{jj}]|+\frac{\mathfrak{M}}{n}\\
&&\qquad\leq\mathfrak{M}\biggl|EX_{11}^3E\frac{1}{n}\operatorname{tr}\bbA_1^{-1}(z_1)\biggr|\max
_{i\neq j}|E(\bbe_i^T\bbA_1^{-1}(z_1)\bbe_j)|+\frac
{\mathfrak{M}}{n}\\
&&\qquad=o(n^{-1/2}).
\end{eqnarray*}
Here we also use the estimate, via (\ref{g29}),
\[
E\bigl|\bigl(\bbe_i^T\bbA_1^{-1}(z_1)\bbe_j-E(\bbe_i^T\bbA_1^{-1}(z_1)\bbe
_j)\bigr)\bigl((\bbA^{-1}(z_1))_{jj}-E(\bbA^{-1}(z_1))_{jj}\bigr)\bigr|=O(n^{-1}).
\]
Thus, the proof of (\ref{a58}) is
complete.
\end{pf}

%s6.2 ###
\subsection{Truncation of the underlying random variables}
\label{trunrv}

To guarantee the results holding under the fourth moment, it is
necessary to truncate and centralize the underlying r.v.'s at an
appropriate rate. As in \cite{b2}, (1.8), one may select a
positive sequence
$\varepsilon_n$ so that
%
%e6.9 ###
%
\begin{equation}\label{slow}
\varepsilon_n\rightarrow0 \quad\mbox{and}\quad
\varepsilon_n^{-4}EX_{11}^4I\bigl(|X_{11}|\geq
\varepsilon_n\sqrt{n}\bigr)\rightarrow0.
\end{equation}
Set
$\widehat{X}_{ij}=X_{ij}I(|X_{ij}|\leq
\varepsilon_n\sqrt{n})-EX_{ij}I(|X_{ij}|\leq
\varepsilon_n\sqrt{n})$ and
$\tbbX_n=\bbX_n-\hbbX_n=(\tilde{X}_{ij})$
with $\hbbX_n=(\hat{X}_{ij})$. Let
$\sigma_n=\sqrt{E|\hat{X}_{11}|^2}$,
$\cbbS_{n}=(n\sigma_n^2)^{-1}
\hbbX_n\hbbX_n^T$ and $\cbbA^{-1}(z)=(\cbbS_n-zI)^{-1}$. Moreover, introduce
$\bar{\check{{\mathbf s}}}=\frac{1}{n}\sum_{j=1}^n\check{{\mathbf
s}}_j$, where $\check{{\mathbf s}}_j$ is the $j$th
column of the matrix $(\sigma_n)^{-1}\hbbX_n$.
\begin{lemma}\label{trunlem}
Assume that $X_{ij}, i=1,\ldots,p, j=1,\ldots,n$ are i.i.d. with
$EX_{11}=0, E|X_{11}|^2=1$ and $E|X_{11}|^4<\infty$, for $z\in
\mathcal{C}_n^+$, we have then
%
%e6.10 ###
%
\begin{equation}\label{applem1}
\sqrt{n}\bigl(\bar{{\mathbf s}}^T\bbA^{-1}(z)\bar{\mathbf s}-\bar{\check
{{\mathbf s}}}{}^T\cbbA^{-1}(z)\bar{\check{{\mathbf s}}}\bigr)\stackrel
{\mathit{i.p.}}\longrightarrow
0,
\end{equation}
where the convergence in probability holds
uniformly for $z\in\mathcal{C}_n^+$. Moreover,
%
%e6.11 ###
%
\begin{equation}
\label{slow1}
\sqrt{n}(\bar{{\mathbf s}}^T\bar{\mathbf s}-\bar{\check{{\mathbf
s}}}{}^T\bar
{\check{{\mathbf s}}})\stackrel{\mathit{i.p.}}\longrightarrow
0.
\end{equation}
\end{lemma}
\begin{pf}
Write
\[
\sqrt{n}\bigl(\bar{{\mathbf s}}^T\bbA^{-1}(z)\bar{\mathbf s}-\bar{\check
{{\mathbf s}}}{}^T\cbbA^{-1}(z)\bar{\check{{\mathbf
s}}}\bigr)=u_{n1}+u_{n2}+u_{n3},
\]
where
\[
u_{n1}=\sqrt{n}[(\bar{{\mathbf s}}-\bar{\check{{\mathbf s}}})^T\bbA
^{-1}(z)\bar{\mathbf s}],
u_{n2}=\sqrt{n}\bigl[\bar{\check{{\mathbf s}}}{}^T\bigl(\bbA^{-1}(z)-\cbbA
^{-1}(z)\bigr)\bar{\mathbf s}\bigr]
\]
and
\[
u_{n3}=\sqrt{n}[\bar{\check{{\mathbf s}}}{}^T\cbbA^{-1}(z)(\bar{\mathbf
s}-\bar{\check{{\mathbf s}}})].
\]

Consider $u_{n1}$ on the $\mathcal{C}_u$ first. It is observed
that
%
%e6.12 ###
%
\begin{eqnarray}\label{g35}
|u_{n1}|&\leq&
\sqrt{n}\|(\bar{{\mathbf s}}-\bar{\check{{\mathbf s}}})^T\|\|\bbA
^{-1}(z)\|
\|\bar{\mathbf s}\|\leq\frac{\sqrt{n}}{v_0}\|(\bar{{\mathbf s}}-\bar
{\check{{\mathbf s}}})^T\|\|\bar{\mathbf s}\|\nonumber\\[-8pt]\\[-8pt]
&\leq&\frac{\sqrt{n}}{v_0}\biggl|1-\frac{1}{\sigma_n}\biggr|\|\bar{{\mathbf s}}\|
^2+\frac{\sqrt{n}}{v_0}\frac{1}{\sigma_n}\|\bar{\tilde{{\mathbf s}}}\|
\|\bar{\mathbf s}\|,\nonumber
\end{eqnarray}
since $
\bar{{\mathbf s}}-\bar{\check{{\mathbf s}}}=(1-\frac{1}{\sigma_n})\bar
{{\mathbf s}}+\frac{1}{\sigma_n}\bar{\tilde{{\mathbf s}}}$
with $\bar{\tilde{{\mathbf s}}}=\sum_{j=1}^n\tilde{{\mathbf s}}_j/n$ and
$\tilde{{\mathbf s}}_j$ being the $j$th column of~$\tbbX_n$. Moreover, it
follows from (\ref{slow}) that
\[
1-\sigma_n^2\leq2EX_{11}^2I\bigl(|X_{11}|\geq
\varepsilon_n\sqrt{n}\bigr)\leq
2\varepsilon_n^{-2}n^{-1}EX_{11}^4I\bigl(|X_{11}|\geq
\varepsilon_n\sqrt{n}\bigr)=o(\varepsilon_n^{2}n^{-1}),
\]
which implies that
%
%e6.13 ###
%
\begin{equation}\label{g36}
\sqrt{n}(1-1/\sigma_n)=\sqrt{n}(\sigma
_n^2-1)/[\sigma_n(1+\sigma_n)]=o(n^{-1/2}).
\end{equation}
On the other hand,
\[
E\|\bar{\tilde{{\mathbf s}}}\|^2=E\Biggl[\sum_{i=1}^p\Biggl|\frac{1}{n}\sum
_{j=1}^n\tilde{X}_{ij}\Biggr|^2\Biggr]=\frac{1}{n^2}\sum_{i=1}^p
\sum_{j=1}^nE\tilde{X}_{ij}^2\leq
\frac{\mathfrak{M}}{n\varepsilon_n^2}EX_{11}^4I\bigl(|X_{11}|\geq
\varepsilon_n\sqrt{n}\bigr),
\]
which, via (\ref{slow}), gives that
%
%e6.14 ###
%
\begin{equation}\label{g37}
\sqrt{n}\|\bar{\tilde{{\mathbf s}}}\|\stackrel{\mathrm{i.p.}}\longrightarrow0.
\end{equation}
In addition, $\|\bar{\mathbf s}\|^2$ is uniformly integrable because
(\ref{7}) remains true for $k=2$ without truncation by a careful
check on its argument. This, together with (\ref{g35})--(\ref{g37}),
ensures that $u_{n1}$ converges in probability to zero uniformly on
$\mathcal{C}_u$.\vspace*{1pt}

Analyze $u_{n2}$ next. Since
$\bbX_n-\sigma_n^{-1}\hbbX_n=(1-\sigma_n^{-1})\bbX_n+\sigma
_n^{-1}\tbbX_n$,
we have
\begin{eqnarray*}
&&|u_{n2}|\leq\sqrt{n}\|\bar{\check{{\mathbf s}}}{}^T\|\|\bbA
^{-1}(z)-\cbbA^{-1}(z)\|\|\bar{\mathbf s}\|\leq
\frac{\sqrt{n}}{v_0^2}\|\bar{\check{{\mathbf s}}}{}^T\|\|\bbA(z)-\cbbA
(z)\|\|\bar{\mathbf s}\|
\\
&&\qquad\leq\frac{1}{v_0^2\sqrt{n}}\|\bar{\check{{\mathbf s}}}{}^T\|\|\bar
{\mathbf s}\|[\|\bbX_n-\sigma_n^{-1}\hbbX_n\|\|\bbX_n^T\|
+\|\sigma_n^{-1}\hbbX_n\|\|\bbX_n^T-\sigma_n^{-1}\hbbX_n^T\|]
\\
&&\qquad\leq\frac{1}{v_0^2\sqrt{n}}\|\bar{\check{{\mathbf s}}}{}^T\|\|\bar
{\mathbf s}\|[(1-\sigma_n^{-1})\|\bbX_n\|\|\bbX_n^T\|+\sigma_n^{-1}\|
\tbbX_n\|\|\bbX_n^T\|
\\
&&\hspace*{66pt}\qquad\quad{}+\|\sigma_n^{-1}\hbbX_n\|(1-\sigma_n^{-1})\|\bbX_n^T\|+\|
\sigma
_n^{-1}\hbbX_n\|\sigma_n^{-1}\|\tbbX_n^T\|].
\end{eqnarray*}

As before, $\|\bar{\check{{\mathbf s}}}\|$ and $\|\bar{\mathbf s}\|$
are uniformly
integrable. Moreover, the spectral norms $\|\bbX_n^T\|/\sqrt{n}$ and
$\|\sigma_n^{-1}\hbbX_n\|/\sqrt{n}$ both converge to
$(1+\sqrt{c})^2$ with probability one by~\cite{y2}. In addition,
$\|\tbbX_n^T\|/\sqrt{nE\tilde{X}_{11}^2}$ converges to
$(1+\sqrt{c})^2$ with probability one. From (\ref{slow}) we have
\[
nE\tilde{X}_{11}^2\leq2\varepsilon_n^{-2}EX_{11}^4I\bigl(|X_{11}|\geq
\varepsilon_n\sqrt{n}\bigr)=O(\varepsilon_n^2),
\]
which, together with
(\ref{g36}), yields that $u_{n2}$ converges in probability to zero
uniformly on $\mathcal{C}_u$.

Clearly, the argument for $u_{n1}$ works for $u_{n3}$ as well.
Moreover, note that $\|\bbA^{-1}(z)\|$ is bounded for
$z\in\mathcal{C}_l,u_l<0$. As for $z\in\mathcal{C}_l,u_l>0$ or
$z\in\mathcal{C}_r$, by \cite{y2} we have
\[
\lim_{n\rightarrow\infty}\min\bigl(u_r-\lambda_{\max}(\bbA
),\lambda_{\min}(\bbA)-u_l\bigr)>0,\qquad
\mbox{a.s.}
\]
and
\[
\lim_{n\rightarrow\infty}\min\bigl(u_r-\lambda_{\max}(\cbbA
),\lambda_{\min}(\cbbA)-u_l\bigr)>0,\qquad
\mbox{a.s.}
\]
Therefore, the above argument for $u_{nj},j=1,2,3$ for
$z\in\mathcal{C}_u$ of course applies to the cases
(1) $z\in\mathcal{C}_l,u_l<0$; (2) $z\in\mathcal{C}_l,u_l>0$; (3)
$z\in\mathcal{C}_r$. Thus, (\ref{applem1}) holds.

Finally, the above argument for (\ref{applem1}) certainly works
for (\ref{slow1}). Thus, the proof is complete.
\end{pf}
\end{appendix}

\section*{Acknowledgments}
The authors would like to thank the editor, an associate editor and a
referee for their constructive comments which have helped to improve
the paper a great deal.

% imsref loaded by lrinkeviciute, 2011-01-13 11:04:01
%
% imsref loaded by lrinkeviciute, 2011-01-14 16:20:40

%
\printaddresses


\begin{thebibliography}{27}
% BibTex style file: ims-number.bst, 2010-01-07

%b1 ###
\bibitem{a1}
\begin{bbook}[mr]
\bauthor{\bsnm{Anderson},~\bfnm{T.~W.}\binits{T.~W.}}
(\byear{1984}).
\btitle{An Introduction to Multivariate Statistical Analysis},
\bedition{2nd} ed.
%  and Mathematical Statistics}.
\bpublisher{Wiley}, \baddress{New York}.
\bid{mr={0771294}}
\end{bbook}
\endbibitem

%b2 ###
\bibitem{b5}
\begin{barticle}[mr]
\bauthor{\bsnm{Bai},~\bfnm{Zhidong}\binits{Z.}} \AND
  \bauthor{\bsnm{Saranadasa},~\bfnm{Hewa}\binits{H.}}
(\byear{1996}).
\btitle{Effect of high dimension: By an example of a two sample problem}.
\bjournal{Statist. Sinica}
\bvolume{6}
\bpages{311--329}.
\bid{mr={1399305}}
\end{barticle}
\endbibitem

%b3 ###
\bibitem{b1}
\begin{barticle}[mr]
\bauthor{\bsnm{Bai},~\bfnm{Z.~D.}\binits{Z.~D.}},
  \bauthor{\bsnm{Miao},~\bfnm{B.~Q.}\binits{B.~Q.}} \AND
  \bauthor{\bsnm{Pan},~\bfnm{G.~M.}\binits{G.~M.}}
(\byear{2007}).
\btitle{On asymptotics of eigenvectors of large sample covariance matrix}.
\bjournal{Ann. Probab.}
\bvolume{35}
\bpages{1532--1572}.
\bid{doi={10.1214/009117906000001079}, mr={2330979}}
\end{barticle}
\endbibitem

%b4 ###
\bibitem{b4}
\begin{barticle}[mr]
\bauthor{\bsnm{Bai},~\bfnm{Z.~D.}\binits{Z.~D.}} \AND
  \bauthor{\bsnm{Silverstein},~\bfnm{Jack~W.}\binits{J.~W.}}
(\byear{1998}).
\btitle{No eigenvalues outside the support of the limiting spectral
  distribution of large-dimensional sample covariance matrices}.
\bjournal{Ann. Probab.}
\bvolume{26}
\bpages{316--345}.
\bid{doi={10.1214/aop/1022855421}, mr={1617051}}
\end{barticle}
\endbibitem

%b5 ###
\bibitem{b2}
\begin{barticle}[mr]
\bauthor{\bsnm{Bai},~\bfnm{Z.~D.}\binits{Z.~D.}} \AND
  \bauthor{\bsnm{Silverstein},~\bfnm{Jack~W.}\binits{J.~W.}}
(\byear{2004}).
\btitle{C{LT} for linear spectral statistics of large-dimensional sample
  covariance matrices}.
\bjournal{Ann. Probab.}
\bvolume{32}
\bpages{553--605}.
\bid{doi={10.1214/aop/1078415845}, mr={2040792}}
\end{barticle}
\endbibitem

%b6 ###
\bibitem{bili2}
\begin{bbook}[mr]
\bauthor{\bsnm{Billingsley},~\bfnm{Patrick}\binits{P.}}
(\byear{1968}).
\btitle{Convergence of Probability Measures}.
\bpublisher{Wiley}, \baddress{New York}.
\bid{mr={0233396}}
\end{bbook}
\endbibitem

%b7 ###
\bibitem{bili}
\begin{bbook}[mr]
\bauthor{\bsnm{Billingsley},~\bfnm{Patrick}\binits{P.}}
(\byear{1995}).
\btitle{Probability and Measure},
\bedition{3rd} ed.
\bpublisher{Wiley}, \baddress{New York}.
\bid{mr={1324786}}
\end{bbook}
\endbibitem

%b8 ###
\bibitem{burk}
\begin{barticle}[mr]
\bauthor{\bsnm{Burkholder},~\bfnm{D.~L.}\binits{D.~L.}}
(\byear{1973}).
\btitle{Distribution function inequalities for martingales}.
\bjournal{Ann. Probab.}
\bvolume{1}
\bpages{19--42}.
\bid{mr={0365692}}
\end{barticle}
\endbibitem

%b9 ###
\bibitem{Di}
\begin{barticle}[mr]
\bauthor{\bsnm{Diaconis},~\bfnm{Persi}\binits{P.}} \AND
  \bauthor{\bsnm{Evans},~\bfnm{Steven~N.}\binits{S.~N.}}
(\byear{2001}).
\btitle{Linear functionals of eigenvalues of random matrices}.
\bjournal{Trans. Amer. Math. Soc.}
\bvolume{353}
\bpages{2615--2633}.
\bid{doi={10.1090/S0002-9947-01-02800-8}, mr={1828463}}
\end{barticle}
\endbibitem

%b10 ###
\bibitem{h}
\begin{barticle}[auto:STB|2010-11-18|09:18:59]
\bauthor{\bsnm{Hotelling},~\bfnm{H.}\binits{H.}}
(\byear{1931}).
\btitle{The generalization of Student's ratio}.
\bjournal{Ann. Math. Statist.}
\bvolume{2}
\bpages{360--378}.
\end{barticle}
\endbibitem

%b11 ###
\bibitem{j}
\begin{barticle}[mr]
\bauthor{\bsnm{Jiang},~\bfnm{Tiefeng}\binits{T.}}
(\byear{2004}).
\btitle{The limiting distributions of eigenvalues of sample correlation
  matrices}.
\bjournal{Sankhy\=a}
\bvolume{66}
\bpages{35--48}.
\bid{mr={2082906}}
\end{barticle}
\endbibitem

%b12 ###
\bibitem{johan}
\begin{barticle}[mr]
\bauthor{\bsnm{Johansson},~\bfnm{Kurt}\binits{K.}}
(\byear{1998}).
\btitle{On fluctuations of eigenvalues of random {H}ermitian matrices}.
\bjournal{Duke Math. J.}
\bvolume{91}
\bpages{151--204}.
\bid{doi={10.1215/S0012-7094-98-09108-6}, mr={1487983}}
\end{barticle}
\endbibitem

%b13 ###
\bibitem{john}
\begin{barticle}[mr]
\bauthor{\bsnm{Johnstone},~\bfnm{Iain~M.}\binits{I.~M.}}
(\byear{2001}).
\btitle{On the distribution of the largest eigenvalue in principal components
  analysis}.
\bjournal{Ann. Statist.}
\bvolume{29}
\bpages{295--327}.
\bid{doi={10.1214/aos/1009210544}, mr={1863961}}
\end{barticle}
\endbibitem

%b14 ###
\bibitem{Jon}
\begin{barticle}[mr]
\bauthor{\bsnm{Jonsson},~\bfnm{Dag}\binits{D.}}
(\byear{1982}).
\btitle{Some limit theorems for the eigenvalues of a sample covariance matrix}.
\bjournal{J. Multivariate Anal.}
\bvolume{12}
\bpages{1--38}.
\bid{doi={10.1016/0047-259X(82)90080-X}, mr={0650926}}
\end{barticle}
\endbibitem

%b15 ###
\bibitem{leh}
\begin{bbook}[mr]
\bauthor{\bsnm{Lehmann},~\bfnm{E.~L.}\binits{E.~L.}} \AND
  \bauthor{\bsnm{Romano},~\bfnm{Joseph~P.}\binits{J.~P.}}
(\byear{2005}).
\btitle{Testing Statistical Hypotheses},
\bedition{3rd} ed.
\bpublisher{Springer}, \baddress{New York}.
\bid{mr={2135927}}
\end{bbook}
\endbibitem

%b16 ###
\bibitem{MP}
\begin{barticle}[auto:STB|2010-11-18|09:18:59]
\bauthor{\bsnm{Mar{\v{c}}enko},~\bfnm{V.~A.}\binits{V.~A.}} \AND
  \bauthor{\bsnm{Pastur},~\bfnm{L.~A.}\binits{L.~A.}}
(\byear{1967}).
\btitle{Distribution for some sets of random matrices}.
\bjournal{Math. USSR-Sb.}
\bvolume{1}
\bpages{457--483}.
\end{barticle}
\endbibitem

%b17 ###
\bibitem{p2}
\begin{barticle}[mr]
\bauthor{\bsnm{Pan},~\bfnm{G.~M.}\binits{G.~M.}} \AND
  \bauthor{\bsnm{Zhou},~\bfnm{W.}\binits{W.}}
(\byear{2008}).
\btitle{Central limit theorem for signal-to-interference ratio of reduced rank
  linear receiver}.
\bjournal{Ann. Appl. Probab.}
\bvolume{18}
\bpages{1232--1270}.
\bid{doi={10.1214/07-AAP477}, mr={2418244}}
\end{barticle}
\endbibitem

%%b18 ###
%(\byear{1984}).
%  covariance matrices}.

%b19 ###
\bibitem{s1}
\begin{barticle}[mr]
\bauthor{\bsnm{Silverstein},~\bfnm{Jack~W.}\binits{J.~W.}}
(\byear{1989}).
\btitle{On the eigenvectors of large-dimensional sample covariance matrices}.
\bjournal{J. Multivariate Anal.}
\bvolume{30}
\bpages{1--16}.
\bid{doi={10.1016/0047-259X(89)90084-5}, mr={1003705}}
\end{barticle}
\endbibitem

%b20 ###
\bibitem{s2}
\begin{barticle}[mr]
\bauthor{\bsnm{Silverstein},~\bfnm{Jack~W.}\binits{J.~W.}}
(\byear{1990}).
\btitle{Weak convergence of random functions defined by the eigenvectors of
  sample covariance matrices}.
\bjournal{Ann. Probab.}
\bvolume{18}
\bpages{1174--1194}.
\bid{mr={1062064}}
\end{barticle}
\endbibitem

%b21 ###
\bibitem{s3}
\begin{barticle}[mr]
\bauthor{\bsnm{Silverstein},~\bfnm{Jack~W.}\binits{J.~W.}}
(\byear{1995}).
\btitle{Strong convergence of the empirical distribution of eigenvalues of
  large-dimensional random matrices}.
\bjournal{J. Multivariate Anal.}
\bvolume{55}
\bpages{331--339}.
\bid{doi={10.1006/jmva.1995.1083}, mr={1370408}}
\end{barticle}
\endbibitem

%b22 ###
\bibitem{ss}
\begin{barticle}[mr]
\bauthor{\bsnm{Sinai},~\bfnm{Ya.}\binits{Y.}} \AND
  \bauthor{\bsnm{Soshnikov},~\bfnm{A.}\binits{A.}}
(\byear{1998}).
\btitle{Central limit theorem for traces of large random symmetric matrices
  with independent matrix elements}.
\bjournal{Bol. Soc. Brasil. Mat. (N.S.)}
\bvolume{29}
\bpages{1--24}.
\bid{doi={10.1007/BF01245866}, mr={1620151}}
\end{barticle}
\endbibitem

%b23 ###
\bibitem{watch}
\begin{barticle}[mr]
\bauthor{\bsnm{Wachter},~\bfnm{Kenneth~W.}\binits{K.~W.}}
(\byear{1978}).
\btitle{The strong limits of random matrix spectra for sample matrices of
  independent elements}.
\bjournal{Ann. Probab.}
\bvolume{6}
\bpages{1--18}.
\bid{mr={0467894}}
\end{barticle}
\endbibitem

%b24 ###
\bibitem{z}
\begin{barticle}[mr]
\bauthor{\bsnm{Xiao},~\bfnm{Han}\binits{H.}} \AND
  \bauthor{\bsnm{Zhou},~\bfnm{Wang}\binits{W.}}
(\byear{2010}).
\btitle{Almost sure limit of the smallest eigenvalue of some sample correlation
  matrices}.
\bjournal{J. Theoret. Probab.}
\bvolume{23}
\bpages{1--20}.
\bid{doi={10.1007/s10959-009-0270-2}, mr={2591901}}
\end{barticle}
\endbibitem

%b25 ###
\bibitem{y1}
\begin{barticle}[mr]
\bauthor{\bsnm{Yin},~\bfnm{Y.~Q.}\binits{Y.~Q.}}
(\byear{1986}).
\btitle{Limiting spectral distribution for a class of random matrices}.
\bjournal{J. Multivariate Anal.}
\bvolume{20}
\bpages{50--68}.
\bid{doi={10.1016/0047-259X(86)90019-9}, mr={0862241}}
\end{barticle}
\endbibitem

%b26 ###
\bibitem{y2}
\begin{barticle}[mr]
\bauthor{\bsnm{Yin},~\bfnm{Y.~Q.}\binits{Y.~Q.}},
  \bauthor{\bsnm{Bai},~\bfnm{Z.~D.}\binits{Z.~D.}} \AND
  \bauthor{\bsnm{Krishnaiah},~\bfnm{P.~R.}\binits{P.~R.}}
(\byear{1988}).
\btitle{On the limit of the largest eigenvalue of the large-dimensional sample
  covariance matrix}.
\bjournal{Probab. Theory Related Fields}
\bvolume{78}
\bpages{509--521}.
\bid{doi={10.1007/BF00353874}, mr={0950344}}
\end{barticle}
\endbibitem

\end{thebibliography}
\end{document}